\documentclass[12pt,oneside,draft]{article}
\usepackage{mathtext}
\usepackage[T2A]{fontenc}
\usepackage{latexsym,amsfonts,amssymb,amsmath,longtable,amsthm}
\usepackage{cite}
\usepackage[pic]{xy}
\usepackage{makeidx,bbm}
\usepackage[unicode,final,hyperindex]{hyperref}
\usepackage[centerlast,small]{caption}
\renewcommand{\le}{\leqslant}
\renewcommand{\ge}{\geqslant}

\begin{document}
\sloppy

\newcounter{thelem}[section]
\newtheorem{lem}{{\scshape Lemma}}[subsection]
\newtheorem{ttt}[lem]{{\scshape Theorem}}
\newtheorem{cor}[lem]{{\scshape Corollary}}
\newtheorem{prp}[lem]{{\scshape Proposition}}
\theoremstyle{remark}
\newtheorem*{rem}{{\slshape Remark}}
\newtheorem{quest}[lem]{{\slshape Problem}}
\theoremstyle{definition}
\newtheorem{df}[lem]{{\bfseries Definition}}
\newtheorem{ex}[lem]{{\slshape Example}}
\renewcommand{\thelem}{\arabic{section}.\arabic{subsection}.\arabic{lem}}


\newcommand{\Sym}{{\mathrm{Sym}}} 
\newcommand{\Alt}{{\mathrm{Alt}}}
\newcommand{\rank}{{\mathrm{rank}}}
\newcommand{\Aut}{{\mathrm{Aut}}}
\newcommand{\Out}{{\mathrm{Out}}}
\newcommand{\Inn}{{\mathrm{Inn}}}
\newcommand{\N}{{\mathbbmss{N}}}    
\newcommand{\Z}{\mathbbmss{Z}}    
\newcommand{\Q}{\mathbbmss{Q}}    
\renewcommand{\C}{\mathbbmss{C}}    
\newcommand{\F}{\mathbbmss{F}}    
\newcommand{\R}{\mathbbmss{R}}    
\newcommand{\ch}{\mathrm{char}} 
\newcommand{\GL}{{\mathrm{GL}}}             
\newcommand{\GU}{{\mathrm{GU}}}     
\newcommand{\GO}{{\mathrm{GO}}}   
\newcommand{\SO}{{\mathrm{SO}}}   
\newcommand{\SL}{{\mathrm{SL}}}    
\newcommand{\PSL}{{\mathrm{PSL}}} 
\newcommand{\PSp}{{\mathrm{PSp}}} 
\newcommand{\PGL}{{\mathrm{PGL}}} 
\newcommand{\SU}{{\mathrm{SU}}}  
\newcommand{\PSU}{{\mathrm{PSU}}}  
\newcommand{\Sp}{{\mathrm{Sp}}}     
\newcommand{\la}{\mathop{\langle}}         
\newcommand{\ra}{\mathop{\rangle}}            
\newcommand{\diag}{{\mathrm{diag}}}       
\newcommand{\block}{{\mathrm{block-diag}}}    
\newcommand{\ov}{\overline}
\renewcommand{\thesubsection}{\arabic{section}.\arabic{subsection}}
\renewcommand{\P}{\mathbf{P}}


\hfill{\normalsize MSC2000 20D25, 20D06, 20E07, 20G40}

\begin{center}
{\bfseries CARTER SUBGROUPS OF FINITE GROUPS}

E.P.Vdovin
\end{center}

{\tiny\tableofcontents}

{\scriptsize ABSTRACT: In the paper it is proven that Carter subgroups of a
finite group are conjugate. A complete classification of Carter subgroups in
finite almost simple groups is also obtained.}
\bigskip


\section{Introduction}

The present paper is a slightly shorten version of doctoral thesis
``Carter subgroups of finite groups''. The results of the thesis were published
in \cite{VdoReductTheorem,VdoCartAlmSimple,VdoCarExists,VdoDAN},
\cite{PreTamVdoCartClassical}, and \cite{TamVdoCartSbgrpsConj}. Although,
in the thesis all papers were rewritten in order to obtain a coherent text.
This rewritten text is given in the present paper.

\subsection{General characteristic of the results of the paper}

We recall that a subgroup of a finite group is called a {\em Carter
subgroup} if it is nilpotent and self-normalizing. By a well-known
result, any finite solvable group contains exactly one conjugacy
class of Carter subgroups (see \cite{CarNilpSelfNormSubgrpsSolvGrps}). If a
group is not assumed to be finite, then Carter subgroups can be even
nonisomorhic. Indeed, if $N_1,N_2$ are two nonisomorhic nilpotent groups, then
they are Carter subgroups in their free product. On the other hand, a finite
nonsolvable group may fail to contain Carter subgroups, the minimal counter
example is the alternating group of degree $5$. Although there is not known a
finite group containing nonconjugate Carter subgroups, and the following
problem, due to R.Carter, is known.

\begin{quest} ({\bfseries Conjugacy Problem})
Are Carter subgroups of a finite group conjugate?
\end{quest}

This problem for several classes of finite groups closed to be simple
was investigated by many authors. In symmetric and alternating groups Carter
subgroups were classified by L.Di Martino and M.C.Tamburini (see
\cite{DT1CartSbgrpsSymAlt}). In every group $G$ such that
$\SL_n(q)\leq G\leq \GL_n(q)$ Carter subgroups were classified by L.Di
Martino and M.C.Tamburini, and in the case $G=\GL_n(q)$, by N.A.Vavilov (see
\cite{DT2CartSbgrpsPGL} and \cite{Va} respectively). For symplectic groups
$\Sp_{2n}(q)$, general unitary groups $\GU_n(q)$, and, for $q$ odd, general
orthogonal groups $\GO_n^{\pm}(q)$ the classification of Carter subgroups was
obtained by L.Di Martino, A.E.Zalessky, and M.C.Tamburini (see
\cite{DTZCartSbgrpsClassGrps}). For some sporadic simple groups Carter
subgroups were found in~\cite{D'ACartSbgrpsSomeSporadic}. In the mentioned
above nonsolvable groups Carter subgroups coinside with the normalizers of
Sylow $2$-sub\-gro\-ups, and hence are conjugate.

A finite group $G$ is called a {\em minimal counter example to  Conjugacy
Problem} or simply a {\em minimal counter example}, if $G$ contains
nonconjugate Carter subgroups, but in every group  $H$, with $\vert
H\vert< \vert G\vert$, Carter subgroups are conjugate. In
\cite{DLTReductionAlmSimpleGrps} F.Dalla Volta, A.Lucchini, and M.C.Tamburini
have proven that a minimal counter example should be almost simple. This result
allows to use the classification of finite simple groups to solve Conjugacy
Problem.

Note that the using of F.Dalla Volta, A.Lucchini, and M.C.Tamburini
result to the classification of Carter subgroups in almost simple groups
essentially depends on the classification of finite simple groups. Indeed, in
order to use the inductive hypothesis that Carter subgroups in every proper
subgroup of a minimal counter example are conjugate, one needs to know that all
almost simple groups of order less than the order of a minimal counter example
are found. To avoid using the classification of finite simple groups we
strengthen the result of F.Dalla Volta, A.Lucchini, and M.C.Tamburini, proving
that if Carter subgroups are conjugate in the group of induced automorphisms of
every non-Abelian composition factor, then they are conjugate in the group.

To inductive description of Carter subgroups in almost simple groups one needs
to know homomorphic images of Carter subgroups and intersections of
Carter subgroups with normal subgroups, i.~e., the answers to the following
problems.

\begin{quest}
Is a homomorphic image of a Carter subgroup again a Carter subgroup?
\end{quest}

\begin{quest}
Is the intersection of a Carter subgroup with a normal subgroup again a
Carter subgroup (of the normal subgroup)?
\end{quest}

The first problem is closely connected with Conjugacy Problem, namely, if
Conjugacy Problem has an affirmative answer, then the first problem also
has an affirmative answer. So we shall solve both of these problems by
considering Carter subgroups in almost simple groups. It is easy to see that
the second problem has a negative answer. Indeed, consider a solvable group
$\Sym_3$ and its normal subgroup of index $2$, the alternating group $\Alt_3$.
Then a Carter subgroup of $\Sym_3$ is a Sylow $2$-sub\-gro\-up, while a Carter
subgroup of $\Alt_3$ is a Sylow $3$-sub\-gro\-up. Thus in the paper some
properties
of Carter subgroups in a group and some of its normal subgroups are found.

The present paper is divided into six sections, including Introduction.
In the introduction we give general results of the paper, and also necessary
definitions and results.

In the second section we prove that Carter subgroups of a finite group are
conjugate, if they are conjugate in the group of induced automorphisms
of every its non-Abelian composition factor, thereby strengthen the results of
F.Dalla Volta, A.Lucchini, and M.C.Tamburini. In the second section we also
obtain some properties of Carter subgroups.

In the third section we consider the problem of conjugacy for elements of prime
order in finite groups of Lie type. At the end of the third section, by
using the results on conjugacy, we obtain the classification of Carter subgroups
in a wide class of almost simple groups.

In the fourth section we introduce the notion of semilinear groups of Lie type
and corresponding semilinear algebraic groups, and transfer the results on the
normalizers of $p$-sub\-gro\-ups and the centralizers of semisimple elements in
groups of Lie type. We also obtain some additional results on the conjugacy of
elements of prime order in these groups.

In the fifth section we complete the classification of Carter subgroups in
almost simple groups and prove that Carter subgroups of almost simple groups
are conjugate. As a corollary we obtain an affirmative answer to Conjugacy
Problem and prove that a homomorphic image of a Carter subgroup is a Carter
subgroup.

In the sixth section we study the problem of existence of a Carter subgroup in a
finite group, give a criterion of existence, and construct an example showing
that the property of containing a Carter subgroup is not preserved under
extensions. More over in the last subsection of sixth section we give tables
with classification of Carter subgroups in almost simple groups.

I am grateful to my scientific adviser, a corresponding member of RAS,
V.D.Mazurov.
В.~Д.~Мазурову. His contribution to my development as a mathematician, an his
continuous support are inestimable. I am also sincerely thankful to Professor
M.C.Tamburini, who has initiated my work in these problems, and has provided a
help during the work. I especially thank Dr. A.V.Vasiliev, Dr. M.A.Grechkoseeva,
Dr. A.V.Zavarnitsine, and Dr. D.O.Revin for very useful discussion of the
paper, that allows to simplify some proofs and to improve inaccuracies and
mistakes. I am also pleased to Professor A.S.Kondratiev for valuable comments,
that improving the final text. I am grateful and wish to honor a serene memory
of Professor Yu.I.Merzlyakov, who awaken my interest to algebra and group
theory.

The work is supported by Russian fond of basic research (projects numbers
99--01--00550, 01--01--06184, 02--01--00495, 02--01--06226 and 05--01--00797),
grants of President RF for young scientists (MK--1455.2005.1 and
MK--3036.2007.1), SB RAS (gran N~29 for young scientists and Integration
project 2006.1.2), and Program ``Universities of Russia'' (project number
UR.04.01.202). A part of the work has been made during my post-doctoral
fellowship in the university of Padua (Italy), and I am grateful to this
university, to all members of algebra chair, and, especially to Professor
F.Menegazzo for support.

\subsection{Notation and results from group theory}

Out notation is standard. If $G$ is a group, then $H\leq
G$\glossary{HleqG@$H\leq G$} and
$H\unlhd G$\glossary{HunlhdG@$H\unlhd G$} mean that $H$ is a subgroup and a
normal subgroup of $G$ respectively. By
$\left|G:H\right|$\glossary{G:H@$\vert G:H\vert$} we denote the index of
$H$ in $G$,
$N_G(H)$\glossary{NG(H)@$N_G(H)$} is the normalizer of $H$ in $G$. If
$H$ is normal in $G$, then by
$G/H$\glossary{G/H@$G/H$} we denote the factor group of $G$ by
$H$. If $M$ is a subset of $G$, then $\langle
M\rangle$\glossary{Mrangle@$\la M\ra$} denotes the subgroup generated by $M$,
$\vert M\vert$ denotes \glossary{Mvert@$\vert M\vert$}  the cardinality
of $M$ (or the order of an element, if there is an element instead of a set).
By $C_G(M)$\glossary{CG(M)@$C_G(M)$} we denote the centralizer of $M$ in
$G$, by $Z(G)$\glossary{Z(G)@$Z(G)$} we denote the center of $G$. The conjugate
of $x$ by an element $y$ in $G$ is written as
$x^y=y^{-1}xy$\glossary{xy@$x^y=y^{-1}xy$}
(${}^yx=yxy^{-1}$\glossary{xyinverse@${}^yx=yxy^{-1}$}), by
$\left[x,y\right]=x^{-1}x^y$\glossary{x,y@$[x,y]=x^{-1}y^{-1}xy$} we denote the
commutator of $x,y$. The symbol $\left[A,B\right]$\glossary{ABcommutant@$[A,B]$}
means the mutual commutant of subgroups $A$ and $B$ of $G$. For groups $A$ and
$B$ the expressions $A\times B$\glossary{AtimesB@$A\times B$}, $A\circ
B$\glossary{AcircB@$A\circ B$}, and $A\rightthreetimes
B$\glossary{ArigthtreetimesB@$A\rightthreetimes B$} mean direct, central, and
semidirect products respectively of $A$ and $B$ with a normal subgroup
$B$. If $A$ and $B$ are subgroups of $G$ such that $A\unlhd B$, then the factor
group $B/A$ is called a {\em section}\index{секция} of $G$. The Fitting
subgroup of $G$ is denoted by $F(G)$\glossary{F(G)@$F(G)$}, the generalized
Fitting subgroup is denoted by~$F^\ast(G)$\glossary{Fast(G)@$F^\ast(G)$}.

The set of Sylow $p$-sub\-gro\-ups of a finite group $G$ we shall denote by
$Syl_p(G)$\glossary{Sylp@$Syl_p(G)$}. If$\varphi$ is a homomorphism of
$G$, $g$ is an element of $G$, then $G^\varphi$, $g^\varphi$ are the images
of $G$ and $g$ under $\varphi$ respectively {Gphup@$G^\varphi$, $g^\varphi$}.
By
$G_\varphi$\glossary{Gphi@$G_\varphi=\{g\in G\vert g^\varphi=g\}$} we denote
the set of stable points of $G$  under an endomorphism $\varphi$. By
$\Aut(G)$\glossary{Aut(G)@$\Aut(G)$}, $\Out(G)$\glossary{Out(G)@$\Out(G)$} and
$\Inn(G)$\glossary{Inn(G)@$\Inn(G)$} we denote the group of all automorphisms,
the group of outer automorphisms, and the group of inner automorphisms
of $G$ respectively. If $G$ is a group, we denote by $\P G$\glossary{PG@$\P G$}
the factor group $G/Z(G)$. An isomorphism $\P G\simeq \Inn(G)$ is known, in
particular, if $Z(G)$ is trivial, then $G\simeq \Inn(G)$, and we may assume
that~$G\leq \Aut(G)$. A finite group $G$ is said to be {\em almost
simple}\index{group!almost simple}, if there exists a finite group
$S$ with $S\leq G\leq \Aut(S)$, i.~e., $F^\ast(G)=S$ is a simple group. For
every positive integer $t$ by $\Z_t$\glossary{Zt@$\Z_t$} we denote a cyclic
group of order~$t$.

If $\pi$ is a set of primes, then by $\pi'$\glossary{piprime@$\pi'$} we denote
its complement in the set of all primes. For every positive integer $n$ by
$\pi(n)$\glossary{pin@$\pi(n)$} we denote the set of prime divisors of $n$, and
by $n_\pi$\glossary{npi@$n_\pi$} we denote the maximal divisor of $n$ such
that~${\pi(n_\pi)\subseteq \pi}$. As usual we denote by
$O_\pi(G)$\glossary{OpidownG@$O_\pi(G)$} the maximal normal $\pi$-sub\-gro\-up
of
$G$, and we denote by
$O^{\pi'}(G)$\glossary{OpiupG@$O^{\pi'}(G)$} the subgroup generated by all
$\pi$-ele\-ments of $G$. If $\pi=\{2\}'$ is a set of all odd primes, then
$O_\pi(G)=O_{2'}(G)$ is denoted by~$O(G)$\glossary{OG@$O(G)$}. If
$g\in G$, then by $g_\pi$\glossary{gpidown@$g_\pi$} we denote the
$\pi$-part of $g$, i.~e., $g_\pi=g^{\vert g\vert_{\pi'}}$.

Let $G$ be a group, $A,B,H$ be subgroups of $G$ and $B$ is normal in $A$. Then
$N_H(A/B)=N_H(A)\cap N_H(B)$\glossary{NH(A/B)@$N_H(A/B)$}. If
$x\in N_H(A/B)$, then $x$ induces an automorphism
$Ba\mapsto B x^{-1}ax$ of $A/B$. Thus there exists a homomorphism of $N_H(A/B)$
into $\Aut(A/B)$. The image of this homomorphism is denoted by
$\Aut_H(A/B)$\glossary{AutH(A/B)@$\Aut_H(A/B)$} and is called a {\em group of
induced automorphisms}\index{group!of induced automorphisms} of $H$ on the
section $A/B$. In particular, if $S=A/B$ is a composition factor of
$G$, then for each subgroup $H\leq G$ the group~$\Aut_H(S)=\Aut_H(A/B)$ is
defined. Note that the structure of $\Aut_H(S)$ depends on the choice of a
composition series. If $A,H$ are subgroups of $G$, then
$\Aut_H(A)=\Aut_H(A/\{e\})$ by definition.

\subsection{Linear algebraic groups}

Necessary information about the structure and properties of linear algebraic
groups can be found in \cite{HuLinearAlgGrps}. Since we consider linear
algebraic groups only, we shall omit the word ``linear'' for brevity.

If $\ov{G}$ is an algebraic group, then by $\ov{G}^0$\glossary{G0@$\ov{G}^0$}
we denote the unit component of $\ov{G}$. An algebraic group is
called {\em semisimple}\index{group!algebraic!semisimple}, if its radical
$R(\ov{G})$\glossary{R(G)@$R(\ov{G})$} is trivial, and an algebraic group is
called {\em reductive}\index{group!algebraic!reductive}, if its unipotent
radical $R_u(\ov{G})$\glossary{Ru(G)@$R_u(\ov{G})$} is trivial (in both cases
an algebraic group is assumed to be finite). A connected semisimple algebraic
group is known (for example, see \cite[Theorem~27.5]{HuLinearAlgGrps}) to be a
central product of connected simple algebraic groups, while a
connected reductive algebraic group  $\ov{G}$ is known to be a central product
of a torus $\ov{S}$ and a semisimple group $\ov{M}$ with
$\ov{S}=Z(\ov{G})^0$,
$\ov{M}=\left[\ov{G},\ov{G}\right]$, and $\ov{S}\cap \ov{M}$ is finite.

If $\ov{G}$ is a connected reductive algebraic group, then let $\ov{T}$ be
its maximal torus (by a {\em torus}\index{torus!of algebraic group} we always
mean a connected diagonalizable ($d$-) group). The dimension of a maximal torus
is called a {\em rank}\index{rank!of algebraic group} of an algebraic group.
By $\Phi(\ov{G})$\glossary{Phi(ovG)@$\Phi(\ov{G})$} the root system of
\index{system!of roots of algebraic group} of $\ov{G}$ with respect to a
maximal torus $\ov{T}$ (it does not depend on the choice of a maximal torus) is
denoted, and $W(\ov{G})\simeq
N_{\ov{G}}(\ov{T})/\ov{T}$\glossary{W(ovG)@$W(\ov{G})$} is the {\em Weyl
group}\index{Weyl group!of algebraic group} of $G$. If $\ov{G}$ is a reductive
group of rank $n$, then the dimension of the centralizer of any its element is
not less than $n$. An element is called {\em regular}\index{element regular} if
the dimension of its centralizer is equal to $n$. In particular, a semisimple
element $s$ is regular, if $C_{\ov{G}}(s)^0$ is a maximal torus of~$\ov{G}$.

Recall that or every root system $\Phi$ there exists a set of roots
$r_1,\ldots,r_n$ such that each root of $\Phi$ can be uniquely written as
$\sum_{i=1}^n\alpha_ir_i$, where all coefficients $\alpha_i$ are integers, and
either nonnegative, or nonpositive. Such a set of roots is called  a
{\em fundamental set}\index{set!of roots fundamental} of $\Phi$, and its
elements are called {\em fundamental
roots}\index{root!fundamental}.  At that a fundamental set is a basis of
$\Z\Phi\otimes_{\Z}\R$. The dimension of
$\Z\Phi\otimes_{\Z}\R$ is called a {\em rank} of
$\Phi$\index{rank!of root system}. Note that the ranks of $\ov{G}$ and of its
root system $\Phi(\ov{G})$ are equal. Below we assume that all fundamental
roots are positive. Then a root $r$ is {\em
positive}\index{root!positive} if and only if it is a linear combination of
fundamental roots with nonnegative coefficients. For a root system $\Phi$ by
$\Phi^+$\glossary{Phiplus@$\Phi^+$}
($\Phi^-$)\glossary{Phiminus@$\Phi^-$} the set of all positive (negative) roots
is denoted. The number $h(r)=\sum_{i=1}^n\alpha_i$\glossary{h(r)@$h(r)$} is
called a {\em height}\index{height of root $h(r)$} of $r=\sum_{i=1}^n
\alpha_i r_i$. In every irreducible root system  $\Phi$ there exists a
unique root of maximal height\index{root!of maximal height}, which is denoted
by~$r_0$\glossary{r0@$r_0$} below. Note that the Weyl group
$W(\Phi)$\glossary{W(Phi)@$W(\Phi)$} of a root system $\Phi$ is generated by
reflections in fundamental roots, which are called  {\em fundamental
reflections}\index{reflection fundamental}. If we denote by
$l(w)$\glossary{l(w)@$l(w)$} the minimal number of multipliers in a
decomposition of $w$ into the product of fundamental reflections, a ({\em
length}\index{length of element of Weyl group $l(w)$}), then there exists a
unique element of maximal length, denoted by~$w_0$\glossary{w0@$w_0$} below,
that is a unique element of the Weyl group mapping all positive roots into
negative roots. In general,  $l(w)$ is equal to $\vert \Phi^- \cap
(\Phi^+)^w\vert$, i.~e., to the number of positive roots, that $w$ maps into
negative roots.

Let  $\ov{G}$ be a connected simple algebraic group,  $\pi$ be its exact
rational representation, $\Gamma_\pi$ be a
lattice\glossary{Gammapi@$\Gamma_\pi$} generated by weights of the
representation $\pi$. By $\Gamma_{ad}$\glossary{Gammaad@$\Gamma_{ad}$} we
denote the lattice generated by the roots of $\Phi$, by $\Gamma_{sc}$ we denote
the lattice \glossary{Gammasc@$\Gamma_{sc}$} generated by the fundamental
weights. The lattices  $\Gamma_{sc}$, $\Gamma_\pi$, and  $\Gamma_{ad}$ do not
depend on the representation of $\ov{G}$, and the following inclusions
$\Gamma_{ad}\leq\Gamma_\pi\leq\Gamma_{sc}$ (see
\cite[31.1]{HuLinearAlgGrps}) hold. Several distinct algebraic groups, which
are called {\em isogenies}\index{isogeny}, is known to exist for a given root
system. They differs by the structure of $\Gamma_\pi$ and the order of the
finite center. If $\Gamma_\pi$ coinsides with $\Gamma_{sc}$, a group $\ov{G}$
is said to be {\em simply connected}\index{group!algebraic!simply connected}, it
is denoted by  $\ov{G}_{sc}$\glossary{Gsc@$\ov{G}_{sc}$}. If
$\Gamma_\pi$ coinsides with $\Gamma_{ad}$, a group  $\ov{G}$ is said to have an
{\em adjoint type}\index{group!algebraic!adjoint}, it is denoted by через
$\ov{G}_{ad}$\glossary{Gad@$\ov{G}_{ad}$}. Every linear algebraic group with
a root system $\Phi$ can be obtained as a factor group of $\ov{G}_{sc}$ by a
subgroup of its center. The center of $\ov{G}_{ad}$  is trivial, and this group
is simple as an abstract group. The factor group  $\Gamma_{sc}/\Gamma_\pi$ is
denoted by~$\Delta(\ov{G})$\glossary{DeltaG@$\Delta(\ov{G})$} and is called a
{\em fundamental group}\index{group!fundamental of algebraic group
$\Delta(\ov{G})$} of~$\ov{G}$. The factor group
$\Gamma_{sc}/\Gamma_{ad}$ depends on the root system $\Phi$ only and is denoted
by $\Delta(\Phi)$\glossary{DeltaPhi@$\Delta(\Phi)$}. Хорошо известно,
The group  $\Delta(\Phi)$ is known to be cyclic, except the root system
$\Phi=D_{2n}$, when $\Delta(D_{2n})=\Z_2\times \Z_2$ is elementary Abelian of
order~$4$.

Let $\ov{B}$ be a Borel subgroup\glossary{Bov@$\ov{B}$},
$\ov{T}\leq \ov{B}$ be a maximal torus, and
$\ov{U}=R_u(\ov{B})$\glossary{Uov@$\ov{U}$} be a maximal connected unipotent
subgroup of $\ov{G}$. There exists a unique Borel subgroup
$\ov{B}^-$\glossary{Bovmin@$\ov{B}^-$} such that $\ov{B}\cap
\ov{B}^-=\ov{T}$, denote by
$\ov{U}^-=R_u(\ov{B}^-)$\glossary{Uovmin@$\ov{U}^-$}. If we fix an order
on $\Phi(\ov{G})$, then each element $u\in \ov{U}$ (respectively
$u\in\ov{U}^-$) can be uniquely written in the form
\begin{equation}\label{canonicalform}
u=\prod_{r\in\Phi^+}x_r(t_r)
\end{equation}
\glossary{xr(t)@$x_r(t_r)$}(respectively
$u=\prod_{r\in\Phi^-}x_r(t_r)$), where roots are taken in given order, elements
$t_r$-s are in the definition field of~$\ov{G}$ and
$\{X_r,r\in\Phi\}$ \glossary{Xr@$X_r$} is a set of {\em 1-dimensional
$\ov{T}$-in\-va\-ri\-ant subgroups}\index{subgroup!1-dimensional
$\ov{T}$-in\-va\-ri\-ant} (a set of {\em root
subgroups}\index{subgroup!root of algebraic group}). The multiplication of
elements from distinct root subgroups is defined by {\em Chevalley commutator
formulae}\index{Chevalley commutator formulae}.
\begin{lem} \label{ChevalleyFormulae}  {\em \cite[5.2.2]{CarSimpleGrpsLieType},
(Chevalley commutator formulae)}
Let $x_r(t)$, $x_s(u)$ be elements from distinct root subgroups  $X_r$ and
$X_s$ respectively,  $r\ne - s$. Then
$$[x_r(t),x_s(u)]=\prod_{ir+js\in\Phi;i,j>0}x_{ir+js}(C_{ijrs}(-t)^iu^j),$$
where the constants $C_{ijrs}$ do not depend on $t$ and~$u$.
\end{lem}
\noindent Substantially this formulae means that the mutual commutant of
$X_r$ and $X_s$ is in the group generated by subgroups $X_{ir+js}$, where
$i,j>0$ and~${ir+js\in\Phi}$.

Let $c_i$ be the coefficient of a fundamental root  $r_i$ in the decomposition
of $r_0$. Primes, dividing $c_i$-s, are called {\em bad}\index{bad prime}
primes. The diagram, obtained from the Dynkin diagram by addition of $-r_0$ and
its connection with other fundamental roots by usual rule, is called an {\em
extended Dynkin diagram}\index{extended Dynkin diagram}. Let $\ov{R}$ be a
(connected) reductive subgroup of maximal rank of a connected simple algebraic
group  $\ov{G}$. As we already noted, in this case
$\ov{R}=\ov{G}_1\circ\ldots\circ \ov{G}_k\circ
Z(\ov{R})^0$, where  $\ov{G}_i$-s are connected simple algebraic groups of
rank less, than the rank of $\ov{G}$. More over, if
$\Phi_1,\ldots,\Phi_k$ are root systems of
$\ov{G}_1,\ldots,\ov{G}_k$ respectively, then $\Phi_1\cup\ldots\cup\Phi_k$ is a
subsystem of $\Phi$. There exists a nice algorithm, due to Borel and de
Siebental \cite{BorSieSubsystems} and, independently, Dynkin
\cite{DynSubsystems}, of determining subsystems of a root system. One needs to
extend the Dynkin diagram to the extended Dynkin diagram, remove some vertices
from it, and repeat the procedure for obtained connected components. Diagrams
obtained in this way are subsystem diagrams and diagram of any subsystem can be
obtained in this way.

In Table \ref{Dyndiagr} we give extended Dynkin diagrams of all
irreducible root systems and coefficients of fundamental roots in the
decomposition of $r_0$. The numberring in  Table~\ref{Dyndiagr} is chosen
as in~\cite{der1CentSemisimpleExcept}.

{\footnotesize
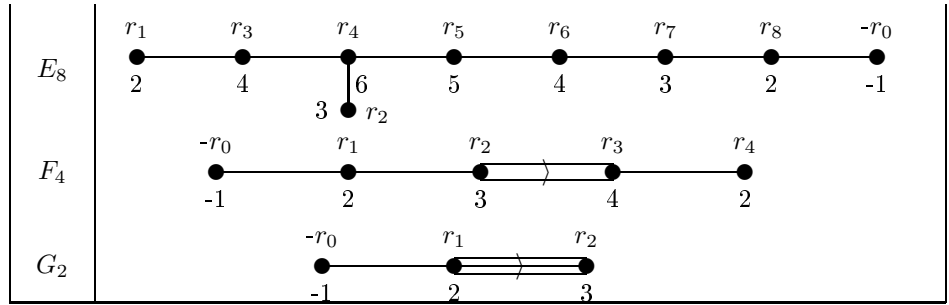
\begin{longtable}
{|c|c|} \caption{Root systems and extended Dynkin diagrams\label{Dyndiagr}} \\
\hline
$\Phi$&Extended Dynkin diagram\\ \hline &\\
\begin{picture}(20,50)\put(10,25){\makebox(0,0){$A_n$}}
\end{picture}
& \begin{picture}(310,50)(-25,-10) \put(0,0){\circle*{6}}
\put(0,0){\line(1,0){50}}
\put(50,0){\circle*{6}} \put(50,0){\line(1,0){50}} \put(100,0){\circle*{6}}
\put(100,0){\line(1,0){50}} \put(152.5,0){\circle*{1}} \put(155,0){\circle*{1}}
\put(157.5,0){\circle*{1}} \put(160,0){\circle*{1}} \put(162.5,0){\circle*{1}}
\put(165,0){\circle*{1}} \put(167.5,0){\circle*{1}} \put(170,0){\circle*{1}}
\put(172.5,0){\circle*{1}} \put(175,0){\circle*{1}} \put(177.5,0){\circle*{1}}
\put(180,0){\circle*{1}} \put(182.5,0){\circle*{1}} \put(185,0){\circle*{1}}
\put(187.5,0){\circle*{1}} \put(190,0){\circle*{1}} \put(192.5,0){\circle*{1}}
\put(195,0){\circle*{1}} \put(197.5,0){\circle*{1}} \put(250,0){\circle*{6}}
\put(200,0){\line(1,0){50}} \put(0,10){\makebox(0,0){$r_1$}}
\put(0,-10){\makebox(0,0){1}} \put(50,10){\makebox(0,0){$r_2$}}
\put(50,-10){\makebox(0,0){1}} \put(100,10){\makebox(0,0){$r_3$}}
\put(100,-10){\makebox(0,0){1}} \put(250,10){\makebox(0,0){$r_n$}}
\put(250,-10){\makebox(0,0){1}} \put(0,0){\line(3,1){125}}
\put(250,0){\line(-3,1){125}} \put(125,40){\circle*{6}}
\put(125,50){\makebox(0,0){$-r_0$}} \put(125,30){\makebox(0,0){-1}}
\end{picture}\\

&\\
\begin{picture}(20,50)\put(10,25){\makebox(0,0){$B_n$}}
\end{picture}&

\begin{picture}(310,60)(-25,-30) \put(0,20){\circle*{6}}
\put(0,20){\line(5,-2){50}} \put(50,0){\circle*{6}} \put(50,0){\line(1,0){50}}
\put(100,0){\circle*{6}} \put(100,0){\line(1,0){25}} \put(147.5,0){\circle*{1}}
\put(150,0){\circle*{1}} \put(152.5,0){\circle*{1}} \put(160,0){\circle*{1}}
\put(157.5,0){\circle*{1}} \put(165,0){\circle*{1}} \put(162.5,0){\circle*{1}}
\put(170,0){\circle*{1}} \put(167.5,0){\circle*{1}} \put(145,0){\circle*{1}}
\put(172.5,0){\circle*{1}} \put(130,0){\circle*{1}} \put(127.5,0){\circle*{1}}
\put(135,0){\circle*{1}} \put(132.5,0){\circle*{1}} \put(140,0){\circle*{1}}
\put(137.5,0){\circle*{1}} \put(147.5,0){\circle*{1}} \put(142.5,0){\circle*{1}}
\put(155,0){\circle*{1}} \put(175,0){\line(1,0){25}} \put(200,0){\circle*{6}}
\put(250,0){\circle*{6}} \put(200,3){\line(1,0){50}}
\put(200,-3){\line(1,0){50}}
\put(0,30){\makebox(0,0){$r_1$}} \put(0,10){\makebox(0,0){1}}
\put(50,10){\makebox(0,0){$r_2$}} \put(50,-10){\makebox(0,0){2}}
\put(100,10){\makebox(0,0){$r_3$}} \put(100,-10){\makebox(0,0){2}}
\put(250,10){\makebox(0,0){$r_n$}} \put(250,-10){\makebox(0,0){2}}
\put(0,-20){\circle*{6}} \put(0,-20){\line(5,2){50}}
\put(0,-10){\makebox(0,0){$-r_0$}} \put(0,-30){\makebox(0,0){-1}}
\put(200,10){\makebox(0,0){$r_{n-1}$}} \put(200,-10){\makebox(0,0){2}}
\put(225,0){\makebox(0,0){$\rangle$}}
\end{picture}\\

\begin{picture}(20,20)\put(10,10){\makebox(0,0){$C_n$}}
\end{picture}&

\begin{picture}(310,40)(-25,-20) \put(0,0){\circle*{6}}
\put(0,3){\line(1,0){50}} \put(0,-3){\line(1,0){50}} \put(50,0){\circle*{6}}
\put(50,0){\line(1,0){50}} \put(100,0){\circle*{6}} \put(100,0){\line(1,0){25}}
\put(147.5,0){\circle*{1}} \put(150,0){\circle*{1}} \put(152.5,0){\circle*{1}}
\put(160,0){\circle*{1}} \put(157.5,0){\circle*{1}} \put(165,0){\circle*{1}}
\put(162.5,0){\circle*{1}} \put(170,0){\circle*{1}} \put(167.5,0){\circle*{1}}
\put(145,0){\circle*{1}} \put(172.5,0){\circle*{1}} \put(130,0){\circle*{1}}
\put(127.5,0){\circle*{1}} \put(135,0){\circle*{1}} \put(132.5,0){\circle*{1}}
\put(140,0){\circle*{1}} \put(137.5,0){\circle*{1}} \put(147.5,0){\circle*{1}}
\put(142.5,0){\circle*{1}} \put(155,0){\circle*{1}} \put(175,0){\line(1,0){25}}
\put(200,0){\circle*{6}} \put(250,0){\circle*{6}} \put(200,3){\line(1,0){50}}
\put(200,-3){\line(1,0){50}} \put(0,10){\makebox(0,0){$-r_0$}}
\put(0,-10){\makebox(0,0){-1}} \put(50,10){\makebox(0,0){$r_1$}}
\put(50,-10){\makebox(0,0){2}} \put(100,10){\makebox(0,0){$r_2$}}
\put(100,-10){\makebox(0,0){2}} \put(250,10){\makebox(0,0){$r_n$}}
\put(250,-10){\makebox(0,0){1}} \put(200,10){\makebox(0,0){$r_{n-1}$}}
\put(200,-10){\makebox(0,0){2}} \put(225,0){\makebox(0,0){$\langle$}}
\put(25,0){\makebox(0,0){$\rangle$}}
\end{picture}\\

\begin{picture}(20,50)\put(10,25){\makebox(0,0){$D_n$}}
\end{picture}&

\begin{picture}(310,60)(-25,-30) \put(0,20){\circle*{6}}
\put(0,20){\line(5,-2){50}} \put(50,0){\circle*{6}} \put(50,0){\line(1,0){50}}
\put(100,0){\circle*{6}} \put(100,0){\line(1,0){25}} \put(147.5,0){\circle*{1}}
\put(150,0){\circle*{1}} \put(152.5,0){\circle*{1}} \put(160,0){\circle*{1}}
\put(157.5,0){\circle*{1}} \put(165,0){\circle*{1}} \put(162.5,0){\circle*{1}}
\put(170,0){\circle*{1}} \put(167.5,0){\circle*{1}} \put(145,0){\circle*{1}}
\put(172.5,0){\circle*{1}} \put(130,0){\circle*{1}} \put(127.5,0){\circle*{1}}
\put(135,0){\circle*{1}} \put(132.5,0){\circle*{1}} \put(140,0){\circle*{1}}
\put(137.5,0){\circle*{1}} \put(147.5,0){\circle*{1}} \put(142.5,0){\circle*{1}}
\put(155,0){\circle*{1}} \put(175,0){\line(1,0){25}} \put(200,0){\circle*{6}}
\put(250,20){\circle*{6}} \put(250,-20){\circle*{6}} \put(200,0){\line(5,2){50}}
\put(200,0){\line(5,-2){50}} \put(0,30){\makebox(0,0){$r_1$}}
\put(0,10){\makebox(0,0){1}} \put(50,10){\makebox(0,0){$r_2$}}
\put(50,-10){\makebox(0,0){2}} \put(100,10){\makebox(0,0){$r_3$}}
\put(100,-10){\makebox(0,0){2}} \put(250,30){\makebox(0,0){$r_{n-1}$}}
\put(250,10){\makebox(0,0){1}} \put(0,-20){\circle*{6}}
\put(0,-20){\line(5,2){50}}
\put(0,-10){\makebox(0,0){$-r_0$}} \put(0,-30){\makebox(0,0){-1}}
\put(200,10){\makebox(0,0){$r_{n-2}$}} \put(200,-10){\makebox(0,0){2}}
\put(250,-10){\makebox(0,0){$r_n$}} \put(250,-30){\makebox(0,0){1}}
\end{picture}\\

\begin{picture}(20,60)\put(10,30){\makebox(0,0){$E_6$}}
\end{picture}&

\begin{picture}(310,60)(-45,-50)
\put(0,0){\line(1,0){50}} \put(50,0){\line(1,0){50}} \put(100,0){\line(1,0){50}}
\put(150,0){\line(1,0){50}} \put(0,0){\circle*{6}} \put(50,0){\circle*{6}}
\put(100,0){\circle*{6}} \put(150,0){\circle*{6}} \put(200,0){\circle*{6}}
\put(100,0){\line(0,-1){20}} \put(100,-20){\line(0,-1){20}}
\put(100,-40){\circle*{6}} \put(100,-20){\circle*{6}}
\put(0,-10){\makebox(0,0){1}}
\put(0,10){\makebox(0,0){$r_1$}} \put(50,-10){\makebox(0,0){2}}
\put(50,10){\makebox(0,0){$r_3$}} \put(105,-10){\makebox(0,0){3}}
\put(100,10){\makebox(0,0){$r_4$}} \put(150,-10){\makebox(0,0){2}}
\put(150,10){\makebox(0,0){$r_5$}} \put(200,-10){\makebox(0,0){1}}
\put(200,10){\makebox(0,0){$r_6$}} \put(90,-20){\makebox(0,0){2}}
\put(112,-22){\makebox(0,0){$r_2$}} \put(90,-40){\makebox(0,0){-1}}
\put(112,-42){\makebox(0,0){-$r_0$}}
\end{picture}\\

\begin{picture}(20,50)\put(10,25){\makebox(0,0){$E_7$}}
\end{picture}&

\begin{picture}(310,40)(-5,-30)
\put(0,0){\line(1,0){50}} \put(50,0){\line(1,0){50}} \put(100,0){\line(1,0){50}}
\put(150,0){\line(1,0){50}} \put(200,0){\line(1,0){50}}
\put(250,0){\line(1,0){50}}
\put(0,0){\circle*{6}} \put(50,0){\circle*{6}} \put(100,0){\circle*{6}}
\put(150,0){\circle*{6}} \put(200,0){\circle*{6}} \put(250,0){\circle*{6}}
\put(300,0){\circle*{6}} \put(150,0){\line(0,-1){20}} \put(150,-20){\circle*{6}}
\put(0,10){\makebox(0,0){-$r_0$}}
\put(0,-10){\makebox(0,0){-1}}\put(50,10){\makebox(0,0){$r_1$}}
\put(50,-10){\makebox(0,0){2}}\put(100,10){\makebox(0,0){$r_3$}}
\put(100,-10){\makebox(0,0){3}}\put(150,10){\makebox(0,0){$r_4$}}
\put(155,-10){\makebox(0,0){4}}\put(200,10){\makebox(0,0){$r_5$}}
\put(200,-10){\makebox(0,0){3}}\put(250,10){\makebox(0,0){$r_6$}}
\put(250,-10){\makebox(0,0){2}}\put(300,10){\makebox(0,0){$r_7$}}
\put(300,-10){\makebox(0,0){1}}\put(161,-22){\makebox(0,0){$r_2$}}
\put(140,-20){\makebox(0,0){2}}
\end{picture}\\

&\\
\begin{picture}(20,50)\put(10,25){\makebox(0,0){$E_8$}}
\end{picture}&

\begin{picture}(310,40)(-10,-30)
\put(0,0){\line(1,0){40}} \put(40,0){\line(1,0){40}} \put(80,0){\line(1,0){40}}
\put(120,0){\line(1,0){40}} \put(160,0){\line(1,0){40}}
\put(200,0){\line(1,0){40}}
\put(0,0){\circle*{6}} \put(40,0){\circle*{6}} \put(80,0){\circle*{6}}
\put(120,0){\circle*{6}} \put(160,0){\circle*{6}} \put(200,0){\circle*{6}}
\put(240,0){\circle*{6}} \put(80,0){\line(0,-1){20}}
\put(80,-20){\circle*{6}}
\put(0,10){\makebox(0,0){$r_1$}}
\put(0,-10){\makebox(0,0){2}}\put(40,10){\makebox(0,0){$r_3$}}
\put(40,-10){\makebox(0,0){4}}\put(80,10){\makebox(0,0){$r_4$}}
\put(85,-10){\makebox(0,0){6}}\put(120,10){\makebox(0,0){$r_5$}}
\put(120,-10){\makebox(0,0){5}}\put(160,10){\makebox(0,0){$r_6$}}
\put(160,-10){\makebox(0,0){4}}\put(200,10){\makebox(0,0){$r_7$}}
\put(200,-10){\makebox(0,0){3}}\put(240,10){\makebox(0,0){$r_8$}}
\put(240,-10){\makebox(0,0){2}}\put(91,-22){\makebox(0,0){$r_2$}}
\put(70,-20){\makebox(0,0){3}} \put(240,0){\line(1,0){40}}
\put(280,0){\circle*{6}}
\put(280,-10){\makebox(0,0){-1}}\put(280,10){\makebox(0,0){-$r_0$}}
\end{picture}\\

\begin{picture}(20,20)\put(10,10){\makebox(0,0){$F_4$}}
\end{picture}&

\begin{picture}(310,20)(-40,-10)
\put(0,0){\line(1,0){50}} \put(50,0){\line(1,0){50}} \put(100,3){\line(1,0){50}}
\put(100,-3){\line(1,0){50}} \put(150,0){\line(1,0){50}} \put(0,0){\circle*{6}}
\put(50,0){\circle*{6}} \put(100,0){\circle*{6}} \put(150,0){\circle*{6}}
\put(200,0){\circle*{6}} \put(0,-10){\makebox(0,0){-1}}
\put(0,10){\makebox(0,0){-$r_0$}} \put(50,-10){\makebox(0,0){2}}
\put(50,10){\makebox(0,0){$r_1$}} \put(100,-10){\makebox(0,0){3}}
\put(100,10){\makebox(0,0){$r_2$}} \put(150,-10){\makebox(0,0){4}}
\put(150,10){\makebox(0,0){$r_3$}} \put(200,-10){\makebox(0,0){2}}
\put(200,10){\makebox(0,0){$r_4$}} \put(125,0){\makebox(0,0){$\rangle$}}
\end{picture}\\

&\\
\begin{picture}(20,20)\put(10,10){\makebox(0,0){$G_2$}}
\end{picture}&

\begin{picture}(310,20)(-80,-10)
\put(0,0){\line(1,0){50}} \put(50,0){\line(1,0){50}} \put(50,3){\line(1,0){50}}
\put(50,-3){\line(1,0){50}} \put(50,0){\line(1,0){50}} \put(0,0){\circle*{6}}
\put(50,0){\circle*{6}} \put(100,0){\circle*{6}}
\put(0,-10){\makebox(0,0){-1}}\put(0,10){\makebox(0,0){-$r_0$}}
\put(50,-10){\makebox(0,0){2}}\put(50,10){\makebox(0,0){$r_1$}}
\put(100,-10){\makebox(0,0){3}}\put(100,10){\makebox(0,0){$r_2$}}
\put(75,0){\makebox(0,0){$\rangle$}}
\end{picture}\\ \hline
\end{longtable}}

For every semisimple element $s\in \ov{G}$, where $\ov{G}$ is a connected
reductive group, the unit component $C_{\ov{G}}(s)^0$ is a reductive subgroup
of maximal rank and  $C_{\ov{G}}(s)/C_{\ov{G}}(s)^0\simeq D\leq
\Delta(\ov{G})$ (see Lemma \ref{FactorCentrByConnectedComponent} below).

\subsection{Structure of finite groups of Lie type}

Our notation and definitions for finite groups of Lie type mainly agree with
that of \cite{CarSimpleGrpsLieType} (except the definition of  finite groups of
Lie type, see below). If $G$ is a finite group of Lie type with the trivial
center (we do not exclude non-simple groups of Lie type, such as $A_1(2)$, all
exceptions are given in \cite[Theorems 11.1.2 and~14.4.1]{CarSimpleGrpsLieType}
and cited i Table \ref{notsimple} below), then
$\widehat{G}$\glossary{Gwidehat@$\widehat{G}$} denoted the group of
inner-diagonal automorphisms of $G$. In view of
\cite[3.2]{Ste2} we have that $\Aut(G)$ is generated by inner-diagonal, field,
and graph automorphisms. Note that the definition of a field and a graph
automorphisms in the present paper is slightly different from the definitions
given in \cite{Ste2}, precise definitions are given in subsection 4.1. Since we
are
assuming that $Z(G)$ is trivial, then $G$ is isomorphic to the group of
its inner automorphisms, and so we may assume that~$G\leq\widehat{G}\leq
\Aut(G)$.

\begin{longtable}{|c|l|}\caption{Groups of Lie type which are not
simple\label{notsimple}}\\ \hline
Group&Properties\\ \hline
$A_1(2)$&Group is solvable\\
$A_1(3)$&Group is solvable\\
$B_2(2)$&$B_2(2)\simeq\mathrm{Sym}_6$\\
$G_2(2)$&$[G_2(2),G_2(2)]\simeq {}^2A_2(3)$\\
${}^2A_2(2)$&Group is solvable\\
${}^2B_2(2)$&Group is solvable\\
${}^2G_2(3)$&$[{}^2G_2(3),{}^2G_2(3)]\simeq A_1(8)$\\
${}^2F_4(2)$&$[{}^2F_4(2),{}^2F_4(2)]$ is the simple Tits group\\ \hline
\end{longtable}

Let $\overline{G}$ be a simple connected algebraic group over an algebraic
closure $\ov{\F}_p$\glossary{Fp@$\ov{\F}_p$} of a finite field of positive
characteristic  $p$. Here $Z(\ov{G})$ can be nontrivial. An endomorphism
$\sigma$ of $\ov{G}$ is called a {\em Frobenius
map}\index{Frobenius map}, if $\ov{G}_\sigma$ is finite, while the kernel of
$\sigma$ is trivial (i.~e., $\sigma$ is an automorphism of $G$ as an abstract
group). Groups $O^{p'}(\ov{G}_\sigma)$ are called {\em canonical finite groups
of Lie type}\index{group!of Lie type! canonical finite}, and every group $G$,
with
$O^{p'}(\ov{G}_\sigma)\leq G\leq \ov{G}_\sigma$, is called a {\em finite group
of Lie type}\index{group!of Lie type!finite}. If $\ov{G}$ is a simple algebraic
group of adjoint type, then we shall say that $G$ also has an {\em adjoint
type}\index{group!of Lie type!finite adjoint}. Note that in
\cite{CarSimpleGrpsLieType} only groups $O^{p'}(\ov{G})$ are called finite
groups of Lie type. But later in
\cite{CarFiniteGrpsLieTypeConjClassesCharacters} R.Carter says that every group
$\ov{G}_\sigma$ is a finite group of Lie type, for every connected reductive
group $\ov{G}$. More over, in \cite{Ca2CentsSemisimpleClassical}
and \cite{der1CentSemisimpleExcept}, without any comment, every group
$G$ with  $O^{p'}(\ov{G}_\sigma)\leq G\leq \ov{G}_\sigma$ is called a finite
group of Lie type. Thus giving the definition of finite groups of Lie type
and of finite canonical groups of Lie type we intend to clarify the situation
here. For example,  $\P \SL_2(3)$ is a canonical finite group of Lie type and
$\P \GL_2(3)$ is a finite group of Lie type. Note that an element of order
$3$ is not conjugate to its inverse in  $\P \SL_2(3)$ and is conjugate to its
inverse in $\P \GL_2(3)$. Since such information about the conjugation
is important in many cases (and is very important and useful in the paper), we
find it reasonable to use such notation.

In general, for a given group of Lie type $G$ (if we consider it as an
abstract group) a corresponding algebraic group is not uniquely
defined. For example, if $G= \P \SL_2(5)\simeq \SL_2(4)$, then $G$ can be
obtained either as $(\SL_2(\ov{\F}_2))_\sigma$, or as $O^{5'}((\P
\GL_2(\ov{\F}_5))_\sigma)$ (for suitable $\sigma$-s). Hence, for every finite
group of Lie type  $G$ we fix (in some way) a corresponding algebraic group
$\ov{G}$ and a Frobenius map~$\sigma$ such that
$O^{p'}(\ov{G}_\sigma)\leq G\leq \ov{G}_\sigma$.

We say that groups ${}^2A_n(q)$, ${}^2D_n(q)$,
${}^2E_6(q)$ are defined over $\F_{q^2}$, groups ${}^3 D_4(q)$
are defined over $\F_{q^3}$, and  the remaining groups are defined over
$\F_q$\index{field!of definition of group of Lie type}. The field $\F_q$ in
all cases is called a {\em base field}\index{field!base of group of Lie type}.
In view of \cite[Lemma~2.5.8]{GorLySol}, if $\ov{G}$ is of adjoint type,
then $\ov{G}_\sigma$ is the group of in\-ner-di\-a\-go\-nal automorphisms of
$O^{p'}(\ov{G}_\sigma)$. If $\ov{G}$ is simply connected, then
$\ov{G}_\sigma=O^{p'}(\ov{G}_\sigma)$ (see
\cite[12.4]{SteLangSteinbergTheorem}). In any case, in view of
\cite[Theorem~2.2.6(g)]{GorLySol} $\ov{G}_\sigma=\ov{T}_\sigma
O^{p'}(\ov{G}_\sigma)$ for every $\sigma$-stab\-le maximal torus $\ov{T}$ of
$\ov{G}$. Let $U\leq \langle X_r\vert r\in \Phi^+\rangle=\ov{U}$ be a
maximal unipotent subgroup of  $G$ (at that $\ov{U}$ is a maximal connected
$\sigma$-stab\-le unipotent subgroup of $\ov{G}$). Then each $u\in U$ can be
uniquely written in form
\eqref{canonicalform}, where elements $t_r$-s are in the definition field
of~$G$. If $O^{p'}(G)$ coinsides with one of the groups
${}^2A_n(q)$, $^2B_2(2^{2n+1})$, ${}^2D_n(q)$, ${}^3D_4(q)$,
${}^2E_6(q)$, ${}^2G_2(3^{2n+1})$, or ${}^2F_4(2^{2n+1})$, then we shall say
that $G$ is {\em twisted}\index{group!of Lie type!finite twisted}, in the
remaining cases $G$ is called {\em split}\index{group!of Lie type!finite split}.
If $O^{p'}(\ov{G}_\sigma)\leq G\leq \ov{G}_\sigma$ is a twisted group of Lie
type and  $r\in\Phi(\ov{G})$, then by
$\bar{r}$\glossary{rbar@$\bar{r}$} we always denote the image of a root $r$
under the symmetry of the root system, corresponding to the graph automorphism
used during the construction of $G$. Sometimes we shall use the notation
$\Phi^\varepsilon(q)$\glossary{Phiepsilon@$\Phi^\varepsilon(q)$},
where  $\varepsilon\in\{+,-\}$, and  $\Phi^+(q)=\Phi(q)$ is a split group of
Lie type with the base field  $\F_q$,
$\Phi^-(q)={}^2\Phi(q)$ is a twisted group of Lie type defined over the field
$\F_{q^2}$ (with the base field~$\F_q$).

Let $\ov{R}$ be a connected $\sigma$-stab\-le subgroup of
$\ov{G}$. Then we may consider $R=G\cap \ov{R}$ and
$N(G,R)=G\cap N_{\ov{G}}(\ov{R})$\glossary{NGR@$N(G,R)$}. Note that
$N(G,R)\not=N_G(R)$ in general, and $N(G,R)$ is called the  {\em algebraic
normalizer}\index{algebraic normalizer $N(G,R)$} of $R$. For example, if we
consider $G=\SL_n(2)$, then the subgroup of diagonal matrices
$H$ of $G$ is trivial, hence $N_G(H)=G$. But
$G=(\SL_n(\ov{\F}_2))_\sigma$, where $\sigma$ is a Frobenius map
$\sigma:(a_{i,j})\mapsto (a_{i,j}^2)$. Then $H=\ov{H}_\sigma$, where
$\ov{H}$ is the subgroup of diagonal matrices in $\SL_n(\ov{\F}_2)$. Thus
$N(G,H)$ is the group of monomial matrices in $G$. We use the term
``algebraic normalizer'' in order to avoid such difficulties and
to make our proofs to be universal.  A group $R$ is called a {\em
torus}\index{torus!of finite group of Lie type} (respectively a
{\em reductive subgroup}\index{subgroup!reductive!of finite group of Lie type},
a {\em parabolic subgroup}\index{subgroup!parabolic!of finite group of
Lie type}, a {\em maximal torus}\index{torus!maximal!of finite group of Lie
type}, a {\em reductive subgroup of maximal
rank}\index{subgroup!reductive of maximal rank!of finite group of Lie type})
if $\ov{R}$ is a torus (respectively a reductive subgroup, a parabolic
subgroup, a maximal torus, a reductive subgroup of maximal rank) of $\ov{G}$.
A maximal $\sigma$-stab\-le torus $\ov{T}$ of $\ov{G}$ such that
$\ov{T}_\sigma$ is a Cartan subgroup of $\ov{G}_\sigma$ is called a {\em
maximal split torus}\index{torus!maximal!split} of~$\ov{G}$.

Assume that a reductive subgroup $\ov{R}$ is
$\sigma$-stab\-le. In view of
\cite[10.10]{SteLangSteinbergTheorem} there exists a
$\sigma$-stab\-le maximal torus $\ov{T}$ of $\ov{R}$. Let
$\ov{G}_{i_1},\ldots,\ov{G}_{i_{j_i}}$ be a $\sigma$-or\-bit of
$\ov{G}_{i_1}$. Consider the induced action of $\sigma$ on the factor group
\begin{equation*}
(\ov{G}_{i_1}\circ\ldots\circ\ov{G}_{i_{j_i}})/Z(\ov{G}_{i_1}
\circ\ldots\circ \ov{G}_{i_{j_i}})\simeq \P
\ov{G}_{i_1}\times\ldots\times\P \ov{G}_{i_{j_i}}.
\end{equation*}
Since $\P \ov{G}_{i_1}\simeq\ldots\simeq \P \ov{G}_{i_{j_i}}$ are simple
(as abstract groups), then $\sigma$ induces a cyclic permutation of the set
$\{\P \ov{G}_{i_1},\ldots, \P \ov{G}_{i_{j_i}}\}$, and we may assume that the
numberring is chosen so that $\P \ov{G}_{i_1}^\sigma=\P \ov{G}_{i_2}$, \ldots,
$\P \ov{G}_{i_{j_i}}^\sigma=\P \ov{G}_{i_1}$. Thus the equality
\begin{multline*}(\P \ov{G}_{i_1}\times\ldots\times\P
\ov{G}_{i_{j_i}})_\sigma=\\
\{x\mid x= g\cdot g^\sigma\cdot\ldots\cdot g^{\sigma^{j_i-1}}\text{ for some
}g\in\P \ov{G}_{i_1}\}_\sigma\simeq (\P \ov{G}_{i_1})_{\sigma^{j_i}}
\end{multline*}
holds. In view of
\cite[10.15]{SteLangSteinbergTheorem} the group $\P \ov{G}_{\sigma^{j_i}}$ is
finite, hence $O^{p'}((\P \ov{G}_{i_1})_{\sigma^{j_i}})$ is a canonical  finite
group of Lie type, probably with the base field larger than the base field
of~$O^{p'}(\ov{G}_\sigma)$.

Let $\ov{B}_{i_1}$ be the preimage of a $\sigma^{j_i}$-stab\-le Borel subgroup
of $\P \ov{G}_{i_1}$ in $\ov{G}_{i_1}$ under the natural epimorphism, and
$\ov{T}_{i_1}$ be a $\sigma^{j_i}$-stab\-le maximal torus of
$\ov{G}_{i_1}$, contained in $\ov{B}_{i_1}$ (their existence follows from
\cite[10.10]{SteLangSteinbergTheorem}). Then from the note at the beginning of
section 11 from \cite{SteLangSteinbergTheorem}, subgroups
$\ov{U}_{i_1}$ and $\ov{U}_{i_1}^-$, generated by
$\ov{T}_{i_1}$-in\-va\-ri\-ant root subgroups, taken over all positive and
negative roots respectively, are also $\sigma^{j_i}$-stab\-le. Since
$\ov{G}_{i_1}$ is a simple algebraic group, then $\ov{G}_{i_1}$ is generated by
subgroups $\ov{U}_{i_1}$ and $\ov{U}_{i_1}^-$. Now
$Z(\ov{G}_{i_1}\circ\ldots\circ \ov{G}_{i_{j_i}})$ consists of semisimple
elements, so the restriction of the natural epimorphism
$\ov{G}_{i_1}\rightarrow \P \ov{G}_{i_1}$ on $\ov{U}_{i_1}$ and
$\ov{U}_{i_1}^-$ is an isomorphism. Therefore, for each $k$ the subgroups
$(\ov{U}_{i_1})^{\sigma^k}$ and  $(\ov{U}_{i_1}^-)^{\sigma^k}$ are maximal
$\sigma^{j_i}$-stab-\-le connected unipotent subgroups of $\ov{G}_{i_k}$ and
they generate~$\ov{G}_{i_k}$.

Thus, $\ov{U}_{i_1}\times
(\ov{U}_{i_1})^\sigma\times\ldots\times(\ov{U}_{i_1})^{\sigma^{j_i-1}}$
and $\ov{U}_{i_1}^-\times
(\ov{U}_{i_1}^-)^\sigma\times\ldots\times(\ov{U}_{i_1}^-)^{\sigma^{j_i-1}}$
are maximal $\sigma$-stab\-le connected unipotent  subgroups of
$\ov{G}_{i_1}\circ\ldots\circ \ov{G}_{i_{j_i}}$ and they generate
$\ov{G}_{i_1}\circ\ldots\circ \ov{G}_{i_{j_i}}$. By
\cite[Corollary~12.3(a)]{SteLangSteinbergTheorem}, we have
\begin{multline*}
O^{p'}((\ov{G}_{i_1}\circ\ldots\circ \ov{G}_{i_{j_i}})_\sigma)=\\
\langle (\ov{U}_{i_1}\times
(\ov{U}_{i_1})^\sigma\times\ldots\times(\ov{U}_{i_1})^{\sigma^{j_i-1}})_\sigma,
(\ov{U}_{i_1}^-\times (\ov{U}_{i_1}^-)^\sigma\times\ldots\times
(\ov{U}_{i_1}^-)^{\sigma^{j_i-1}})_\sigma\rangle \simeq\\ \langle
(\ov{U}_{i_1})_{\sigma^{j_i}},(\ov{U}_{i_1}^-)_{\sigma^{j_i}}\rangle=
O^{p'}((\ov{G}_{i_1})_{\sigma^{j_i}}).
\end{multline*}
By \cite[11.6 and Corollary~12.3]{SteLangSteinbergTheorem}, the
group $\langle
(\ov{U}_{i_1})_{\sigma^{j_i}},(\ov{U}_{i_1}^-)_{\sigma^{j_i}}\rangle$
is a canonical finite group of Lie type. More over, from the above
arguments it follows that the groups  $ \langle
(\ov{U}_{i_1})_{\sigma^{j_i}},(\ov{U}_{i_1}^-)_{\sigma^{j_i}}\rangle/
Z( \langle
(\ov{U}_{i_1})_{\sigma^{j_i}},(\ov{U}_{i_1}^-)_{\sigma^{j_i}}\rangle)$
and $O^{p'}((\P\ov{G}_{i_1})_{\sigma^{j_i}})$ are isomorphic.
Denoting $O^{p'}((\ov{G}_{i_1}\circ\ldots\circ
\ov{G}_{i_{j_i}})_\sigma)$ by $G_i$, we obtain that $G_i$ is a
canonical finite group of Lie type for all~$i$. Subgroups $G_i$-s of
$O^{p'}(\ov{G}_\sigma)$, appearing in this way, are called {\em
subsystem subgroups}\index{subgroup!subsystem}
of~$O^{p'}(\ov{G}_\sigma)$.

Since
$\ov{G}_{i_1}\circ\ldots\circ\ov{G}_{i_{j_i}}$ is a
$\sigma$-stab\-le subgroup, then
$\ov{G}_{i_1}\circ\ldots\circ\ov{G}_{i_{j_i}}\cap \ov{T}$ is a
$\sigma$-stab\-le maximal torus of
$\ov{G}_{i_1}\circ\ldots\circ\ov{G}_{i_{j_i}}$. Therefore we may assume that
for each $\sigma$-or\-bit
$\{\ov{G}_{i_1},\ldots,\ov{G}_{i_{j_i}}\}$, the intersection
$\ov{T}\cap \ov{G}_{i_1}\circ\ldots\circ\ov{G}_{i_{j_i}}$ is a maximal
$\sigma$-stab\-le torus of
$\ov{G}_{i_1}\circ\ldots\circ\ov{G}_{i_{j_i}}$. Then
$\ov{R}_\sigma=\ov{T}_\sigma(G_1\circ\ldots\circ G_m)$ and
$\ov{T}_\sigma$ normalizes each of~{$G_i$-s}.

For a  $\sigma$-or\-bit $\{\ov{G}_{i_1},\ldots,\ov{G}_{i_{j_i}}\}$ of
$\ov{G}_{i_1}$, where $G_i=O^{p'}((\ov{G}_{i_1}\circ\ldots\circ
\ov{G}_{i_{j_i}})_\sigma)$, consider $\Aut_{\ov{R}_\sigma}(G_i)$. Since
$G_1\circ\ldots\circ G_{i-1}\circ G_{i+1}\circ\ldots\circ
G_k\circ\ov{Z}_\sigma\leq C_{\ov{R}_\sigma}(G_i)$, we have that
$\Aut_{\ov{R}_\sigma}(G_i)\simeq \left(\ov{T}_\sigma G_i\right)/
Z\left(\ov{T}_\sigma G_i\right).$ From \cite[Proposition~2.6.2]{GorLySol}
it follows that automorphisms induced by $\ov{T}_\sigma$ on
$G_i$, are diagonal. Therefore, the inclusions $\P G_i\leq
\Aut_{\ov{R}_\sigma}(G_i)\leq \widehat{\P G_i}$ hold, in particular,
$\Aut_{\ov{R}_\sigma}(G_i)$ is a finite group of Lie type.

Now consider the case, when $\ov{L}\unlhd \ov{H}\leq \ov{G}$, where $\ov{L}$
and $\ov{H}$ are  $\sigma$-stab\-le and closed. Clearly $\sigma$ induces an
action on $\ov{H}/\ov{L}$ and, if $\ov{L}$ is connected, then
Lang-Ste\-i\-n\-berg Theorem (Lemma \ref{LangSteinbergTheorem}) implies
$(\ov{H}/\ov{L})_\sigma=\ov{H}_\sigma/\ov{L}_\sigma$. Let $\ov{R}$ be a
$\sigma$-stab\-le connected reductive subgroup of maximal rank (in particular,
$\ov{R}$ can be a maximal torus) of $G$. Since groups
$N_{\ov{G}}(\ov{R})/\ov{R}$ and $N_W(W_{\ov{R}})/W_{\ov{R}}$ are isomorphic,
where $W$ is the Weyl group of $\ov{G}$, $W_{\ov{R}}$ is the Weyl group of
$\ov{R}$ (and it is a subgroup of $W$), we obtain an induced action of
$\sigma$ on $N_W(W_{\ov{R}})/W_{\ov{R}}$, and we say that
$w_1\equiv w_2$, for
$w_1,w_2\in N_W(W_{\ov{R}})/W_{\ov{R}}$, if there exists an element  $w\in
N_W(W_{\ov{R}})/W_{\ov{R}}$, satisfying to the equality
$w_1=w^{-1}w_2w^\sigma$. Let
$Cl(\ov{G}_\sigma,\ov{R})$\glossary{CLGsigmaR@$Cl(\ov{G}_\sigma,\ov{R})$} be
the set of  $\ov{G}_\sigma$-con\-ju\-ga\-ted classes of
$\sigma$-stab\-le subgroups $\ov{R}^g$, where $g\in\ov{G}$. Then
$Cl(\ov{G}_\sigma, \ov{R})$ is in 1-1 correspondence with the set of
$\sigma$-con\-ju\-ga\-te classes
$Cl(N_W(W_{\ov{R}})/W_{\ov{R}},\sigma)$. If $w$ is an element of
$N_W(W_{\ov{R}})/W_{\ov{R}}$, and $(\ov{R}^g)_\sigma$ corresponds to the
$\sigma$-con\-ju\-ga\-te class of $w$, then $(\ov{R}^g)_\sigma$ is said to be
obtained by {\em twisting} of the group  $\ov{R}$ by the element
$w\sigma$. Further $(\ov{R}^g)_\sigma\simeq \ov{R}_{\sigma w}$. The
construction of twisting is known and is given, for example, in
\cite{Car5CentSemisimpleLie} with all necessary results. When $\ov{H}=\ov{T}$
is a $\sigma$-stab\-le maximal torus and
$W=N_G(T)/T$, then by
\cite[Proposition~3.3.6]{CarFiniteGrpsLieTypeConjClassesCharacters},
\begin{equation}\label{twisted}
\left(\frac {N_G(T_w)}
{T_w}\right)_\sigma=\frac {(N_G(T_w))_\sigma}
{(T_w)_\sigma }\simeq C_{W,\sigma}(w)=\{x\in W\mid \sigma(x)
wx^{-1}= w\}.
\end{equation}

Now assume that the group $\ov{R}$ is a $\sigma$-stab\-le parabolic subgroup of
$\ov{G}$ and $\ov{U}$ is its unipotent radical. Then it contains a connected
reductive subgroup $\ov{L}$ such that $\ov{R}/\ov{U}\simeq \ov{L}$. A subgroup
$\ov{L}$ is called a {\em Levi factor}\index{Levi factor} of $\ov{R}$. More
over, if $\ov{Z}=Z(\ov{L})^0$, then $\ov{L}=C_{\ov{G}}(\ov{Z})$ (see
\cite[30.2]{HuLinearAlgGrps}). Let $R(\ov{R})$ be the radical of $\ov{R}$. Then
it is a $\sigma$-stab\-le connected solvable subgroup, hence by
\cite[10.10]{SteLangSteinbergTheorem}, it contains a $\sigma$-stab\-le maximal
torus $\ov{Z}$. Now $C_{\ov{G}}(\ov{Z})=C_{\ov{R}}(\ov{Z})$ is a
$\sigma$-stab\-le Levi factor of $\ov{R}$. Thus each
$\sigma$-stab\-le parabolic subgroup of $\ov{G}$ contains a
$\sigma$-stab\-le Levi factor $\ov{L}$ and $\ov{L}$ is a connected reductive
subgroup of maximal rank of~$\ov{G}$.

\subsection{Known results}

In this section we recall some structure results that will be often used below.

\begin{lem}\label{CentrOfSemisimpleClassic} {\em
\cite[Theorem~2.2]{Hu2ConjClasses}}
Let $\ov{G}$ be a connected reductive algebraic group, $s\in \ov{G}$ be a
semisimple element of $\ov{G}$ and  $\ov{T}$ be a maximal torus of
$\ov{G}$, containing~$s$.

\noindent Then  $C_{\ov{G}}(s)^0$ is a reductive subgroup of maximal rank of
$\ov{G}$. The centralizer  $C_G(s)$ is generated by a torus $\ov{T}$, those
$\ov{T}$-ro\-ot subgroups $X_r$, for which $s^r=e$ and representatives
$n_w$ of elements $w\in W$, which commute with $s$. Further
$C_{\ov{G}}(s)^0$ is generated by the torus $\ov{T}$, those
$\ov{T}$-ro\-ot subgroups $X_r$, for which $s^r=e$, and each unipotent element
centralizing $s$, is in~${C_{\ov{G}}(s)^0}$.
\end{lem}

\begin{lem}\label{FactorCentrByConnectedComponent} {\em
\cite[Proposition~2.10]{Hu2ConjClasses}}
Let $\ov{G}$ be a simple algebraic group and $s$ be its semisimple element
of finite order.

\noindent Then the factor group $C_{\ov{G}}(s)/C_{\ov{G}}(s)^0$ is isomorphic
to a subgroup of the fundamental group~$\Delta(\ov{G})$. In particular, if
$\ov{G}$ is simply connected, then $C_{\ov{G}}(s)$ is connected.
\end{lem}

\begin{lem}\label{LangSteinbergTheorem}{\em
\cite[Theorem~10.1]{SteLangSteinbergTheorem}} \index{theorem!Lang-Steinberg}
Let $\ov{G}$ be a connected algebraic group and $\sigma$ be a Frobenius map.

\noindent Then the map $x\mapsto x^{-1}x^\sigma$ is surjective.
\end{lem}

The following lemma is known as Bo\-rel-Tits theorem.

\begin{lem}\label{Borel-Titsclassic}\index{theorem!Borel-Tits}
Let $X$ be a subgroup of a finite group of Lie type $G$ such that $O_p(X)$
is nontrivial.

\noindent Then there exists a $\sigma$-stab\-le parabolic subgroup
$\overline{P}$ of $\overline{G}$ such that $X\leq \overline{P}$ and $O_p(X)\leq
R_u(\overline{P})$.
\end{lem}

\begin{proof}
Define $U_0=O_p(X)$, $N_0=N_{\overline{G}}(U_0)$. Then
$U_i=U\cdot R_u(N_{i-1})$ and $N_i=N_{\overline{G}}(U_i)$. Clearly $U_i$, $N_i$
are $\sigma$-stab\-le for all $i$. In view of
\cite[Proposition~30.3]{HuLinearAlgGrps}, the chain of subgroups
$N_0\leq N_1\leq \ldots\leq N_k\leq\ldots$ is finite and $\overline{P}=\cup_i
N_i$ is a proper parabolic subgroup. Clearly $\overline{P}$ is
$\sigma$-stab\-le.
\end{proof}

\begin{lem}\label{Hartley-Shute}
{\em (Hartley-Shute lemma\index{лемма Хартли-Шюта}
\cite[Lemma~2.2]{HartShute})} Let $G$ be a finite canonical adjoint group of
Lie type with the definition field  $\F_q$. Let $H$ be a Cartan subgroup of
$G$ and $s\in \F_q$. If $G$ is twisted and  $r=\bar{r}$, then assume also that
$s$ is in the base field of~$G$.

\noindent Then there exists an element $h(\chi)\in H$ such that
$\chi(r)=s$, except the following cases, when  $h(\chi)$ can be chosen so that
$\chi(r)$ would have indicated values:
\begin{itemize}
\item[{\em (a)}] $G=A_1(q)$, $\chi(r)=s^2$;
\item[{\em (b)}] $G=C_n(q)$, $r$ is a long root,
$\chi(r)=s^2$;
\item[{\em (c)}] $G={}^2A_2(q)$, $r\not=\bar{r}$, $\chi(r)=s^3$;
\item[{\em (d)}] $G={}^2 A_3(q)$, $r\not=\bar{r}$, $\chi(r)=s^2$;
\item[{\em (e)}] $G={}^2D_n(q)$, $r\not=\bar{r}$, $\chi(r)=s^2$;
\item[{\em (f)}] $G={}^2G_2(3^{2n+1})$, $r=a$ or $r=3a+b$, where $a$ is a
short,  $b$ is a long fundamental roots, $\chi(r)=s^2$.
\end{itemize}
\end{lem}

\begin{ttt}\label{CarterInClassicalGroups} {\em
\cite[Theorem~1.1]{PreTamVdoCartClassical}}
Let $q=p^\alpha$, where $p$ is a prime, and assume $G=\Sp_n(q)$, or
$\mathrm{SO}^\varepsilon_n(q)\leq G\leq \GO^\varepsilon_n(q)$, where $q$ is
odd, or $\SU_n(q)\leq G\leq\GU_n(q)$. If $G$ admits a a Carter
subgroup $K$, then either $K$ is the normalizer of a Sylow $2$-sub\-gro\-up
of $G$, or one of the following  holds:
\begin{itemize}
\item[{\em (a)}] $G\in \{\Sp_2(3), \ \SL_2(3),\ 2. \SU_2(3)\}$ and $K$ is the
normalizer of a Sylow $3$-sub\-gro\-up of~$G$;
\item[{\em (b)}] $G=\GU_3(2)$ has  order $2^3\cdot 3^4$, and $K$ has
order~${2\cdot 3^2}$.
\end{itemize}
Moreover, if  $G$ is orthogonal, $K$ is a $2$-group, except possibly
when~${G=\mathrm{SO}_2^\varepsilon(q)}$.
\end{ttt}

\section{Conjugacy criterion for Carter subgroups}

\subsection{Main results of this section}

\begin{df}
A finite group $G$ is said to satisfy condition
{\bfseries (C)}\index{Conjugacy condition {\bfseries
(C)}}\glossary{C@{\bfseries(C)}}, if, for {\slshape every} non-Abelian
composition factor $S$ of
{\slshape every} composition series of $G$
and for every its nilpotent subgroup $N$, Carter subgroups of
$\langle \mathrm{Aut}_N(S),S\rangle$ are conjugate (in particular, they may not
exist).
\end{df}

\begin{lem}\label{CarterSubgroupsInCompFactorsUnderHomomorphism}
Let  $H$ be a normal subgroup of a finite group $G$, $B\lhd A\leq G$ and
$S=(A/H)/(B/H)$ is a composition factor of $G/H$, and  $L\leq G$.

\noindent Then $\Aut_L(A/B)\simeq \Aut_{LH/H}((A/H)/(B/H))$.
\end{lem}

\begin{proof}
Since $H\leq B$, then $H\leq C_G(A/B)$, so we may assume that $L=LH$. Further
more we may assume that $L\leq N_G(A)\cap N_G(B)$ and
$G=LA$. Then the action on $A/B$ given by
$x:Ba\mapsto Bx^{-1}ax$ coincides with the action on $(A/H)/(B/H)$ given by
$xH:BaH\mapsto Bx^{-1}axH$, and the lemma follows.
\end{proof}

The following lemma is known.

\begin{lem}\label{nilpotentExists} Let $G$ be a finite group, $H$ be a normal
subgroup of $G$ and $\overline{N}$ be a nilpotent subgroup of
$\overline{G}=G/H$.

\noindent Then there exists a nilpotent subgroup $N$ of
$G$ such that~${NH/H=\overline{N}}$.
\end{lem}

\begin{proof}
Clearly we may assume that $G/H=\overline{N}$. There exists a subgroup
$U$ of $G$ such that $UH=G$. Choose a subgroup of minimal order with this
property. Then $U\cap H$ is contained in the Frattini subgroup $F$ of $U$.
Indeed, if there exists a maximal subgroup $M$ of $U$, not containing
$U\cap H$, then clearly $MH=G$, which contradicts the minimality of
$U$. Thus the group $U/F$ is nilpotent, hence $U$ is nilpotent and~$N=U$.
\end{proof}

By Lemmas \ref{CarterSubgroupsInCompFactorsUnderHomomorphism} and
\ref{nilpotentExists} it follows that, if a finite group $G$ satisfies
{\bfseries(C)}, then for every its normal subgroup $N$ and solvable subgroup
$H$, groups $G/N$ and $HN$ satisfy~{\bfseries(C)}.

In this section we prove that if $G$ satisfies {\bfseries
(C)}, then its Carter subgroups are conjugate. More precisely, the following
theorem will be proven.

\begin{ttt}\label{ConjugacyCriterion}\index{Conjugacy criterion}
If a finite group~$G$ satisfies {\em\bfseries(C)}, then Carter subgroups of~$G$
are conjugate.
\end{ttt}

Below in subsections 2.2,~2.3 we are assuming that $X$ is a counter example to
Theorem \ref{ConjugacyCriterion} of minimal order, i.~e., that $X$ is a finite
group satisfying~{\bfseries(C)}, and  $X$ contains nonconjugate Carter
subgroups, but Carter subgroups in every group~$M$  of order less than~$\vert
X\vert$, satisfying~{\bfseries (C)}, are conjugate.

\subsection{Preliminary results}

Recall that $X$ is a counter example to Theorem \ref{ConjugacyCriterion} of
minimal order.

\begin{lem}\label{FactorCarterX}
Let $G$ be a finite group satisfying {\em\bfseries(C)},
 $\vert G\vert\leqslant \vert X \vert$, and $H$ be a Carter subgroup of~$G$.

\noindent If $N$ is a normal subgroup of~$G$, then $HN/N$ is a Carter subgroup
of~$G/N$.
\end{lem}

\begin{proof}
Since $HN/N$ is nilpotent, we have just to prove that it is
self-nor\-ma\-li\-zing in
$G/N$. Clearly, this is true if $G=HN$. So, assume $M=HN<G$ (note that by
Lemmas  \ref{CarterSubgroupsInCompFactorsUnderHomomorphism} and
\ref{nilpotentExists} the group  $M$ satisfies {\bfseries(C)}). By the
minimality
of $X$, $M^x=M$, $x\in G$, implies $H^x=H^m$ for some $m\in M$. It follows
$xm^{-1}\in N_G(H)=H$ and $x\in M$. This proves that $HN/N$ is nilpotent and
self-normalizing in~$G/N$.
\end{proof}

\begin{lem}\label{PropertiesOfX}
Let $B$ be a minimal normal subgroup of $X$ and $H,K$ be non-conjugate Carter
subgroups of~$X$.
\begin{itemize}
\item[{\rm (1)}] $B$ is non-soluble.
\item[{\rm (2)}]  $X=BH=BK$.
\item[{\rm (3)}] $B$ is the unique minimal normal subgroup of~$X$.
\end{itemize}
\end{lem}

\begin{proof}
(1) We give a proof by contradiction. Assume that $B$ is soluble and let
$\pi:X\rightarrow X/B$ be the canonical homomorphism. Then $H^\pi$ and $K^\pi$
are Carter subgroups of $X/B$, by Lemma \ref{FactorCarterX}. By the
minimality of $X$, there  exists $\bar{x}=Bx$ such that
$(K^\pi)^{\bar{x}}=H^\pi$. It follows $K^x\leq BH$. Since $BH$ is soluble,
$K^x$ is conjugate to $H$ in $BH$, hence $K$ is conjugate to $H$ in $X$, a
contradiction.

(2) Assume that $BH<X$. By Lemma~\ref{FactorCarterX} and the minimality of
$X$, $BH/B$ and $BK/B$ are conjugate in $X/B$: so there exists $x\in X$ such
that $K^x\leq BH$. It follows that $K^x$ is conjugate to $H$ in $BH$, hence $K$
is conjugate to $H$ in $X$, a contradiction.

(3) Suppose that $M$ is a minimal normal subgroup of $X$ different from $B$.
By~(1), $M$ is non-soluble. On the other hand, $MB/B\simeq M$ is a subgroup of
the nilpotent group $X/B\simeq H/H\cap B$, a contradiction.
\end{proof}

\begin{lem}\label{CarterInGroupOFInducedAutomorphismsWithoutCenter}
Let $K$ be a Carter subgroup of a finite group $G$. Assume that there exists a
normal subgroup $B=T_1\times\ldots\times T_k$ of $G$ such that $G=KB$,
$Z(T_i)=\{e\}$, and $T_i$ is not decomposable into direct product of its proper
subgroups for all~$i$.

\noindent Then $\mathrm{Aut}_K(T_i)$ is a  Carter subgroup of~$\langle
\mathrm{Aut}_K(T_i), T_i\rangle$.
\end{lem}

\begin{proof}
Assume that our statement is false and $G$ is a counterexample with $k$
minimal, then $k>1$. Since each group $T_i$ has trivial center and is not
decomposable into direct product of proper subgroups, a corollary of
Krull-Remak-Shmidt theorem \cite[3.3.10]{Rob} implies that the action by
conjugation of $G$ on the set $\{T_1,\ldots,T_k\}$ induces permutations of this
set. Clearly, $G$ acts transitively, by conjugation, on the set
$\Omega:=\{T_1,\ldots, T_k\}$. We may
assume that the $T_j$-s are indexed so that $G$ acts
primitively on the set $\{\Delta_1,\ldots,\Delta_p\}$, $p>1$, where for
each~$i$: $$\Delta_i:=\{T_{1+(i-1)l},\ldots,T_{il}\},\ \ \ \ \ k=pl.$$
Denote by $\varphi:G\rightarrow \mathrm{Sym}_p$  the induced permutation
representation. Clearly, $B\leq \mathrm{ker}\ \varphi$, so that
$G^\varphi=(BK)^\varphi=K^\varphi$ is a primitive nilpotent subgroup of
$\mathrm{Sym}_p$. Hence $p$ is prime and $G^\varphi$ is a cyclic group of order
$p$. In
particular, $Y:=ker\ \varphi$ coincides with the stabilizer of any $\Delta_i$,
so that $\varphi$ is permutationally equivalent to the representation of $G$ on
the right cosets of $Y$.  For each $i=1,\ldots,p$, let
$S_i=T_{1+(i-1)l}\times\ldots\times T_{il}$. Then $Y=N_G(S_i)$ and
$B=S_1\times\ldots\times S_p$. Consider $\xi:Y\rightarrow \mathrm{Aut}_Y(S_1)$,
let
$A=Y^\xi$, $S=S_1^\xi$. Clearly $S$ is a normal subgroup of $A$; moreover, $S$
is isomorphic to $S_1$, since $S_1$ has trivial center. On the other hand, for
each $i\not=1$, $S_i\leq \mathrm{ker}\ \xi$, since $S_i$ centralizes~$S_1$.

Denote by $A\wr \Z_p$ the wreath product of~$A$ and a cyclic group
$\Z_p$, and let
$\{x_1=e,\ldots,x_p\}$ be a right transversal of $Y$. Then the map
$\eta:G\rightarrow A\wr \Z_p$ such that, for each $x\in G$:
$$x\mapsto \left(\left(x_1xx^{-1}_{1^{x^\varphi}}
\right)^\xi,\ldots,\left(x_pxx^{-1}_{p^{x^\varphi}}
\right)^\xi\right)x^\varphi$$ is a homomorphism. Clearly $Y^\eta$ is a
subdirect product of the base subgroup $A^p$ and $$S_1^\eta=
\{(s,1,\ldots,1)\vert s\in S\}, B^\eta=\{(s_1,\ldots,s_p)\vert s_i\in S\}\leq
Y^\eta.$$ Moreover, $\mathrm{ker}\ \eta=C_G(B)=\{e\}$, so we may identify $G$
with
$G^\eta$. We choose $h\in K\setminus Y$. Then $$G=\langle Y,h\rangle,\ h^p\in
Y,\ K=(Y\cap K)\langle h\rangle$$ and we may assume
$$ h=(a_1,a_2,\dots,a_p)\pi, \ a_i\in A,
\quad \pi=(1,2,\dots,p)\in \Z_p. $$ For each
$i$, $1\leqslant i\leqslant p$, let $\psi_i:A^p\rightarrow A$ be the canonical
projection
and let $K_i:=(K\cap Y)^{\psi_i}$. Clearly, $Y^{\psi_i}=A$. Moreover, for each
$i\geqslant 2$, $K_i=K_1^{h^{i-1}}=K_1^{a_1\ldots a_{i-1}}$ since $h$
normalizes~${Y\cap K}$. Let $N:=(K_1\times\ldots\times K_p)\cap Y.$ $N$ is
normalized by $K$, since $K=(N\cap K)\langle h\rangle$ and $K_i^h=K_{i+1\pmod
p}.$ We claim that  $K_1$ is a Carter subgroup of $A$. Assume $n_1\in
N_A(K_1)\setminus K_1$. From $Y=(Y\cap K)B$, it follows $n_1=h_1s$, $h_1\in
K_1$, $s\in N_S(K_1)\setminus K_1$. Let $b:=(s,s^{a_1},\ldots,s^{a_1\ldots
a_{p-1}})\in B.$ Then $b$ normalizes $N$, for:
$$K_i^b=K_i^{s^{a_1\ldots
a_{i-1}}}=K_1^{a_1\ldots a_{i-1}s^{a_1\ldots
a_{i-1}}}=K_1^{sa_1\ldots a_{i-1}}=K_1^{a_1\ldots a_{i-1}} =K_i.$$ Now
$[b,h^{-1}]:=b^{-1}hbh^{-1}\in Y$ is such that: $$[b,h^{-1}]^{\psi_i}=1\text{
if }i\not=p, [b,h^{-1}]^{\psi_p}=[s,(a_1\cdot\ldots\cdot
a_p)^{-1}]^{a_1\cdot\ldots\cdot
a_{p-1}},$$ where $a_1\cdot\ldots\cdot a_p=(h^p)^{\psi_1}\in K_1$. Since $s\in
N_S(K_1)$,
it follows $$[s,(a_1\cdot\ldots\cdot a_p)^{-1}]\in K_1,\
[s,(a_1\cdot\ldots\cdot a_p)^{-1}]^{a_1\cdot\ldots\cdot a_{p-1}}\in K_p.$$ So
$[b,h^{-1}]\in N$
and $b\in  N_G(N\langle h\rangle)$. But $K\leq N\langle h\rangle$, implies
$N_G(N\langle h\rangle)=N\langle h\rangle$. Indeed, if $g\in N_G(N\langle
h\rangle)$, then $K^g$ is a Carter subgroup of $N\langle h\rangle$. But
$N\langle h\rangle$ is soluble, hence there exists $y\in N\langle h\rangle$
with $K^g=K^y$. Now $K$ is a Carter subgroup of $G$, thus $gy^{-1}\in K$
and $g\in N\langle h\rangle$. Therefore $b\in N, s\in K_1$, i.~e., $n_1\in K_1$,
a contradiction.

Now $A=K_1(T_1\times\ldots\times T_l)$ and $l<k$. By induction we have that
$\mathrm{Aut}_{K_1}(T_1)$ is a Carter subgroup of
$\langle\mathrm{Aut}_{K_1}(T_1),T_1\rangle$.
In view of our construction, $\mathrm{Aut}_K(T_1)=\mathrm{Aut}_{K_1}(T_1)$ and
the lemma follows.
\end{proof}

\subsection{Proof of Theorem \ref{ConjugacyCriterion}}

Recall that $B=T_1\times\dots\times T_k$, where $T_i\simeq T$ is a
non-Abelian simple group. What remains to prove is $k=1$.
Осталось доказать, что $k=1$.  In the notations of the proof of
Lemma~\ref{CarterInGroupOFInducedAutomorphismsWithoutCenter} we
have shown that $H_1$ is a Carter subgroup of~$A$. If $k>1$, then
$\vert A\vert<\vert X\vert$ and $A$ satisfies {\bfseries(C)}. So each
$K_i$ is conjugate with~$K_1$ in~$A$ and $N_A(K_i)=K_i$,
$i=1,\dots,p$. It follows easily that $N$ is a Carter subgroup of~$Y$.
Let $y:=(y_1,\dots,y_p)\in N_Y(N)$. From $N^{\psi_i}=K_i$ we have
$y_i\in N_A(K_i)=K_i$ for each~$i$, hence $y\in N$.

We have seen that, to each Carter subgroup $K$ of $X$ we can associate a Carter
subgroup $N=N_K$ of $Y$, such that $K$ normalizes $N_K$. Clearly,
$N_K\not=\{e\}$, otherwise $X$ would have order $p$. So let $H$ be a Carter
subgroup of $X$, not conjugate to $K$, and let $N_H$ be the Carter subgroup of
$Y$ corresponding to $H$. If $k>1$, then $Y$ is a proper subgroup of $X$ and $Y$
satisfies {\bfseries (C)}. By the minimality of~$X$ we obtain that $N_H$ and
$N_K$ are conjugate in $Y$, and we may assume that $N_H=N_K$. Then $HN_H=KN_H$
is solvable, hence, the subgroups $H$ and $K$ are conjugate. This contradiction
completes the proof of Theorem~\ref{ConjugacyCriterion}.

\subsection{Some properties of Carter subgroups}

Here we shall prove some lemmas that will be useful in studying Carter
subgroups in finite groups, in particular, in almost simple groups.

\begin{lem}\label{HomImageOfCarter}
Let $K$ be a Carter and $N$ be a normal subgroups of a finite group $G$. Assume
that $KN$ satisfies {\em{\bfseries (C)}} (this condition holds if either
$G$ satisfies {\em\bfseries (C)} or $N$ is solvable) or~${KN=G}$.

\noindent Then $KN/N$ is a Carter subgroup of~$G/N$.
\end{lem}

\begin{proof}
If $KN=G$, then the statement is evident. Assume that $KN\not=G$, i.~e.,
$KN$ satisfies {\bfseries (C)}. Consider $x\in G$ and assume that $xN\leq
N_{G/N}(KN/N)$. Therefore $x\in N_G(KN)$. We have that $K^x$ is a Carter
subgroup of $KN$. Since  $KN$ satisfies {\bfseries (C)}, we obtain that its
Carter subgroups are conjugate. Thus there exists  $y\in KN$ such that
$K^y=K^x$. Since $K$ is a Carter subgroup of $G$, it follows that $xy^{-1}\in
N_G(K)=K$ and~${x\in KN}$.
\end{proof}

\begin{lem}\label{power}
Let $K$ be a Carter subgroup of a finite group~$G$. Assume also that
$e\not= z\in Z(K)$ and $C_G(z)$  satisfies~{\em\bfseries(C)}.
\begin{itemize}
\item[{\em (1)}] Every subgroup $Y$ which contains $K$ and satisfies
{\em\bfseries(C)}, is self-normalizing in~$G$.
\item[{\em(2)}] No conjugate of $z$ in $G$, except $z$,
lies in $Z(G)$.
\item[{\em(3)}] If $H$ is a Carter subgroup of $G$, non-conjugate
to $K$, then $z$ is not conjugate to any
element in the center of $H$.
\end{itemize}
In particular the centralizer $C_G(z)$ is self-normalizing
in $G$, and $z$ is not conjugate to any power~${z^k\not= z}$.
\end{lem}

\begin{proof}
(1) Take $x\in N_G(Y)$. Then $K^x$ is a Carter subgroup of $Y$. By Theorem
\ref{ConjugacyCriterion} Carter subgroups of $Y$ are conjugate. Therefore there
exists $y\in Y$ with
$K^x=K^y$. Hence $$xy^{-1}\in N_G(K)=K\leq Y \text{ and }x\in Y.$$

(2) Assume $z^{x^{-1}}\in Z(K)$ for some $x\in G$. Then
$z$ belongs to the center of $\langle G,G^x\rangle\leq C_G(z)$. Since $C_G(z)$
satisfies {\bfseries(C)}, there exists $y\in C_G(z)$ such that $K^x=K^y$. From
$xy^{-1}\in
C_G(z)$, we get $z^{xy^{-1}}=z$ hence $z^x=z^y=z$. We conclude $z^{x^{-1}}=z$.

(3) If our claim is false, substituting $H$
with some conjugate $H^x$ (if necessary), we may assume
$z\in Z(K)\cap Z(H)$, i.~e. $z\in Z(\langle K,H\rangle)\leq C_G(z)$. Again since
$C_G(z)$ satisfies~{\bfseries(C)}, there exists $y\in C_G(z)$ such that
$H=K^y$; a contradiction.
\end{proof}

Note that for every known finite simple group~$G$
(and hence almost simple, since the group of outer automorphisms is soluble)
and for all elements $z\in  G$ of prime order we see that composition factors
of~$C_G(z)$ are known simple groups. Indeed, for sporadic groups this statement
can be checked by using~\cite{ATLAS}. Composition factors of
$C_{A_n}(z)$ are alternating groups. If $G$ is a finite simple group of Lie
type over a field of characteristic~$p$ and $(\vert z\vert,p)=1$, then $z$ is
semisimple and all composition factors of $C_G(z)$ are finite groups of
Lie type. If  $\vert z\vert=p$ and $p$ is a good prime for~$G$, then by
Theorems 1.2 and~1.4 from \cite{SeitzConjUnipElts}, all composition
factors of~$C_G(z)$ are finite groups of Lie type. From papers of several
authors it follows that in case when $p$ is a bad prime for a finite adjoint
group of Lie type~$G$, all composition factors of the centralizer of an
element of order $p$ are known finite simple groups. Therefore, if we are
classifying Carter subgroups of an almost simple group~$A$, then by induction
we may assume that $C_A(z)$ satisfies {\bfseries(C)} for all elements
$z\in A$ of prime order.

\begin{lem}\label{CritSyl2Carter}
Let $Q$ be a Sylow $2$-sub\-gro\-up of a finite group~$G$.

\noindent Then $G$ contains a Carter subgroup $K$, satisfying $Q\leq K$, if
and only if ${N_G(Q)=QC_G(Q)}$.
\end{lem}

\begin{proof}
Assume that $G$ contains a Carter subgroup $K$, satisfying $Q\leq
K$. Since $K$ is nilpotent, it follows that $Q$ is normal in $K$ and
$K\leq QC_G(Q)\unlhd N_G(Q)$. By Feit-Thomp\-son theorem (see
\cite{ftOddOrder}), we obtain that $N_G(Q)$ is solvable. Thus by Lemma
\ref{power}(1) we have that
$QC_G(Q)$ is self-nor\-ma\-li\-zing in $N_G(Q)$, so~${N_G(Q)=Q C_G(Q)}$.

Assume now that $N_G(Q)=QC_G(Q)$, i.~e., the equality
$N_G(Q)=Q\times O(C_G(Q))$ holds. Since $O(C_G(S))$ is of odd order, it is
solvable. Therefore it contains a Carter subgroup $K_1$. Consider a nilpotent
subgroup $K=Q\times K_1$ of $G$. Assume that $x\in N_G(K)$, then $x\in
N_G(Q)$. But $K$ is a Carter subgroup of $N_G(Q)$, hence $x\in K$ and $K$ is a
Carter subgroup of~$G$.
\end{proof}

\begin{df}
A finite group $G$ is said to satisfy {\bfseries
(ESyl2)}\index{condition!Sylow $2$-sub\-gro\-up is contained in Carter
subgroup {\bfseries
(ESyl2)}}\glossary{ESyl2@{\bfseries(ESyl2)}}, if for its Sylow $2$-sub\-gro\-up
$Q$ the equality $N_G(Q)=QC_G(Q)$ holds. In other words,  $G$ satisfies
{\bfseries(ESyl2)}, if every element of odd order, normalizing a Sylow
$2$-sub\-gro\-up $Q$ of $G$, centralizes~$Q$.
\end{df}

\begin{lem}\label{Syl2centrcomposit}
Let $Q$ be a Sylow $2$-sub\-gro\-up of a finite group $G$ and $x$ be an element
of odd order from $N_G(Q)$. Assume that there exist normal subgroups
$G_1,\ldots,G_k$ of $G$ such that
$G_1\cap\ldots\cap G_k\cap Q\leq Z(N_G(Q))$ and $x$ centralizes $Q$ modulo
$G_i$ for all~$i$.

\noindent Then $x$ centralizes $Q$.  In particular, if $G/G_i$ satisfies
{\em\bfseries (ESyl2)} for all $i$, then $G$
satisfies~{\em\bfseries(ESyl2)}.
\end{lem}

\begin{proof}
Consider the normal series $Q\unrhd Q_1\unrhd \ldots\unrhd Q_k\unrhd
Q_{k+1}=\{e\}$, where $Q_i=Q\cap(G_1\cap\ldots\cap G_i)$. The conditions
of the lemma imply that $x$ centralizes each factor $Q_{i-1}/Q_i$. Since
$x$ is an element of odd order, this implies that $x$ centralizes~$Q$.
\end{proof}

\begin{lem}\label{InhBy2-ext}
Let $H$ be a subgroup of a finite group $G$ such that $\vert G:H\vert=2^t$, $H$
satisfies {\em \bfseries (ESyl2)}, and each element of odd order of $G$ is in
$H$ (this condition is  evidently equivalent to the subnormality
of~$H$).

\noindent Then $G$ satisfies~{\em\bfseries (ESyl2)}.
\end{lem}

\begin{proof}
Let $Q$ be a Sylow $2$-sub\-gro\-up of $G$ such that $Q\cap H$ is a Sylow
$2$-sub\-gro\-up of $H$. Consider an element $x\in N_G(Q)$ of odd order. Since
$x\in H$, then $x\in N_H(Q)\leq N_H(Q\cap H)=(Q\cap H)\times O(N_H(Q\cap H))$,
i.~e.,  $x\in O(N_H(Q\cap H))$. Thus the set of elements of odd order in
$N_G(Q)$ forms a subgroup $R=O(N_H(Q\cap H))\cap N_G(Q)$ of $N_G(Q)$. Clearly
$R$ is normal in $N_G(Q)$, therefore $R=O(N_G(Q))$. On the other hand, $Q$ is
normal in $N_G(Q)$ by definition and $Q\cap R=\{e\}$, whence $N_G(Q)=Q\times
O(N_G(Q))$.
\end{proof}

By using the result of this section, we shall improve the definition of minimal
counter example.

\begin{df}\label{MinimalCounterExample}
A finite almost simple group $A$ is called a {\em minimal counter
example}\index{minimal counter example}, if it contains nonconjugate Carter
subgroups, but Carter subgroups of every almost simple group, of order less
than $\vert A\vert$ with simple socle being a  {\slshape known} simple group,
are conjugate.
\end{df}

\section{Conjugacy in simple groups}

\subsection{Brief review of results of the section}

Recall that in view of Lemma \ref{power} none element from the center of a
Carter subgroup can be conjugate to its nontrivial power (if the centralizer of
the element satisfies {\bfseries(C)}). Thus if we would be able to prove that
each element of prime order $r$ of $G$ is conjugate to its nontrivial power and
at the same time its centralizer satisfies {\bfseries(C)}, then we may state
that order of a Carter subgroup (if it exists) is not divisible by~$r$.

In this section we obtain the information on the conjugacy of elements of
prime order in finite simple groups and, by using this information, we obtain a
description of Carter subgroups in a wide class of almost simple groups.
Actually, in almost simple groups, distinct from $A_n^\varepsilon(q)$
($\varepsilon=\pm$), Carter subgroups should be $2$-gro\-ups, as made clear
below. The results can be formulated as a list of almost simple groups $A$-s,
that cannot be minimal counter example (see Theorem~\ref{main}). This list is
summarized in the Table  \ref{AlmSimpleNotCounterEx}, where
 ${\rm Field}(S)$\glossary{FieldS@${\rm Field}(S)$} stands for the group
generated by field and inner-diagonal automorphisms of a finite group of Lie
type~$S$.

\begin{longtable}{|c|c|}\caption{Finite almost simple groups, which are not
minimal counter examples}\label{AlmSimpleNotCounterEx}\\  \hline
Soc($A$)$=G$&\ Conditions for $A$\\
\hline alternating, sporadic;&\\
$A_1(p^{t})$, $B_\ell(p^{t})$,\
$C_\ell(p^{t})$, $t$ is even if $p=3$;&\\
${^2B_2}(2^{2n+1})$,$ G_2(p^{t})$, $F_4(p^{t})$, ${^2F_4}(2^{2n+1})$;&
\\
$E_7(p^{t})$, $p\not=3$; $E_8(p^{t})$, $p\not=3,5$&none\\
${^3D_4(p^{t})}$, $D_{2\ell}(p^{t})$,  ${^2D_{2\ell}}(p^{t})$,
& \\
 $t$ is even if $p=3$ in the last $2$ cases and,& \\
 if $G=D_4(p^t)$, $|({\rm Field}(G)\cap A):(\widehat
G\cap A)|_{2'}>1$ &\\
\hline $B_\ell(3^t)$, $C_\ell(3^t)$, $D_{2\ell}(3^t)$,
${^3D_4(3^{t})}$,  ${{^2D_{2\ell}}(3^{t})}$, &\\
${D_{2\ell+1}(r^{t})}$, ${^2D_{2\ell+1}}(r^{t})$,
${^2G_2(3^{2n+1})}$, & $A=G$\\
$E_6(p^{t})$, ${^2E_6}(p^{t})$,
$E_7(3^t)$, $E_8(3^t)$, $E_8(5^t)$&\\
\hline
\end{longtable}

In particular,  $A$ cannot be simple (case $A=A_\ell^\varepsilon(q)$ is
excluded by Theorem~\ref{CarterInClassicalGroups}).

\subsection{Preliminary results}

\begin{lem}\label{CentrOfInvolution}
Let $\ov{G}$ be a simple connected algebraic group over a field of
characteristic $p$,  $t$ be an element of order $r$ of $\ov{G}$, not divisible
by~$p$.

\noindent Then $C_{\ov{G}}(t)/C_{\ov{G}}(t)^0$ is a $\pi(r)$-gro\-up.
\end{lem}

\begin{proof}
Since $p$ does not divide $r$, then $t$ is semisimple.  By Lemma
\ref{CentrOfSemisimpleClassic}, $C_{\ov{G}}(t)^0$ is a connected reductive
subgroup of maximal rank of $\ov{G}$ and every $p$-ele\-ment of $C_{\ov{G}}(t)$
is contained in $C_{\ov{G}}(t)^0$. Assume that a prime $s\not\in\pi(r)$ divides
order $\vert C_{\ov{G}}(t)/(C_{\ov{G}}(t)^0)\vert$. Then $s\not=p$ and
$C_{\ov{G}}(t)$ contains an element $x$ of order $s^k$ such that
$x\not\in C_{\ov{G}}(t)^0$. Since  $x,t$ commute, we have that $x\cdot t$ is a
semisimple element of $\ov{G}$ (of order $r\cdots^k$). Therefore, there exists a
maximal torus $\ov{T}$ of $\ov{G}$ containing $x\cdot t$. Then $(xt)^r=x^r\in
\ov{T}$. Since $(s,r)=1$, there exists $m$ such that
$rm\equiv 1\pmod{s^k}$, thus $(x^r)^m=x\in \ov{T}$. Since
$xt, x\in \ov{T}$, then $t\in \ov{T}$, so $\ov{T}\leq
C_{\ov{G}}(t)^0$, hence  $x\in C_{\ov{G}}(t)^0$; a contradiction.
\end{proof}

\begin{lem}\label{1} Let  $s\in \ov{G}$ be a semisimple element of order $r$
such that $(r,\Delta(\ov{G}))=1$.

\noindent Then $C_{\ov{G}}(s)$ is connected. In particular, it follows that for
every Frobenius map $\sigma$ of $\ov{G}$, two semisimple elements $s,s^\prime\in
\ov{G}_\sigma$ are conjugate in $\ov{G}_\sigma$ if and only if they are
conjugate in~$\ov{G}$.
\end{lem}

\begin{proof}
Follows from Lemmas \ref{FactorCentrByConnectedComponent}
and~\ref{CentrOfInvolution}.
\end{proof}

The following lemma plays an important role, since it shows that a
semisimple element of odd prime order is usually conjugate to its
inverse.

\begin{lem}\label{omnibus}
Let $G=O^{r^\prime}(\ov{G}_\sigma)$, $\ov{G}$ has an adjoint type and the root
system of $\ov{G}$ has type distinct from~${A_\ell\ (\ell >1),D_{2\ell+1},E_6}$.

\noindent Then each semisimple element of odd order  $s\in
\widehat{G}$ is conjugate to its inverse by an element of~$G$.
\end{lem}

\begin{proof}
There exists some $\sigma$-stab\-le maximal torus
$\ov{T}$ of $\ov{G}$ with $s\in \ov{T}$.  $\ov{T}$ is generated by the set
$\{h_\alpha(\lambda)\mid \alpha\in \Phi, \lambda\in \ov{\F}_p^\ast\}$ and
the factor group $N_{\ov{G}}(\ov{T})/\ov{T}$ is isomorphic to the Weyl group
$W$ of $\ov{G}$. If $w\in W$ and $n_w$ is a preimage of $w$ under the natural
epimorphism $N_{\ov{G}}(\ov{T})\to W$, then
$h_\alpha(\lambda)^{n_w}=h_{\alpha^w}(\lambda)$.  Now let $w_0$ be
the unique involution of $W$ such that $w_0(\Phi^+)=\Phi^-$ and
let $n_0$ be a preimage of $w_0$. Since we are assuming
$\Phi\ne A_\ell\ (\ell > 1)$, $D_{2\ell+1}$, and $E_6$, we have
$\alpha^{w_0}=-\alpha$ for all $\alpha\in \Phi$, hence
$h_\alpha(\lambda)^{n_0}=h_{-\alpha}(\lambda)=
h_\alpha(\lambda)^{-1}$.  We conclude that $s^{n_0}=s^{-1}$, i.e.
that $s$ is conjugate to $s^{-1}$ in $\ov{G}$. Thus, by the previous
Lemma, $s$ and $s^{-1}$ are conjugate in $\ov{G}_\sigma$. Finally,
from $\ov{G}_\sigma=\ov{T}_\sigma G$, we conclude that $s$ and $s^{-1}$ are
conjugate in~$G$.
\end{proof}

\begin{lem}\label{3} Let  $\ov{C}$ be a connected reductive subgroup of
maximal rank of $\ov{G}$. Denote by $W$ and $W_{\ov{C}}$ the Weyl groups of
$\ov{G}$ and $\ov{C}$ respectively, by $W_{\ov{C}}^\perp$ the subgroup of
$W$, generated by reflections in roots orthogonal to all roots from
$\Phi(\ov{C})$, and by $\Delta_{\ov{C}}$ the Dynkin diagram of $\ov{C}$.
Then:
\begin{itemize}
\item[{\em (a)}] $N_W(W_{\ov{C}})/(W_{\ov{C}}\times W_{\ov{C}}^\perp)\simeq
{\Aut}_W(\Delta_{\ov{C}})$;
\item[{\em (b)}]\enskip  $N_{\ov{G}}(\ov{C})/\ov{C}\simeq
N_W(W_{\ov{C}})/W_{\ov{C}}$.
\end{itemize}

Let $G=O^{p^\prime}(\ov{G}_\sigma)$ be split or one of the groups
${^2A_\ell(p^{t})},$ ${^2D_{2\ell+1}(p^{t})}$, ${^2E_6}(p^{t})$. If $s\in
G$ is a
semisimple element such that $C_{\ov{G}}(s)$ is connected and
$N_G(C_{\ov{G}}(s))>C_{\ov{G}}(s)$, then~${N_G(C_G(s))>C_G(s)}$.
\end{lem}

\begin{proof}
Point (а) may be found in \cite[Proposition~4]{Car5CentSemisimpleLie}.  As to
point (б), let $\ov{T}$  be a maximal torus of $\ov{G}$ contained in
$\ov{C}$, so that we may assume $W=N_{\ov{G}}(\ov{T})/\ov{T}$ and
$W_{\ov{C}}=N_{\ov{C}}(\ov{T})/\ov{T}$. All maximal tori of $\ov{C}$ are
conjugate in $\ov{C}$,  since $\ov{C}$ is connected.  It follows  easily
that $N_{\ov{G}}(\ov{C})=\ov{C} N_{N_{\ov{G}}(\ov{T})}(\ov{C})$.  Moreover it is
shown in \cite[Proposition~5]{Car5CentSemisimpleLie} that
$N_{N_{\ov{G}}(\ov{T})}(\ov{C})=N_{N_{\ov{G}}(\ov{T})}(N_{\ov{C}}(\ov{T})).$
Hence,
\begin{equation}\label{two}
\frac {N_{\ov{G}}(\ov{C})}{\ov{C}}=
\frac{\ov{C}N_{N_{\ov{G}}(\ov{T})}(N_{\ov{C}}(\ov{T}))}{\ov{C}}\simeq
\frac{N_{N_{\ov{G}(\ov{T})}(N_{\ov{C}}(\ov{T}))}}{N_{\ov{C}}(\ov{T})}\simeq
\frac{N_W(W_{\ov{C}})}{W_{\ov{C}}}. \end{equation}

Now, let $G=O^{p^\prime}(\ov{G}_\sigma)$ be as in the statement,
and set $\ov{C}=C_{\ov{G}}(s)$. Write $\sigma=\tau \varphi$, where
$\tau$ is the graph automorphism of $\ov{G}$ induced by a symmetry
$\rho$ of the Dynkin diagram of $\Phi=\Phi(\ov{G})$ and
$\varphi$ is a field automorphism. Now let $\tau$ be the isometry which
extends $\rho$ on Euclidean space  $\R\otimes_{\Z}\Z\Phi$. If $\ov{T}_1$ is a
$\sigma$-stab\-le maximal split torus of $\ov{G}$, then for each
$x\in N_{\ov{G}}(\ov{T}_1)/\ov{T}_1$, we have
$x^\sigma={^\tau x}$  (considering $N_{\ov{G}}(\ov{T}_1)/\ov{T}_1=W_1$ as a
group of isometries of $\R\otimes_{\Z}\Z\Phi$). Thus if $G$ is split, i.~e.,
$\rho=\tau=e$, then $\sigma$ acts trivially on $W_1$.
If $G$  is twisted, hence of type  $A_\ell$, $D_{2\ell+1}$, or $E_6$, it is
possible to show directly, that ${-\tau}\in W_1$. Thus we may twist $\ov{T}_1$
by $-\tau$, obtaining the $\sigma$-stab\-le torus $(\ov{T}_1)_{{-\tau}}$.
By {equation} \eqref{twisted}:
$$\frac {(N_{\ov{G}}((\ov{T}_1)_{{-\tau}})))_\sigma}
{((\ov{T}_1)_{{-\tau}}))_\sigma }\simeq
C_{{W_1},\sigma}(-\tau)=
\{x\in {W_1}\mid {^\tau x}(-\tau) x^{-1}= {-\tau}\}={W_1}.$$
Let $\{X_\alpha\mid \alpha\in\Phi\}$ be the set of $\ov{T}_1$-ro\-ot subgroups
and set $\ov{C}_1=\la \ov{T}_1, X_\alpha\mid \alpha\in \Phi(\ov{C})\ra$. Since
$\Phi(\ov{C})$ is $\sigma$-in\-va\-ri\-ant it follows that
$\ov{C}_1$  is $\sigma$-stab\-le. Moreover, since
$\tau(\Phi(\ov{C}))=\Phi(\ov{C})$, we have that
$-\tau\in N_{W_1}(W_{\ov{C}_1})$. By
\cite[Proposition 1 and~2]{Car5CentSemisimpleLie}, it follows that there exists
$(\ov{C}_1)_{-\tau}$ obtained from $\ov{C}_1$ by twisting with $-\tau$.  Up to
conjugation in $G$ we may assume that
$(\ov{T}_1)_{-\tau}\leq (\ov{C}_1)_{-\tau}$. Define
$\ov{T}_0=\ov{T}_1$ and $\ov{C}_0=\ov{C}_1$ if $G$  is split, and
$\ov{T}_0=(\ov{T}_1)_{-\tau}$ and $\ov{C}_0=(\ov{C}_1)_{-\tau}$ if $G$ is
twisted.

Since $\Phi(\ov{C})=\Phi(\ov{C}_0)$, there exists $g\in \ov{G}$, such that
${^g\ov{C}_0}=C$  and ${^g\ov{T}_0}=\ov{T}$. It follows that
$\dot w=g^{-1}\sigma (g)\in N_{\ov{G}}(\ov{C}_0)\cap N_{\ov{G}}(\ov{T}_0)$.
So the image $w$ of $\dot w$ in $W_0=N_{\ov{G}}(\ov{T}_0)/\ov{T}_0$
belongs to~${N_{W_0}(W_{C_0})}$.

From $\ov{G}_\sigma=\ov{T}_\sigma G$ it follows
$(N_{\ov{G}}(\ov{C}))_\sigma=N_{\ov{G}_\sigma}(\ov{C})=\ov{T}_\sigma
N_G(\ov{C})$. Hence we are done if we can show that it is nontrivial the
group
$\frac{(N_{\ov{G}}(\ov{C}))_\sigma}{\ov{C}_\sigma}=\frac{\ov{T}_\sigma
N_G(\ov{C})}{\ov{T}_\sigma C_G(s)} \simeq
\frac{N_G(\ov{C})}{C_G(s)}$,  which is a subgroup
of~$\frac{N_G(C_G(s))}{C_G(s)}.$

Using {equation} \eqref{two} we get
\begin{equation*}
\frac{(N_{\ov{G}}(\ov{C}))_\sigma}{\ov{C}_\sigma}\simeq\frac{(N_{N_{\ov{G}}(\ov{
T } ) } (N_{\ov{C}}(\ov{T}))/\ov{T})_\sigma}{
(N_{\ov{C}}(\ov{T})/\ov{T})_\sigma }
\simeq \frac{N_{N_{\ov{G}}(\ov{T})}(N_{\ov{C}}(\ov{T}))/\ov{T}\cap
(N_{\ov{G}}(\ov{T})/\ov{T})_\sigma}
{N_{\ov{C}}(\ov{T})/\ov{T}\cap
(N_{\ov{G}}(\ov{T})/\ov{T})_\sigma}.\end{equation*}
By our choice of $\ov{T}_0$ we have
$\frac {N_{\ov{G}}(\ov{T}_0)}{\ov{T}_0}=\left(\frac
{N_{\ov{G}}(\ov{T}_0)}{\ov{T}_0}\right)_\sigma $, i.~e. $\sigma$  acts trivially
on the finite group $\frac {N_{\ov{G}}(\ov{T}_0)}{\ov{T}_0}$. Now, if $w\in
W_{\ov{C}_0}$, by \cite[Proposition~1]{Car5CentSemisimpleLie}  we may
assume $w=e$, $\ov{T}=\ov{T}_0$, $\ov{C}=\ov{C}_0$. It follows
$N_{\ov{G}}(\ov{T})/\ov{T}=(N_{\ov{G}}(\ov{T})/\ov{T})_\sigma$,
hence
$$\frac{(N_{\ov{G}}(\ov{C}))_\sigma}{\ov{C}_\sigma}\simeq
\frac{N_{W}(W_{\ov{C}})}{W_{\ov{C}}}\simeq \frac{N_{\ov{G}}(\ov{C})}{\ov{C}}$$
which is non-trivial by assumption.
Finally assume that $w\not\in W_{\ov{C}_0}$, i.~e.
$\dot w=g^{-1}\sigma(g)\not\in \ov{C}_0.$ It follows that  $^g\dot
w=\sigma(g)g^{-1}\not\in \ov{C}$, i.~e.
${^g\dot w}\ov{T}\not\in N_{\ov{C}}(\ov{T})/\ov{T}$.
On the other hand  $^g\dot w\in N_{\ov{G}}(\ov{C})\cap
N_{\ov{G}}(\ov{T})$.
Moreover, since $\sigma$ acts trivially on
$N_{\ov{G}}(\ov{T}_0)/\ov{T}_0$, we have that $\sigma(\dot w
\ov{T}_0)=\dot w\ov{T}_0$, i.~e.
$\sigma(g)^{-1}g\sigma(g)^{-1}\sigma^2(g)=t_0\in \ov{T}_0$. Hence
$^gt_0=t\in \ov{T}$ and $^{\sigma(g)g^{-1}}t= (^g\dot w)^{-1}\sigma(^g\dot w)\in
\ov{T}$. It follows that $\sigma({^g\dot w}\ov{T})= {^g\dot
w}\ov{T}$.
So, if $w\not\in W_{\ov{C}_0}$, we conclude that
${^g\dot w}\ov{T}$ maps onto a non trivial element of  the group
$\frac{N_{N_{\ov{G}}(\ov{T})}(N_{\ov{C}}(\ov{T}))/\ov{T}\cap
(N_{\ov{G}}(\ov{T})/\ov{T})_\sigma} {N_{\ov{C}}(\ov{T})/\ov{T}\cap
(N_{\ov{G}}(\ov{T})/\ov{T})_\sigma}.$
\end{proof}

The rest of this subsection is devoted to unipotent elements
in groups of Lie type.

\begin{lem}\label{Seitz} Let  $G=O^{p^\prime}(\ov{G}_\sigma)$ be a finite group
of Lie type with the base field $\F_{p^t}$, with $p$ odd.
If $p=3$, suppose $t$ even.  Assume further that
$\Phi(\ov{G})\not=G_2,F_4,E_6,E_7,E_8$ if $p=3$,  and
$\Phi(\ov{G})\not=E_8$ if $p=5$.

\noindent Then every unipotent element $u$ of order $p$ is conjugate in $G$ to
some power~$u^k\not= u$.\end{lem}

\begin{proof}
Under our assumptions $p$ is a good prime. By point
({\it i}) of  \cite[Theorem~1.4]{SeitzConjUnipElts}, there exists a closed
$\sigma$-stab\-le subgroup $A_1(\ov{\F}_p)$ of $\ov{G}$ such that $u\in
A_1(\ov{\F}_p)$. Clearly $O^{p^\prime }((A_1(\ov{\F}_p))_\sigma)$
is isomorphic either to $\SL_2(p^{tm})$, or to $\P\SL_2(p^{tm})$,
for some positive integer $m>0$. Up to conjugation inside
$A_1(\ov{\F}_p)$,  we may assume $u=\left(\begin{array}{cc}
1&\zeta\\ 0&1\\ \end{array}\right)$ (or is equal to the projective image of
this matrix) for some $\zeta\in \F_{p^{tm}}$.   Under our assumptions,
there exist $\eta\in \F_{p^t}$ such that  $1\not=
\eta^2=k\in \F_p$.  Let $x$ be the matrix $\left(\begin{array}{cc}
\eta^{-1}&0\\ 0&\eta\\
\end{array}\right)$ or its projective image. Then
$x\in G$, and $u$,
$u^x=\left(\begin{array}{cc} 1&k\zeta\\ 0&1\\
\end{array}\right)=u^k$ are conjugate in~$G$.
\end{proof}

\begin{lem}\label{G2}
Let $u\in G=G_2(3^t)$ be an element of
order $3$.

\noindent Then $u$ is conjugate to  $u^{-1}$ in~$G$.
\end{lem}

\begin{proof}
By \cite[Proposition~6.4]{En}  there exist $9$ unipotent conjugacy classes in
$G$. All of them may be found in Table  \ref{UnipotentG2}, where $\alpha,\beta$
 denote respectively a short and a long fundamental root of  $G_2$,
$\zeta$ is an element of $\F_{3^t}$ such that the polynomial
$x^3-x+\zeta$ is irreducible in $\F_{3^t}[x]$ and $\eta$ is a non-square of
$F_{3^t}$. Since $|x_1|=9$ and $x_2$, $x_3$  are conjugate to $x_1$
in $G_2(\ov{\F}_3)$, we only need to verify that $x_4,x_5,x_6,x_7,x_8$  are
conjugate to their inverses. Using the formulae $x_\beta (u)^{h_\alpha
(t)}= x_\beta (t^{\frac {2(\alpha,\beta)}{(\alpha,\alpha)}}u)$ for each $\alpha
,\beta\in \Phi$ (see \cite[Proposition~6.4.1]{CarSimpleGrpsLieType}), we
get: $x_6^{h_\alpha(-1)}=x_6^{-1}$,
$x_8^{h_\beta(-1)}=x_8^{-1}$, $x_4^{h_\beta(-1)}=x_4^{-1}$, and
$x_5^{h_\beta (-1)}= x_5^{-1}$. Finally  $|C_K(x_7)|\not=
|C_K(x_i)|$ for all $i\not=7$:  thus also   $x_7$ is conjugate to its inverse.
\end{proof}

\begin{longtable}{|l||l|}\caption{Unipotent classes
in $G_2(q), \ q=3^t$.\label{UnipotentG2}}\\ \hline
representative $x$&$|C_K(x)|$\\ \hline $x_0=1$&$q^6(q^2-1)(q^6-1)$\\
$x_1=x_\alpha(1)x_\beta(1)$&$3q^2$\\
$x_2=x_\alpha(1)x_\beta(1)x_{3\alpha+\beta}(\zeta)$&$3q^2$\\
$x_3=x_\alpha(1)x_\beta(1)x_{3\alpha+\beta}(-\zeta)$&$3q^2$\\
$x_4=x_{\alpha+\beta}(1)x_{3\alpha+\beta}(1)$&$2q^4$\\
$x_5=x_{\alpha+\beta}(1)x_{3\alpha+\beta}(\eta)$&$2q^4$\\
$x_6=x_{2\alpha+\beta}(1)$&$q^6(q^2-1)$\\
$x_7=x_{2\alpha+\beta}(1)x_{3\alpha+2\beta}(1)$&$q^6$\\
$x_8=x_{3\alpha+2\beta}(1)$&$q^6(q^2-1)$\\ \hline
\end{longtable}

\begin{lem}\label{F4}
Let $u\in G=F_4(3^t)$ be an element of order~$3$.

\noindent Then $u$ is conjugate to $u^{-1}$ in~$G$.
\end{lem}

\begin{proof}
By \cite[Table~6]{Sh} there exist $28$ unipotent conjugacy classes of $G$. All
of them may be found in Table \ref{UnipotentF4}.  Recall that in an Euclidean
$4$-di\-men\-si\-o\-nal space with orthonormal base
$\varepsilon_1,\varepsilon_2,\varepsilon_3,\varepsilon_4$ all
roots of $F_4$ may be written as $\{\pm
\varepsilon_i\pm\varepsilon_j,\pm\varepsilon_i,\frac {1}{2}
(\pm\varepsilon_1\pm\varepsilon_2\pm\varepsilon_3\pm\varepsilon_4)\}$.
In Table \ref{UnipotentF4} the symbols $\pm i\pm j$, $\pm i$, and
$\pm1\pm2\pm3\pm4$ denote the roots $\pm\varepsilon_i
\pm\varepsilon_j$, $\pm\varepsilon_i$, and $\frac {1}{2}
(\pm\varepsilon_1\pm\varepsilon_2\pm\varepsilon_3\pm\varepsilon_4)$
respectively, $\eta$ is a fixed non-square element of $\F_{3^t}$,
$\xi$ is a fixed element of $\F_{3^t}$ such that $x^2+\xi
x+\eta$ is an irreducible polynomial in $\F_{3^t}[x]$, $\zeta$ is a fixed
element of $\F_{3^t}$ such that $x^3-x+\zeta$ is an
irreducible polynomial in $\F_{3^t}[x]$.   By using \cite[Table~7]{Sh} one may
easily verify that $|x_9|=|x_{10}|>3,|x_i|>3$ for all
$i\geq 12$.  Indeed, by \cite[Table~7]{Sh} we have that elements $x_9$ and
$x_{10}$ are conjugate in $F_4(\ov{\F}_3)$. They also are conjugate to an
element $c_7=x_{r_1}(1)x_{r_2}(1)x_{r_3}(1)$, where the roots $r_1,r_2$
and $r_3$ я are fundamental roots in a root system of type $A_3$.  But it is
evident, that $|c_7|>3$.  In all cases when $|x_i|>3$, we proceed in a similar
way. In the remaining cases one can see that $|C_K(x_i)|\not=|C_K(x_j)|$  for
all $i\not=j$. So if $|x_i|=3$, then $i=1,2,3,4,5,6,7,8,11$, and $x_i$ is
conjugate to its inverse under~$G$.
\end{proof}

{\footnotesize
\begin{longtable}{|l||l|}\caption{Unipotent classes $F_4(q),\
q=3^t$\label{UnipotentF4}}\\ \hline
representative $x$&$|C_K(x)|$\\ \hline $x_0=1$&
$|K|$\\
$x_1=x_{1+2}(1)$&$q^{24}(q^2-1)(q^4-1)(q^6-1)$\\
$x_2=x_{1-2}(1)x_{1+2}(-1)$&$2q^{21}(q^2-1)(q^3-1)(q^4-1)$\\
$x_3=x_{1-2}(1)x_{1+2}(-\eta)$&$2q^{21}(q^2-1)(q^3+1)(q^4-1)$\\
$x_4=x_2(1)x_{3+4}(1)$&$q^{20}(q^2-1)^2$\\
$x_5=x_{2-3}(1)x_4(1)x_{2+3}(1)$&$2q^{17}(q^2-1)(q^3-1)$\\
$x_6=x_{2-3}(1)x_4(1)x_{2+3}(\eta)$&$2q^{17}(q^2-1)(q^3+1)$\\
$x_7=x_2(1)x_{1-2+3+4}(1)$&$q^{14}(q^2-1)(q^6-1)$\\
$x_8=x_{2-3}(1)x_4(1)x_{1-2}(1)$&$q^{16}(q^2-1)$\\
$x_9=x_{2-3}(1)x_{3-4}(1)x_{3+4}(-1)$&$2q^{12}(q^2-1)^2$\\
$x_{10}=x_{2-3}(1)x_{3-4}(1)x_{3+4}(-\eta)$&$2q^{12}(q^4-1)$\\
$x_{11}=x_{2+3}(1)x_{1+2-3-4}(1)x_{1-2+3+4}(1)$&$q^{14}(q^2-1)$\\
$x_{12}=x_{2-3}(1)x_4(1)x_{1-4}(1)$&$2q^{12}(q^2-1)$\\
$x_{13}=x_{2-3}(1)x_4(1)x_{1-4}(\eta)$&$2q^{12}(q^2-1)$\\
$x_{14}=x_{2-4}(1)x_{3+4}(1)x_{1-2}(-1)x_{1-3}(-1)$&$24q^{12}$\\
$x_{15}=x_{2-4}(1)x_{3+4}(1)x_{1-2}(-\eta)x_{1-3}(-1)$&$8q^{12}$\\
$x_{16}=x_{2-4}(1)x_{2+4}(-\eta)x_{1-2+3+4}(1)x_{1-3}(-1)$&$4q^{12}$\\
$x_{17}=x_{2-4}(1)x_{3+4}(1)x_{1-2-3+4}(1)x_{1-2}(-\eta)x_{1-3}(\xi)$&$4q^{
12}$\\
$x_{18}=x_{2}(1)x_{3+4}(1)x_{1-2+3-4}(1)x_{1-2}(-1)x_{1-3}(\zeta)$&$3q^{12}$\\
$x_{19}x_{2-3}(1)x_{3-4}(1)x_4(1)$&$q^8(q^2-1)$\\
$x_{20}=x_2(1)x_{3+4}(1)x_{1-2-3-4}(1)$&$q^8(q^2-1)$\\
$x_{21}=x_{2-4}(1)x_3(1)x_{2+4}(1)x_{1-2-3+4}(1)$&$2q^8$\\
$x_{22}=x_{2-4}(1)x_3(1)x_{2+4}(\eta)x_{1-2-3+4}(1)$&$2q^8$\\
$x_{23}=x_{2-3}(1)x_{3-4}(1)x_4(1)x_{1-2}(1)$&$2q^6$\\
$x_{24}=x_{2-3}(1)x_{3-4}(1)x_4(1)x_{1-2}(\eta)$&$2q^6$\\
$x_{25}=x_{2-3}(1)x_{3-4}(1)x_4(1)x_{1-2-3-4}(1)$&$3q^4$\\
$x_{26}=x_{2-3}(1)x_{3-4}(1)x_4(1)x_{1-2-3-4}(1)x_{1-2+3+4}(\zeta)$&$3q^4$\\
$x_{27}=x_{2-3}(1)x_{3-4}(1)x_4(1)x_{1-2-3-4}(1)x_{1-2+3+4}(-\zeta)$&$3q^4$\\
\hline\end{longtable}}


\begin{lem}\label{E6}
Let $u\in G$ be an element of order $3$, where $G=E_6(3^t)$
or $G={^2E_6}(3^{t})$ is a canonical finite group of Lie type.

\noindent Then $u$ is conjugate to $u^{-1}$ in~$G$.
\end{lem}

\begin{proof}
Let $\ov{G}$ and $\sigma$ be such that $G=O^{p^\prime}(\ov{G})$. Since
the characteristic equals $3$, we have that $Z(\ov{G}_{sc})=1$. So we may
assume $\ov{G}=\ov{G}_{sc}$ to be universal. Thus  $\ov{G}$ is simply connected
and $G=\ov{G}_\sigma$. We assemble the information from
\cite[Lemmas 4.2, 4.3, 4.4, and Theorem~4.13]{Mi} on conjugacy classes of
unipotent elements of $G$  in Table \ref{UnipotentE6}.
In Table \ref{UnipotentE6} we substitute the root $\alpha_1
r_1+\alpha_2r_2+\alpha_3r_3+\alpha_4r_4+\alpha_5r_5+ \alpha_6r_6$, where
$r_1,r_2,r_3,r_4,r_5,r_6$ form a fundamental system
of~$E_6$, by the 6-tuple $\alpha_1\alpha_2\alpha_3\alpha_4\alpha_5\alpha_6$
of its coefficients.

Note that if $n\geqslant 3$ and $r_1,r_2,\ldots,r_n$ are fundamental
roots of root system of type $A_n$, then $|x_{r_1}(1)x_{r_2}(1)\ldots
x_{r_n}(1)|>3$.  By using this fact we obtain, that $|x_4|>3$,
$|x_7|>3$, $|x_8|>3$, $|x_i|>3$, where $i\geqslant 10$, $i\not=12,16$. Thus we
have to consider remaining cases only. We have that
$x_1^{h_{r_1}(\lambda)}=x_1^{-1}$, where $\lambda$ is a square root of $-1$ in
$\ov{\F}_3$.  For each $x\in \ov{G}$ denote by  ${\rm Ccl}(x)$  its conjugacy
class in $\ov{G}$. Since
$C_{\ov{G}}(x_1)=C_{\ov{G}}(x_1)^0$, from \cite[Theorem~8.5]{Hu2ConjClasses}
we have that for every Frobenius map $\sigma$ and for every $x\in {\rm
Ccl}(x_1) \cap \ov{G}_\sigma$, the elements $x$ and $x^{-1}$ are conjugate under
$\ov{G}_\sigma$. So if $x\in {\rm Ccl}(x_1) \cap G$, then $x$ is
conjugate to its inverse.

{\footnotesize
\begin{longtable}{|l||c|}\caption{Unipotent classes in
$E_6(\ov{\F}_3)$\label{UnipotentE6}}\\ \hline
representative $x$& $C= C_G(x)$\\
&$|C:C^0|$\\ \hline $x_1=x_{10000}(1)$&$1$\\
$x_2=x_{100000}(1)x_{001000}(1)$&$2$\\
$x_3=x_{100000}(1)x_{000100}(1)$&$1$\\
$x_4=x_{100000}(1)x_{001000}(1)x_{000100}(1)$&$1$\\
$x_5=x_{100000}(1)x_{001000}(1)x_{000010}(1)$&$1$\\
$x_6=x_{100000}(1)x_{000100}(1)x_{000001}(1)$&$1$\\
$x_7=x_{100000}(1)x_{001000}(1)x_{000100}(1)x_{000010}(1)$&$1$\\
$x_8=x_{100000}(1)x_{001000}(1)x_{000100}(1)x_{000001}(1)$& $1$\\
$x_9=x_{100000}(1)x_{001000}(1)x_{000010}(1)x_{000001}(1)$&$1$\\
$x_{10}=x_{100000}(1)x_{001000}(1)x_{010000}(1)x_{000010}(1)$&$1$\\
$x_{11}=x_{100000}(1)x_{001000}(1)x_{000100}(1)x_{010000}(1)
x_{000001}(1)$&$1$\\ $x_{12}=x_{100000}(1)x_{001000}(1)
x_{000010}(1)x_{000001}(1)x_{010000}(1)$&$1$\\
$x_{13}=x_{100000}(1)x_{001000}(1)x_{000100}(1)x_{000010}(1)
x_{000001}(1)$&$1$\\ $x_{14}=x_{010000}(1)x_{001000}(1)
x_{000100}(1)x_{000010}(1)$&$1$\\ $x_{15}=x_{010000}(1)x_{001000}(1)
x_{000100}(1)x_{010110}(1)$&$6$\\ $x_{16}=x_{000001}(1)x_{000010}(1)
x_{001000}(1)x_{010000}(1)$&$1$\\ $x_{17}=x_{010000}(1)x_{001000}(1)
x_{000010}(1)x_{101100}(1)$&$1$\\ $x_{18}=x_{000010}(1)x_{000100}(1)
x_{001000}(1)x_{100000}(1)x_{000001}(1)x_{111111}(1)$&$2$\\
$x_{19}=x_{010000}(1)x_{000100}(1)x_{000010}(1)x_{000001}(1)
x_{101000}(1)x_{001110}(1)$&$1$\\ $x_{20}=x_{100000}(1)x_{010000}(1)
x_{001000}(1)x_{000100}(1)x_{000010}(1)x_{000001}(1)$&$3$\\ \hline
\end{longtable}}

For the other $x_i$-s such that $3$, with $i\not=2$, we proceed in
the same way. We are left with $x_2$. By \cite[Theorem~8.5]{Hu2ConjClasses} we
have that,  for every Frobenius map $\sigma$, ${\rm Ccl}(x_2) \cap
\ov{G}_\sigma$ consists of two conjugacy classes of $G=\ov{G}_\sigma$.
Assume first that $G=E_6(3^t)$. Then, by \cite[Lemmas~4.2 and~4.4]{Mi}
we have, that if $x\in {\rm Ccl}(x_2) \cap G$, then $x$ is conjugate under
$G$ either to $y_1=x_{100000}(1)x_{001000}(1)$, or to
$y_2=x_{100000}(1)x_{001000}(1)x_{000001}(1)x_{122321}(\eta)$, where
$\eta$ is a nonsquare in $\F_{3^t}$. By \cite[Lemma~4.2]{Mi}
$|C_G(y_1)|=2q^{26}(q^2-1)^2(q^3-1)^2$, by \cite[Lemma~4.4]{Mi}
$|C_G(y_2)|=2y^{26}(q^4-1)(q^6-1)$. For $i=1,2$, let
${\rm Ccl}_G(y_i)$ be the conjugacy class of $y_i$ in $G$. Since
$|C_G(y_1)|\not=|C_G(y_2)|$ we have, that $y_i$ is conjugate
to its inverse under $G$ for $i=1,2$. So if $x\in {\rm Ccl}_G(y_1) $, or
$x\in {\rm Ccl}_G(y_2)$, then $x$ is conjugate to its inverse under~$G$.
Now assume that $G={^2E_6(3^{t})}$ and denote $E_6(3^{2t})$ by $G_1$. Then
$G=(G_1)_\tau$ for some graph automorphism $\tau$ of $G_1$.
There exists a Frobenius map $\sigma$ such that
$G_1=\ov{G}_\sigma$, $G=\ov{G}_{\sigma\tau}$
(see~\cite[(7-2)]{GorLyoLocalStructure}). Let ${\rm Ccl}_1$ and
${\rm Ccl}_2$  be two conjugacy classes of $G_1$ contained in ${\rm
Ccl}(x_2)\cap G_1$. We prove that every $x\in {\rm Ccl}_i$, $i=1,2$, is
conjugate to $x^{-1}$ в $G_1$. Since
${\rm Ccl}(x_2)\cap G$ consists of two conjugacy classes of $G$, we have that
${\rm Ccl}_1\cap G$  consists of one conjugacy class
and ${\rm Ccl}_2\cap G$ consists of one conjugacy class. So, every $x\in {\rm
Ccl}_i\cap G$, $i=1,2$  is conjugate to its inverse under~$G$.
\end{proof}

\begin{lem}\label{centUH}
Let $O^{p'}(\overline{G}_\sigma)\leq G\leq\overline{G}_\sigma$ be a finite
adjoint group of Lie type over a field of odd characteristic $p$ and the root
system $\Phi$ of $\overline{G}$ is one of the following:
$A_n$ $(n\ge2)$, $D_n$ $(n\ge
4)$, $B_n$ $(n\ge 3)$, $G_2$, $F_4$, $E_6$, $E_7$ or $E_8$; and $G\not\simeq
{}^2G_2(3^{2n+1})$. Let $U$ be a maximal unipotent subgroup of $G$, $H$ be a
Cartan subgroup of $G$, normalizing $U$, and $Q$ is a Sylow
$2$-sub\-gro\-up of~$H$.

\noindent Then $C_U(Q)=\{e\}$.
\end{lem}

\begin{proof}
Clearly we enough to prove the lemma for the case
$G=O^{p'}(\ov{G}_\sigma)=O^{p'}(G)$, i.~e., we may assume that $G$ is a
canonical adjoint group of Lie type.

First assume that $G$ is split. Assume that
$C_U(Q)\not=\{e\}$ and $u\in C_U(Q)\setminus\{e\}$. Consider decomposition
\eqref{canonicalform} of $u=\prod_{r\in\Phi^+}x_r(t_r),$ where $t_r$ are
from the definition field $\F_q$ of $G$. In view of
\cite[Theorem~5.3.3(ii)]{CarSimpleGrpsLieType} this decomposition is
unique. Since for every $h(\chi)\in H$, $r\in\Phi$, $t\in
\F_q$ the formulae $h(\chi)x_r(t)h(\chi)^{-1}=x_r(\chi(r)t)$ holds
(see \cite[p.~100]{CarSimpleGrpsLieType}), then we obtain that each multiplier
$x_r(t_r)$ in decomposition \eqref{canonicalform} of $u$ is in
$C_U(Q)$. So we may assume that $u=x_r(t)$ for some $r\in\Phi^+$ and
$t\in \F_q^\ast$. Under our restriction on $\Phi$, by Hart\-ley-Shute lemma
\ref{Hartley-Shute}, there exists $h(\chi)\in H$ such that $\chi(r)=-1$. Since
$h(\chi)^2=h(\chi^2)$ (see \cite[p.~98]{CarSimpleGrpsLieType}), then we have
that $\chi^2(r)=1$, i.~e., $\vert h(\chi)^2\vert<\vert
h(\chi)\vert$. Hence, $\vert h(\chi)\vert$ is even and we may write
$h(\chi)=h_2\cdot h_{2'}=h(\chi_1)\cdot h(\chi_2)$, a decomposition of
$h(\chi)$ as a product of its $2$- and $2'$- parts. Now
$\chi(r)=\chi_1(r)\cdot\chi_2(r)$, therefore $\chi_1(r)=-1$ and $\chi_2(r)=1$.
Thus $h(\chi_1)x_r(t)h(\chi_1)^{-1}=x_r(-t)\not=x_r(t)$. Since
$h(\chi_1)\in Q$, the obtained equation contradicts to the choice of~${x_r(t)\in
C_U(Q)}$.

Assume that $G\simeq {}^2A_n(q)$, $G\simeq{}^2D_n(q)$, or
$G\simeq {}^2E_6(q)$, then $\Phi(\ov{G})$ equals $A_n$, $D_n$ and
$E_6$ respectively. Denote by $\bar{r}$ the image of $r$ of
$\Phi$ under the corresponding symmetry. In terms of
\cite{CarSimpleGrpsLieType}, the root system $\Phi(\ov{G})$ is expressible as a
union of equivalency classes $\Psi_i$-s, while each $\Psi_i$ has type either
$A_1$, or $A_1\times A_1$, or $A_2$. In view of
\cite[Proposition~13.6.1]{CarSimpleGrpsLieType}, the equality
$U=\prod_i X_{\Psi_i}$ holds, where $$X_{\Psi_i}=\{x_r(t)\mid t\in
\F_q\},$$ if $\Psi_i=\{r\}$ has type $A_1$ (here $r=\bar{r}$);
$$X_{\Psi_i}=\{x_r(t)x_{\bar{r}}(t^q)\mid t\in \F_{q^2}\},$$ if
$\Psi_i=\{r,\bar{r}\}$ has type $A_1\times A_1$ (here $r\not=\bar{r}$, and
$r+\bar{r}\not\in\Phi(\ov{G})$);
$$X_{\Psi_i}=\{x_r(t)x_{\bar{r}}(t^q)x_{r+\bar{r}}(u)\mid t\in \F_{q^2},
u+u^q=-N_{r,\bar{r}} t t^q\},$$ if $\Psi_i=\{r,\bar{r},r+\bar{r}\}$ has type
$A_2$ (here $r\not=\bar{r}$ and $r+\bar{r}\in\Phi(\ov{G})$). Now if
$h(\chi)$ is an element of $H$, then the following equalities hold  (see
\cite[p.~263]{CarSimpleGrpsLieType}):
$$h(\chi)x_r(t)h(\chi)^{-1}=x_r(\chi(r)t),$$ if $r=\bar{r}$ and $\Psi_i=\{r\}$
has type $A_1$;
$$h(\chi)x_r(t)x_{\bar{r}}(t^q)h(\chi)^{-1}=
x_r(\chi(r)t)x_{\bar{r}}(\chi(\bar{r }) t^q),$$ if $r\not=\bar{r}$,
$r+\bar{r}\not\in\Phi(\ov{G})$ and $\Psi_i=\{r,\bar{r}\}$ has type $A_1\times
A_1$; $$h(\chi)x_r(t)x_{\bar{r}}(t^q)x_{r+\bar{r}}(u)h(\chi)^{-1}=
x_r(\chi(r)t) x_{\bar{r}}(\chi(\bar{r}) t^q)x_{r+\bar{r}}(\chi(r+ \bar{r})u),$$
if $r\not=\bar{r}$, $r+\bar{r}\in\Phi(\ov{G})$ and
$\Psi_i=\{r,\bar{r},r+\bar{r}\}$ has type~$A_2$.

Let $u$ be a nontrivial element from $C_U(Q)$. Then $u$ contains a nontrivial
multiplier from $X_{\Psi_i}$ for some $i$. In view of uniqueness of
decomposition into the product $\prod_iX_{\Psi_i}$ (see
\cite[Proposition~13.6.1]{CarSimpleGrpsLieType}) we may assume that $u\in
X_\Psi$.

Assume that $\Psi$ has type $A_1$, i.~e., $u=x_r(t)$, $t\in
\F_q$, $r=\bar{r}$. In view of Hart\-ley-Shute lemma \ref{Hartley-Shute}, for
each $s\in \F_q$ there exists $h(\chi)\in H$ such that $\chi(r)=s$. Take $s=-1$.
Then there exists  $h(\chi)\in H$ such that $\chi(r)=-1$. Since
$h(\chi)^2=h(\chi^2)$ (see formulae on p. 98 from \cite{CarSimpleGrpsLieType}),
then we have that $\chi^2(r)=1$, i.~e. $\vert h(\chi)^2\vert<\vert
h(\chi)\vert$. Hence order $\vert h(\chi)\vert$ is even and we may write
$h(\chi)=h_2\cdot h_{2'}=h(\chi_1)\cdot h(\chi_2)$, a decomposition of
$h(\chi)$ into the product of its $2$- and $2'$- parts. Now
$\chi(r)=\chi_1(r)\cdot\chi_2(r)$, therefore $\chi_1(r)=-1$ and
$\chi_2(r)=1$. Thus $h(\chi_1)x_r(t)h(\chi_1)^{-1}=x_r(-t)\not=x_r(t)$. So the
case $u=x_r(t)$ and $\Psi=\{r\}$ has type $A_1$ is impossible.

Assume that  $\Psi=\{r,\bar{r}\}$ has type $A_1\times A_1$. By Hart\-ley-Shute
lemma \ref{Hartley-Shute} for every $s\in \F_{q^2}$
there exists $h(\chi)\in H$ such that $\chi(r)=s^2$. Since there exists
$s\in \F_{q^2}$ such that $s^2=-1$, then there exists
$h(\chi)\in H$ such that $\chi(r)=-1$. As above $h(\chi)$ can
be written as $h(\chi_1)\cdot h(\chi_2)$, a product of its $2$- and
$2'$- parts. Then $\chi_1(r)=-1$, so $$h(\chi_1)x_r(t)
x_{\bar{r}}(t^q)h(\chi_1)^{-1}=x_r(-t) x_{\bar{r}}(-t^q)\not =
x_r(t) x_{\bar{r}}(t^q).$$ Thus the case $u=x_r(t)
x_{\bar{r}}(t^q)$ and $\Psi=\{r,\bar{r}\}$  has type $A_1\times A_1$
is impossible.

Assume that $\Psi=\{r,\bar{r},r+\bar{r}\}$ has type
$A_2$. By Hurt\-ley-Shute lemma \ref{Hartley-Shute}, for each $s\in
\F_{q^2}$ there exists $h(\chi)\in H$ such that $\chi(r)=s^3$. Choose $s=-1$,
then there exists $h(\chi)\in H$ such that $\chi(r)=-1$. Again
$h(\chi)=h(\chi_1)\cdot h(\chi_2)$ is expressible as the product of its
$2$- and $2'$- parts and~${\chi_1(r)\not=1}$. Then
\begin{multline*}
h(\chi_1)x_r(t)
x_{\bar{r}}(t^q)x_{r+\bar{r}}(u)h(\chi_1)^{-1}=\\ x_r(-t)
x_{\bar{r}}(\chi_1(-t^q) x_{r+\bar{r}}(\chi_1(r+\bar{r}) u)\not =\\
x_r(t)
x_{\bar{r}}(t^q)x_{r+\bar{r}}(u)
\end{multline*}
for $t\not=0$. If $t=0$, then choose $s$ so that $s^2=-1$.  Then
$\chi_1(r+\bar{r})=-1$ and, as above, we obtain the inequality.
Hence this case is impossible.

Assume at last that $G\simeq{}^3D_4(q)$. In terms from
\cite{CarSimpleGrpsLieType}, a root system $\Phi(\ov{G})$ is expressible as a
union of equivalency classes $\Psi_i$, when each $\Psi_i$ has type either
$A_1$, or $A_1\times A_1\times A_1$. In view of
\cite[Proposition~13.6.1]{CarSimpleGrpsLieType}, the equality
$U=\prod_i X_{\Psi_i}$ holds, where $$X_{\Psi_i}=\{x_r(t)\mid t\in
\F_q\},$$ if $\Psi_i=\{r\}$ has type $A_1$ (here $r=\bar{r}$);
$$X_{\Psi_i}=\{x_r(t)x_{\bar{r}}(t^q)x_{\bar{\bar{r}}}(t^{q^2})\mid t\in
\F_{q^3}\},$$ if $\Psi_i=\{r,\bar{r},\bar{\bar{r}}\}$ has type
$A_1\times A_1\times A_1$ (here $r\not=\bar{r}$ and
$r+\bar{r}\not\in\Phi(\ov{G})$). In both cases, by Hart\-ley-Shute lemma
\ref{Hartley-Shute}, there exists $h(\chi)\in H$ such that
$\chi(r)=-1$. As above we may assume that $h(\chi)$ is a
$2$-ele\-ment, i.~e. $h(\chi)\in Q$ and $h(\chi)$ does not centralizes
nonidentical elements from~$X_{\Psi_i}$, and the statement of Lemma
\ref{centUH} follows in the last case.
\end{proof}

\begin{lem} \label{unordthree}
In the notations of Lemma {\em \ref{centUH}}, with $p$ odd,  let $K$ be a
Carter subgroup of $G$ such that $|K|=2^a p^b$.

\noindent Then $a>0$. More precisely, up to conjugation, $O_p(K)\leq
C_U(Q)$. In particular, under the assumptions of Lemma {\em
\ref{centUH}}, $K$ is a $2$-gro\-up.
\end{lem}

\begin{proof}
The condition $a=0$ would imply $K=U$. But $U$ is normalized by $H$ which is
non-trivial as $p$ is odd and $G$ is simple. Thus $a>0$. Now, assume $b>0$.
By Bo\-rel-Tits theorem (Lemma \ref{Borel-Titsclassic}), $K$ i contained in a
proper parabolic subgroup $P$ of $G$ and $O_p(K)\leq O_p(P)$. Since $P=L
O_p(P)$, where  $L$ is a Levi factor of $P$,  from Lemma \ref{HomImageOfCarter}
it follows that $KO_p(P)/O_p(P)\cong O_2(K)$ is a Carter subgroup of
$P/O_p(P)\cong L$. Thus $O_2(K)$ is a Sylow $2$-sub\-gro\-up of $L$. But $L$
contains $H$, therefore we may assume that  $Q\leq K$.  It follows
that~${O_p(K)\leq C_U(Q)}$.
\end{proof}

\begin{lem}\label{notl}
Let  $G$ be a non-Abelian simple group not of Lie
type.

\noindent Then every element $z$ of odd order is conjugate to some~${z^k\not=
z}$.
\end{lem}

\begin{proof}
By the classification of finite simple groups, $G$ is either alternating, or
sporadic. Our claim can be checked
directly in first case, and using
the description of the conjugacy classes given in  \cite{ATLAS} in the second
case.
\end{proof}

\subsection[Almost simple groups]{Almost simple groups which are not
minimal counter examples}

In this subsection $A$ denotes a minimal almost simple group that is a minimal
counter example (see definition \ref{MinimalCounterExample}). If $G$ is a group
of Lie type, we denote by ${\rm Field}(G)$ the subgroup of $\Aut(G)$ generated
by inner, diagonal, and field automorphisms. If $G$ is a simple group which is
not of Lie type, we set ${\widehat{G}=G}$ to unify notations. More over for
each $x\in G$ we assume that composition factors of the centralizer $C_G(s)$
are known simple groups, and so $C_G(s)$ satisfies {\bfseries(C)}. As we noted
in subsection 2.4, this assumption is always true. We say it here in order to
emphasize that all results do not depend on the classification of finite simple
groups.

\begin{lem}\label{teor}
Let $A$ be a minimal counter example and $G=F^\ast(A)$. Assume that for every
element $z\in \widehat{G}$ of odd prime order, $z$ is conjugate to some
$z^k\not=z$ in~$G$.

\noindent Then $A$ is not a minimal counter example if one of the following
holds:
\begin{itemize}
\item[{\em (a)}] $\vert A:\widehat{G}\cap A\vert$ is a $2$-po\-w\-er;
\item[{\em (b)}]  $|\widehat{G}:(\widehat{G}\cap A)|$ is a $2$-po\-w\-er and ,
if $\Phi(\ov{G})$ has type $D_4$, then $|({\rm
Field}(G)\cap ):(\widehat{G}\cap A)|_{2'}>1$;
\item[{\em (c)}] for every odd prime $r$ and every Sylow $r$-sub\-gro\-up $R$ of $A$,
either $R\cap G$ has no complement in $R$, or all such complements are
conjugate in~$A$.
\end{itemize}
\end{lem}

\begin{proof}
Let $K,H$ be nonconjugate Carter subgroups of $A$. Note that by Lemma
\ref{power}(b) it follows that $K\cap \widehat{G}$ and $H\cap \widehat{G}$ are
$2$-gro\-ups. We prove (c) first, that we show that (a) and (b) follows
from~(c).

(c) By Theorem \ref{ConjugacyCriterion} and Lemma \ref{HomImageOfCarter}, we
obtain that $KG/G=HG/G=A/G$. In particular, if $r$ is a prime divisor of $\vert
A/G\vert$, then $r$ divides both $\vert K\vert$ and $\vert H\vert$. By Lemma
\ref{power} and by conditions of this lemma, it follows that $K\cap
\widehat{G}$ and $H\cap \widehat{G}$ do not contain elements of odd prime
order, i.~e., are $2$-gro\-ups. If $R\cap G$ has no a complement in $R$ we get
a contradiction immediately, if all such complements are conjugate in $A$, we
obtain a contradiction with Lemma \ref{power}(c). Thus we obtain that $\vert
A/G\vert$ is a $2$-po\-w\-er, hence $K$ and $H$ are $2$-gro\-ups, that is
impossible.

Now (a) evidently follows from (c). As to (b), then it also follows from (c),
by using the conjugacy of complements, that follows by
Lemma~\ref{ConjAutomorphisms}.
\end{proof}

Note that all non-Abelian composition factors of the centralizer of every
element of the alternating group $\Alt_n$ are alternating groups of lower
degree. So Lemmas \ref{notl} and \ref{teor} and induction by $n$ imply
immediately that Carter subgroups of $\Aut(\Alt_n)$ with $n\ge5$ either are
Sylow $2$-sub\-gro\-ups or do not exist. The same statement holds for sporadic
groups. Thus the following statement is true.

\begin{lem}\label{notlCarter}
Let $S$ be a finite non-Abelian simple group, that is either sporadic, or
alternating.

\noindent Then, for every subgroup $A$ of $\Aut(S)$, a Carter subgroup either
does not exist, or is a Sylow $2$-sub\-gro\-up.
\end{lem}

\begin{ttt}\label{orto}
Let $G$ be a finite adjoint group of Lie type such that $G=
\P\Omega_{2(2\ell+1)}^\pm (p^t)$, and assume that  $\ell\geqslant 2$.

\noindent Then $G$ is not a minimal counter example.
\end{ttt}

\begin{proof}
Assume that our statement is false. Then $G$ contains a Carter subgroup $K$,
that is not a $2$-gro\-up. Let $s\in Z(K)$ be an element of odd prime order
$r$. Then we may assume that $s$ is semisimple, except, probably, the case,
when $p\not=2$ and $\vert K\vert= 2^a p^b$. But this is impossible in view of
Lemmas \ref{centUH} and \ref{unordthree}. Hence $s$ is semisimple and from
$K\leq C_G(s)$ it follows that $C_G(s)$ is self-nor\-ma\-li\-zing in $G$ (see
Lemma~\ref{power}(a)). Now let  $\ov{G}=\Omega_{2(2\ell+1)}(\ov{\F}_p)$ and
$\sigma$ be such that $\ov{G}_\sigma= \Omega_{2(2\ell+1)}^\pm (p^t)$. More over
set $K_0$ to be equal to  the preimage of $K$ in $\ov{G}_\sigma$. Clearly $K_0$
is a Carter subgroup of $\ov{G}_\sigma$ and we may identify  $s$ with its
preimage in $\ov{G}_\sigma$, since the center of $\ov{G}_\sigma$ has order $2$
or $4$. Since $\vert s\vert$ is odd, Lemma  \ref{CentrOfInvolution} implies
that $C=C_{\ov{G}}(s)$ is a connected reductive subgroup of maximal rank of
$\ov{G}$ (see Lemmas \ref{CentrOfSemisimpleClassic} and
\ref{FactorCentrByConnectedComponent}). More over $C$ is a proper subgroup of
$\ov{G}$, since $s\notin Z(\ov{G})$. By Lemma \ref{3} the group $N_W(W_C)/W_C$
is isomorphic to $N_{\ov{G}}(C)/C$. By using the description of $N_W(W_C)/W_C$,
given in \cite[Proposition~10]{Ca2CentsSemisimpleClassical} and Lemma \ref{3},
we conclude that $N_G(C_G(s))/C_G(s)$ is trivial only if $W_C^\perp$ and ${\rm
Aut}_W(\Delta_C)$ are both trivial. From assumption $\ell\geq 2$ it follows
that this occurs precisely when $m_1=0$ and $m_{2\ell+1}=1$ (in the notations
from \cite{Ca2CentsSemisimpleClassical}). In this case
$C=A_{2\ell}(\ov{\F}_p)*S$, where $S$ is a 1-di\-men\-si\-o\-nal torus. By
using the fact that $\ov{G}$ contains exactly one class of connected reductive
subgroups isomorphic to $C$, and assuming that $\ov{G}$ preserves the bilinear
form induced by $J=\left(\begin{array}{cc}0&I\\ I&0\\
\end{array}\right)$, we may identify $C$ with the image of $\GL_{2\ell
+1}(\ov{\F}_p)$ under a monomorphism $\varphi$ such that
$$A\mapsto \left(\begin{array}{cc}A&0\\ 0&(A^{-1})^t\\
\end{array}\right) \ .$$
By Lang-Ste\-in\-berg theorem (Lemma \ref{LangSteinbergTheorem}), we may assume
that either $C_\sigma = \varphi(\GL_{2\ell +1}(p^t))$, or $C_\sigma=\varphi
(\GU_{2\ell +1}(p^{2t}))$. Since $K_0$ is a Carter subgroup of $C_\sigma$ and
$\ell\geq 2$, by \cite{DT2CartSbgrpsPGL}, \cite{DTZCartSbgrpsClassGrps}, and
Theorem \ref{CarterInClassicalGroups} it follows that $K_0$ is the normalizer
of a Sylow $2$-sub\-gro\-up $P$ of $C_\sigma$, and either $p^t=2$ (and
$C_\sigma=\varphi(\GL_{2\ell +1}(p^t))$), or $p$ is odd. From $s\in
Z(C_\sigma)$ it follows that $r=\vert s\vert$ divides $p^t-1$ if
$C_\sigma\simeq \GL_{2\ell +1}(p^t)$, and that $r$ divides $p^t+1$ if
$C_\sigma\simeq \GU_{2\ell +1}(p^{2t})$. In particular $p$ is odd. By using
known structure of normalizers of Sylow $2$-sub\-gro\-ups in classical groups
(see \cite{KondNormalizers} and \cite{CarFongNormSylow2}), we may assume that
$K_0$ is a subgroup of the following group
$$L=\left\{ \left(\begin{array}{cc}B&0\\ 0&\beta\\
\end{array}\right)^\varphi\mid B\in \GL_{2\ell}(p^{t})\ ,\
\beta\in \ov{\F}_q^\ast\right\}\ \text{if } C_\sigma\simeq \GL_{2\ell
+1}(p^{t})$$
$$L=\left\{\left(\begin{array}{cc}B&0\\ 0&\beta\\
\end{array}\right)^\varphi\mid B\in \GU_{2\ell}(p^{2t})\ ,\
\beta^{p^{t}+1}=1\right\}\ \text{if } C_\sigma\simeq \GU_{2\ell +1}(p^{2t}).$$
As we noted above, there exists $y\in L$ such that
$y=\left(\begin{array}{cc}I_{2\ell}&0\\
0&\gamma\\
\end{array}\right)^\varphi$ where $\gamma$ has order $r$. Since $y$ is in the
center of $L$, it is also in the center of $K_0$. Thus $K_0\leq
C_{C_{G_\sigma}(s)}(y)= \big (C_C(y)\big )_\sigma.$ From isomorphism $C\simeq
\GL_{2\ell +1}(\ov{\F}_p)$ it follows that $C_C(y)$ is a connected reductive
$\sigma$-in\-va\-ri\-ant subgroup of maximal rank of $\ov{G}$. Thus, in view of
above mentioned result by Carter
\cite[Proposition~10]{Ca2CentsSemisimpleClassical}, $\big (C_C(y)\big )_\sigma$
is self-nor\-ma\-li\-zing in $\ov{G}_\sigma$ only if $C_C(y)$ is conjugate to
$C$. But $\dim(C_C(y))< \dim(C)$, since $y$ is not in the center of $C$. Thus
$(C_C(y))_\sigma$ is not self-normalizing in  $\ov{G}_\sigma$. Since
$Z(\ov{G})\leq C_C(y)$, it follows that the factor group
$(C_C(y))_\sigma/(Z(\ov{G}))_\sigma$ is not self-normalizing in $\ov{G}_\sigma
/(Z(\ov{G}))_\sigma =G$. Thus we have obtained a contradiction with
Lemma~\ref{power}(a), since $K$ is contained in
$(C_C(y))_\sigma/(Z(G))_\sigma$, and  $(C_C(y))_\sigma/(Z(G))_\sigma$
satisfies~{\bfseries(C)}.
\end{proof}

\begin{ttt} \label{care6}
Let $E_6^\varepsilon(p^t)\leq G\leq \widehat{E_6^\varepsilon(p^{t})}$. Then 
$G$ is not a minimal counter example.
\end{ttt}

\begin{proof}
Assume that our claim is false. Then, by Lemma \ref{power}(c), $G$ admits a
Carter subgroup $K$,  which does not
contain any Sylow $2$-sub\-gro\-up of $G$.  In particular $K$ is not a
$2$-gro\-up. Let $s\in Z(K)$ have odd prime order $r$.  By Lemmas
\ref{Seitz}, \ref{E6}, and \ref{power}, $p$  does not divide $\vert
K\vert$. Hence $s$ is semisimple and $K$ is contained in $C_G(s)$,
which, in virtue of Lemma \ref{power}(a),  is self-normalizing.
If $|s|\not=3$, then, by Lemma \ref{1}, it follows that $C_{\ov{G}}(s)$ is
connected. If $|s|=3$, then, by Lemma \ref{FactorCentrByConnectedComponent}, it
follows that $|C:C^0|$ divides $\Delta=3$. Direct calculations by using
\cite{der1CentSemisimpleExcept} and \cite{Mi} show that $C_G(s)$ is not
self-nor\-ma\-li\-zing, if $|s|=3$. Therefore we may assume that $|s|\not=3$
and $C_{\ov{G}}(s)$ is connected. Since $C_G(s)$ is self-nor\-ma\-li\-zing,
Lemma \ref{3} shows that $C=C_{\ov{G}}(s)$  is self-normalizing as well. By
\cite{Mi}, we obtain that $C$  is self-normalizing if and only if
$C=A_4(\ov{\F}_p)\circ A_1(\ov{\F}_p)\circ S$, or
$C=D_5(\ov{\F}_p)\circ S$, where $S$ is a $1$-di\-men\-si\-o\-nal torus
of~$\ov{G}$.

If $C=A_4(\ov{\F}_p)\circ A_1(\ov{\F}_p)\circ S$, then like in proof of
Theorem \ref{orto}, we may find an element $y\in
Z(K)$ such that $|y|=r$ and $C_G(\langle s\rangle\times \langle y\rangle)$
is not self-nor\-ma\-li\-zing; a contradiction with Lemma~\ref{power}.

So, assume that $C=D_5(\ov{\F}_p)\circ S$.  Then $C_G(s)=C\cap G=HL$,
where $H$ is a Cartan subgroup of $G$ and $L=O^{p^\prime}(C_G(s))$ is
either $D_5(p^t)$ or ${^2D_5(p^{t})}$.   Since $|\hat L:L|$ divides $4$, then 
$$O_{2'}(H)=(O_{2'}(H)\cap Z(C_G(s)))\times (O_{2'}(H)\cap L).$$
Denoting by $Q$  a Sylow $2$-sub\-gro\-up of $C_G(s)$, we claim that 
$N_{C_G(s)}(Q)=QZ(C_G(s)).$  Indeed, let $x$ be an element of
$N_{C_G(s)}(Q)$. From $H=O_2(H)\times O_{2^\prime}(H)$ and $C_G(s)=HL$,
we can write $x=h_1zl$ with $h_1\in O_2(H)$, $z\in
O_{2'}(H)\cap Z(C_G(s))$, $l\in L$.  We may clearly assume
$O_2(H)\leq Q$: thus $l\in N_{C_G(s)}(Q)$. From $L$  normal in
$C_G(s)$, it follows $l\in N_{L}(Q\cap L)$.  By
\cite{KondNormalizers}, $N_L(Q\cap L)=Q\cap L$, so, $l\in Q$. We conclude that
$N_{C_G(s)}(Q)=QZ(C_G(s))$ is nilpotent, hence a Carter subgroup
of  $C_G(s)$. Since $C_G(s)<G$, all Carter subgroups in $C_G(s)$ are conjugate.
Therefore, up to conjugation, $K=N_{C_G(s)}(Q)$. In virtue of the formula 
$|(C)_\sigma|=|M_\sigma|\cdot|(Z(C)^0)_\sigma|$, where
$M_\sigma=L$ in our notation (see \cite{Ca2CentsSemisimpleClassical}), we have
that $|G:C_G(s)|$  is odd, so  $Q$ is a Sylow $2$-sub\-gro\-up of $G$,
a contradiction.
\end{proof}

Our results are summarized in the following theorem. 

\begin{ttt} \label{main}
An almost simple group $A$, with socle с цоколем $G$ is not a minimal counter
example in the following cases:

\begin{itemize}
\item[{\em (a)}]  $G$ is alternating, sporadic, or one of the following
groups: $A_1(p^t)$, $B_\ell(p^{t})$, $C_\ell(p^{t})$, where
$t$ is even if $p=3$; ${^2B_2}(2^{2n+1})$, $G_2(p^{t})$,
$F_4(p^{t})$, ${^2F_4}(2^{2n+1})$, ${}^3D_4(q)$; $E_7(p^{t})$, where
$p\not=3$; $E_8(p^{t})$, where $p\not=3,5$, $D_{2\ell}(p^{t})$,
${^3D_4(p^{t})}$, ${^2D_{2\ell}}(p^{t})$,
where $t$ is even if $p=3$ and, more over, if $G=D_4(p^t)$, then 
$|({\rm Field}(G)\cap A): (\widehat G\cap A)|_{2'}>1$;
\item[{\em (b)}]\enskip $A$ is one of the following groups:
$B_\ell(3^t)$,  $D_{2\ell}(3^t)$, ${^2D_{2\ell}}(3^{t})$,
${D_{2\ell+1}(p^{t})}$, ${^2D_{2\ell+1}}(r^{t})$, ${^3D_4(3^{t})}$,
${^2G_2}(3^{2n+1})$, $E_6^\varepsilon(r^{t})$,
$\widehat{E_6^\varepsilon(r^{t})}$, $E_7(3^t)$, $E_8(3^t)$,
$E_8(5^t)$, $C_\ell(3^t)$;
\end{itemize}
In particular, no simple group, can be a minimal counter example. More over,
if each almost simple group with known simple normal subgroup satisfies 
{\em\bfseries (C)}, then in all above mentioned groups a Carter subgroup (if
exists) contains a Sylow $2$-sub\-gro\-up.
\end{ttt}

\begin{proof}
(a) We claim that every element $z\in \widehat G$ of prime odd order is
conjugate, under $G$, to some power $z^k\not=z$. When $G$ is alternating or
sporadic this is true  by Lemma~\ref{notl},  and when $G$ is of Lie type and 
$z$ is semisimple, this is true by Lemma~\ref{omnibus}. On the other
hand, when $z$ is unipotent (hence $p$ is odd), our claim follows
from Lemmas \ref{G2}, \ref{F4}  if $G=G_2(3^t),\
F_4(3^t)$ and from Lemma \ref{Seitz} in the remaining cases. Finally, if
$G\simeq {}^3D_4(q)$, then by \cite[Theorem~1.2(vi)]{TiepZalRealConjClasses}
each element of $G$ is conjugate to its inverse. Thus (a) follows from
Lemma~\ref{teor}, since for all groups under consideration we have either that
$|\widehat G:G|$ is a power of $2$, and so by Lemma 
\ref{ConjAutomorphisms} all complements of odd order are conjugate, or that 
$|A:A\cap \widehat G|$ is a power of $2$ (see 
\cite{ATLAS}, for example).

(b) Our statement follows from the results obtained in
\cite{DTZCartSbgrpsClassGrps} and Theorem \ref{CarterInClassicalGroups}, when 
$G=B_2(3^t)\simeq C_2(3^t)$ or
$G=C_\ell(3^{t})$, and from Theorems \ref{orto} and \ref{care6}, when $G$ is
one of the groups  $D_{2\ell+1}^\varepsilon(p^{t})$,
$E_6^\varepsilon(p^t)$ or $\widehat{E_6^\varepsilon(p^{t})}$. So assume that we
are in the remaining cases. Every semisimple element $z\in \widehat G$  of prime
odd order is conjugate to some $z^{-1}$ by Lemma  \ref{omnibus}. Thus, in
characteristic $2$ a Carter subgroup $K$ of $G$ can only be a Sylow 
$2$-sub\-gro\-up and, in odd characteristic, $K$ can only have order $2^ap^b$.
If
$G\not={^2G_2}(3^{2n+1})$, then the assumptions of Lemma~\ref{centUH} are
satisfied and, by using Lemma~\ref{unordthree}, we conclude that $K$  is again a
$2$-gro\-up.

Now assume $G={^2G_2}(3^{2n+1})$ (here $n\ge 1$). Then $\vert
K\vert=2^a3^b$. Since the normalizer of a Sylow $2$-sub\-gro\-up of 
$G$ contains an element of order $7$ (see \cite{levnuzReeGrps}), then we obtain
that $b>0$. By Lemma \ref{Borel-Titsclassic}, $K$ is contained in a proper
parabolic subgroup $P$ of $G$. Since Lie rank of $G$ equal $1$, then $P$ is a
Borel subgroup,  i.~e. $P=U\leftthreetimes H$, where $H$ is a Cartan
subgroup and $U$ is a maximal unipotent subgroup of $G$. Since $P$ is solvable,
it satisfies {\bfseries(C)} and by Lemma \ref{HomImageOfCarter}, $KU/U$ is a
Carter subgroup of $P/U\simeq H$. But for $n\ge 1$ the subgroup  $H$ contains
an element of odd order,  so $K$ contains an element of prime odd order. A
contradiction with $\vert K\vert=2^a 3^b$.
\end{proof}

Note that after proving the statement that for every known finite simple group
$S$ and a nilpotent subgroup $N\leq \Aut(S)$, Carter subgroups of $\la N,S\ra$
are conjugate, Theorem \ref{main} would imply that Carter subgroups in the
groups mentioned in the theorem should contain a Sylow  $2$-sub\-gro\-up. By
Lemma \ref{CritSyl2Carter} this is possible only if the normalizer 
of a Sylow $2$-sub\-gro\-up $Q$ in  $A$ satisfies $N_A(Q)=Q C_A(Q)$, i.~e. if
and only if $A$ satisfies {\bfseries (ESyl2)}. In \cite{KondNormalizers} and
subsequent results of the present paper, simple groups
satisfying~{\bfseries (ESyl2)} are completely determined. More over Lemmas
\ref{Esyl2InnDiagExtension} and \ref{ESyl2InhFieldGraph} allow to
``lift'' the property {\bfseries(ESyl2)} from a simple group to an almost
simple group. Thus a complete classification of Carter subgroups in groups 
mentioned in Theorem~\ref{main} is known.

\section{Semilinear groups of Lie type}

In this section we shall give a definition of semilinear groups of Lie type and
generalize results about the structure of finite groups of Lie type for them.
We need this theory to finding Carter subgroups in extensions of groups of Lie 
type by field, graph, or graph-fi\-eld automorphisms in section 5. In the last
subsection of this paragraph we shall consider the existence of Carter subgroup
in semilinear groups, either containing a Sylow $2$-sub\-gro\-up, or contained
in the normalizer of a Borel subgroup.

\subsection{Basic definitions}
Now we define some overgroups of finite groups of Lie type. We first give a
more detailed description of a Frobenius map $\sigma$. Note that all maps in
this section are automorphisms, if  $\ov{G}$ is considered as an abstract
group, and they are endomorphisms, if it is considered as an algebraic group.
Since we use the maps to construct connected automorphisms of finite groups and
groups over algebraically closed field, we find it appropriate to call all maps
in this section by automorphisms. Let $\overline{G}$ be a connected simple
linear algebraic group of adjoint type over the algebraic closure 
$\ov{\F}_p$ of a finite field of positive characteristic $p$. Below, if we do
not say opposite, we shall consider groups of adjoint type. Choose a Borel
subgroup $\overline{B}$ of $\overline{G}$, let
$\overline{U}=R_u(\overline{B})$ be the unipotent radical of 
$\overline{B}$. There exists a Borel subgroup $\overline{B}^-$, satisfying 
$\overline{B}\cap\overline{B}^-=\overline{T}$, where
$\overline{T}$ is a maximal torus of $\overline{B}$ (hence of $\overline{G}$).
We partially duplicate the notations and the definitions of subsection
1.3 here. Let $\Phi$ be the root system of $\overline{G}$ and let
$\{X_r\mid r\in\Phi^+\}$ be the set of $\overline{T}$-in\-va\-ri\-ant
$1$-di\-men\-si\-o\-nal subgroups of $\overline{U}$. Each $X_r$ is isomorphic
to the additive group of $\ov{\F}_p$, so each element of $X_r$ can be written as
$x_r(t)$, where $t$ is the image of $x_r(t)$ under this isomorphism. Denote by 
$\overline{U}^-=R_u(\overline{B}^-)$ the unipotent radical of $\overline{B}^-$.
Define like above $\overline{T}$-in\-va\-ri\-ant $1$-di\-men\-si\-o\-nal
subgroups $\{\overline{X}_r\mid r\in\Phi^-\}$ of $\overline{U}^-$. Then 
$\overline{G}=\langle \overline{U},\overline{U}^-\rangle$. Let 
$\bar{\varphi}$ be a field automorphism of $\overline{G}$ (as an abstract group)
and $\bar{\gamma}$ be a graph automorphism of $\overline{G}$. An automorphism
$\bar{\varphi}$ is known to may be chosen so that it acts by
$x_r(t)^{\bar{\varphi}}=x_r(t^p)$ (see \cite[12.2]{CarSimpleGrpsLieType} and
\cite[1.7]{CarFiniteGrpsLieTypeConjClassesCharacters}, for example). In view of
\cite[Propositions~12.2.3 and~12.3.3]{CarSimpleGrpsLieType}, we may choose
$\bar{\gamma}$ so that it acts by 
$x_r(t)^{\bar{\gamma}}=x_{\bar{r}}(t)$, if $\Phi$ has no roots of distinct
length, and by $x_r(t)^{\bar{\gamma}}=x_{\bar{r}}(t^{\lambda_r})$ for suitable
$\lambda_r\in\{1,2,3\}$, if $\Phi$ has roots of distinct length. Recall that
$\bar{r}$ is the image of $r$ under the symmetry $\rho$ (corresponding to 
$\bar{\gamma}$) of a root system $\Phi$. In both cases we can write 
$x_r(t)^{\bar{\gamma}}=x_{\bar{r}}(t^{\lambda_r})$, where
$\lambda_r\in\{1,2,3\}$. From this formulas it is evident that
$\bar{\varphi}\cdot \bar{\gamma}=\bar{\gamma}\cdot\bar{\varphi}$. Let 
$n_r(t)=x_r(t)x_{-r}(-t^{-1})x_r(t)$ and $\overline{N}=\langle
n_r(t)\mid r\in\Phi, t\in\ov{\F}_p^\ast\rangle$. Let
$h_r(t)=n_r(t)n_r(-1)$ and $\overline{H}=\langle h_r(t)\mid r\in\Phi,
t\in\ov{\F}_p^\ast\rangle$. By \cite[Chapters~6 and~7]{CarSimpleGrpsLieType},
$\overline{H}$ is a maximal torus of $\overline{G}$,
$\overline{N}=N_{\overline{G}}(\overline{H})$, and subgroups
$\overline{X}_r$ are root subgroups with respect to 
$\overline{H}$. So we can substitute $\overline{T}$ by 
$\overline{H}$ and assume that $\overline{T}$ is $\bar{\varphi}$- and
$\bar{\gamma}$- invariant under our choice. More over,  $\bar{\varphi}$ induces
a trivial automorphism of~$\overline{N}/\overline{H}$. Note that
$\ov{H}\leq \ov{B}\cap \ov{B}^-$, therefore $\ov{H}=\ov{T}$.

An automorphism $\bar{\varphi}^k$ with $k\in\N$ is called a {\em
classical Frobenius automorphism}\index{automorphism!Frobenius!classical}. We
shall call an automorphism $\sigma$ by a  {\em Frobenius
automorphism}\index{automorphism!Frobenius}, if $\sigma$ is conjugate under
$\ov{G}$ with $\bar{\gamma}^\epsilon\bar{\varphi}^k,\epsilon\in\{0,1\}, k\in
\N$. By Lang-Ste\-in\-berg theorem (Lemma \ref{LangSteinbergTheorem}) it
follows that for every $\bar{g}\in\ov{G}$ elements $\sigma$ and $\sigma\bar{g}$
are conjugate under $\ov{G}$. Thus by \cite[11.6]{SteLangSteinbergTheorem}, we
have that a Frobenius map, defined in subsection 1.4, coinsides with a
Frobenius automorphism, defined here.

Now fix $\overline{G}$, $\bar{\varphi}$, $\bar{\gamma}$, and
$\sigma=\bar{\gamma}^\epsilon\bar{\varphi}^k$; and assume that 
$\vert\bar\gamma\vert\le 2$, i.~e., we do not consider a triality automorphism
of a group $\ov{G}$ with root system $\Phi(\ov{G})=D_4$. Set
$B=\overline{B}_\sigma$, $H=\overline{H}_\sigma$, and
$U=\overline{U}_\sigma$. Since $\overline{B}, \overline{H},$ and
$\overline{U}$ are $\bar{\varphi}$- and $\bar{\gamma}$- invariant, they give us
respectively a Borel subgroup, a Cartan subgroup, and a maximal unipotent
subgroup (a Sylow $p$-sub\-gro\-up) of $\overline{G}_\sigma$ (for more
details see~\cite[1.7--1.9]{CarFiniteGrpsLieTypeConjClassesCharacters} or
\cite[Chapter~2]{GorLySol}).

Assume that $\epsilon=0$, i.~e.,  $O^{p'}(\overline{G}_\sigma)$ is not twisted
(is split). Then  $U=\langle X_r\mid r\in\Phi^+\rangle$, where $X_r$ is
isomorphic to the additive group of $\F_{p^k}=\F_{q}$, and each element of
$X_r$ can be written as $x_r(t), t\in \F_{q}$. Set also
$U^-=\overline{U}^-_\sigma$. As for $U$, we may write 
$U^-=\langle X_r\mid r\in \Phi^-\rangle$ and each element of
$X_r$ can be written as $x_r(t),t\in \F_{q}$. Now we can define an automorphism 
$\varphi$ as a restriction of $\bar{\varphi}$ on $\overline{G}_\sigma$ and
an automorphism $\gamma$ as a restriction of $\bar{\gamma}$ on
$\overline{G}_\sigma$. By definition the equalities $x_r(t)^\varphi=x_r(t^p)$
and $x_r(t)^\gamma=x_{\bar{r}}(t^{\lambda_r})$ hold for all $r\in\Phi$
(see the definition of $\bar{\gamma}$ above). Define an automorphism $\zeta$
of $\overline{G}_\sigma$ by
$\zeta=\gamma^\varepsilon\varphi^\ell,$\glossary{zeta@$\zeta$}
$\varphi^\ell\not=e$, $\varepsilon\in\{0,1\}$, and an automorphism 
$\bar{\zeta}$\glossary{zetabar@$\bar\zeta$} of  $\overline{G}$ by
$\bar\zeta=\bar\gamma^\varepsilon\cdot\bar\varphi^\ell$. Choose a 
$\zeta$-in\-va\-ri\-ant subgroup $G$ with $O^{p'}(\overline{G}_\sigma)\leq
G\leq\overline{G}_\sigma$. Note that if the root system $\Phi$ of $\ov{G}$ is
not equal to $D_{2n}$, then  $\ov{G}_\sigma/(O^{p'}(\ov{G}_\sigma))$ is cyclic.
Thus ro most groups and automorphisms, except groups of type $D_{2n}$ over a
field of odd characteristic, every subgroup $G$ of $\overline{G}_\sigma$, with 
$O^{p'}(\overline{G}_\sigma)\leq G\leq \overline{G}_\sigma$, is 
$\gamma$- and $\varphi$- invariant. Define $\Gamma
G$\glossary{GammaG@$\Gamma G$} as a set of subgroups of type $\langle
G, \zeta g\rangle\leq \ov{G}_\sigma\leftthreetimes\langle\zeta
\rangle$, where $g\in \ov{G}_\sigma$, $\langle\zeta g\rangle\cap
\ov{G}_\sigma\leq G$; and~${\Gamma
\overline{G}}$\glossary{GammaovG@$\Gamma\ov{G}$} as a set of subgroups of type
${\overline{G}\leftthreetimes \langle\bar\zeta\rangle}$. Following
\cite[Definition~2.5.13]{GorLySol}, an automorphism 
$\zeta$ is called a {\em field}\index{automorphism!of split group of
Lie type!field} automorphism, if $\varepsilon=0$, i.~e., $\zeta=\varphi^\ell$
and is called a {\em graph-fi\-eld}\index{automorphism!of split group of Lie
type!graph-fi\-eld} automorphism in the remaining cases (recall that we
are assuming~${\varphi^\ell\not=e}$).

Now assume that $\epsilon=1$, i.~e., $O^{p'}(\overline{G}_\sigma)$ is twisted.
Then $U=\overline{U}_\sigma$ и $U^-=\overline{U}^-_\sigma$. Define 
$\varphi$ on $U^\pm$ as a restriction of $\bar{\varphi}$ on $U^\pm$. Since
$O^{p'}(\overline{G}_\sigma)=\langle U^+,U^-\rangle$, we obtain an automorphism 
$\varphi$ of $O^{p'}(\overline{G}_\sigma)$. Consider
$\zeta=\varphi^\ell\not=e$\glossary{zeta@$\zeta$}, and let $G$ be a 
$\zeta$-in\-va\-ri\-ant group with 
$O^{p'}(\overline{G}_\sigma)\leq G\leq \overline{G}_\sigma$. Then 
$\bar{\zeta}=\bar{\varphi}^\ell$\glossary{zetabar@$\bar\zeta$} is an
automorphism of $\overline{G}$. Define $\Gamma G$\glossary{GammaG@$\Gamma G$}
as a set of subgroups of type $\langle G, \zeta g\rangle\leq
\ov{G}_\sigma\leftthreetimes\langle\zeta\rangle$, where $g\in
\ov{G}_\sigma$, $\langle\zeta g\rangle\cap \ov{G}_\sigma\leq G$;
and~${\Gamma\overline{G}}$\glossary{GammaovG@$\Gamma\ov{G}$} as a set of
subgroups of type~${\ov{G}\leftthreetimes\langle \bar{\zeta}\rangle}$. Following
\cite[Definition~2.5.13]{GorLySol}, we  say that  $\zeta$
is a {\em field}\index{automorphism!of twisted group of Lie type!field}
automorphism, if $\vert \zeta\vert $ is not divisible by 
$\vert\gamma\vert$ (this definition is used also in the case, when 
$\vert\gamma\vert=3$ and $\ov{G}_\sigma\simeq {}^3D_4(q)$), and that 
$\zeta$ is a {\em graph}\index{automorphism!of twisted group of Lie
type!graph} automorphism in the remaining cases.

Groups from above defined set $\Gamma G$ are called {\em semilinear finite
groups of Lie type}\index{group!finite of Lie type!semilinear}  (they are
also called {\em semilinear canonical finite groups of Lie
type}\index{group!finite of Lie type!canonical semilinear}, if
$G=O^{p'}(\ov{G}_\sigma)$), while groups from the set $\Gamma
\overline{G}$ are called {\em semilinear
algebraic groups}\index{group!algebraic!semilinear}. Note that 
$\Gamma\ov{G}$ cannot be defined without $\Gamma G$, since we need to know that 
$\varphi^\ell\not=e$. If $G$ is written in the notations from 
\cite{CarSimpleGrpsLieType}, i.~e., 
$O^{p'}(G)=G=A_n(q)$ or $O^{p'}(G)=G={^2A_n(q)}$ etc., then we shall write 
$\Gamma G$ as $\Gamma A_n(q)$, $\Gamma
{}^2A_n(q)$ etc.

Consider $A\in \Gamma G$ and $x\in A\setminus G$. Then $x=\zeta^k y$
for some $k\in\N$ and $y\in \ov{G}_\sigma$. Define $\bar{x}$ to be equal to
$\bar{\zeta}^ky$. Conversely, if $\bar{x}=\bar{\zeta}^k y$ for some $y\in
\ov{G}_\sigma$,
$\zeta^k\not=e$, and $\langle \zeta^k y\rangle\cap \ov{G}_\sigma\leq
G$, define $x$ to be equal to~$\zeta^k y$. Note that we need not to assume that
$\bar{x}\notin \overline{G}$, since~$\vert\bar{\zeta}\vert=\infty$. If 
$x\in G$, set~${\bar{x}=x}$.

\begin{lem}\label{equivnormalizer}
In above notations consider a subgroup $X$ of $G$. An element 
$x$ normalizes $X$ if and only if $\bar{x}$ normalizes $X$ as a subgroup
of~$\overline{G}$.
\end{lem}

\begin{proof}
Since $\zeta$ is a restriction of $\bar{\zeta}$ on $G$, our statement is
evident.
\end{proof}

Let $X_1$ be a subgroup of $A\in \Gamma G$. Then $X_1$ is generated by a normal
subgroup  $X=X_1\cap G$ and an element $x=\zeta^ky$. By Lemma 
\ref{equivnormalizer}, we may consider the subgroup 
$\overline{X}_1=\langle \bar{x},X\rangle$ of
$\ov{G}\leftthreetimes\langle\bar\zeta\rangle$. Now we find in reasonable 
to explain, why we use so complicated notations and definitions. We have that
order of  $\zeta$ is always finite, but order of $\bar\zeta$ is always
infinite. Thus, even if $Z(\ov{G})$ is trivial, we cannot consider 
$G\leftthreetimes\langle\bar\zeta\rangle$ as a subgroup of 
$\Aut(\ov{G})$. Therefore we need to define in a some way (one possible way is
just given) the connection between elements from $\Aut(G)$ and elements
from $\Aut(\ov{G})$, in order to use the machinery of linear algebraic groups.

Let $\overline{R}$ be a $\sigma$-stab\-le maximal torus
(respectively a reductive subgroup of maximal rank, a parabolic subgroup) of
$\overline{G}$, and an element $y\in N_{\ov{G}\leftthreetimes\langle
\bar\zeta\rangle}(\overline{R})$, is chosen so that there exists
$x\in \langle G,\zeta g\rangle$ with $y=\bar{x}$. Then $R_1=\langle
x,\overline{R}\cap G\rangle$ is called a {\em maximal torus}
(respectively a  {\em reductive subgroup of maximal rank}, a  {\em
parabolic subgroup}) of~$\langle G,\zeta g\rangle$.

\subsection{Translation of basic results}

\begin{lem}\label{parabolic}
Let $M=\langle x, X\rangle$, where $X=M\cap G\unlhd M$ is a subgroup of 
$\langle G,\zeta g\rangle$ such that $O_p(X)$ is nontrivial. Then there exists
a proper $\sigma$- and $\bar{x}$- invariant parabolic subgroup $\overline{P}$ of
$\overline{G}$ such that  $X\leq \overline{P}$ and~${O_p(X)\leq
R_u(\overline{P})}$.
\end{lem}

\begin{proof}
Define $U_0=O_p(X)$, $N_0=N_{\overline{G}}(U_0)$ and by induction 
$U_i=U_0 R_u(N_{i-1})$ and $N_i=N_{\overline{G}}(U_i)$. Clearly $U_i$,
$N_i$ are $\bar{x}$- and $\sigma$- invariant for all $i$. By 
\cite[Proposition~30.3]{HuLinearAlgGrps}, the chain of subgroups 
$N_0\leq N_1\leq \ldots\leq N_k\leq\ldots$ is finite and 
$\overline{P}=\cup_i N_i$ is a proper parabolic subgroup of $\ov{G}$. Clearly 
$\overline{P}$ is $\sigma$- and $\bar{x}$- invariant, $X\leq \ov{P}$
and~$O_p(X)\leq R_u(\ov{P})$.
\end{proof}

\begin{lem}\label{InvolutionsAndTori}
Let $O^{p'}(\ov{G}_\sigma)\leq G\leq\ov{G}_\sigma$ be a finite adjoint group of
Lie type with a base field of characteristic $p$ and order $q$. Assume also
that $O^{p'}(G)$ is not isomorphic to ${}^2D_{2n}(q)$, ${}^3D_4(q)$,
${}^2B_2(2^{2n+1})$, ${}^2G_2(3^{2n+1})$, ${}^2F_4(2^{2n+1}).$ Then there
exists a maximal $\sigma$-stab\-le torus $\ov{T}$ of $\ov{G}$ such that
\begin{itemize}
\item[{\em (a)}] $(N_{\ov{G}}(\ov{T})/\ov{T})_\sigma\simeq
(N_{\ov{G}}(\ov{T}))_\sigma/(\ov{T}_\sigma)=N(\ov{G}_\sigma,\ov{T}_\sigma)/\ov{T
}_\sigma \simeq W$, where $W$ is the Weyl group of~$\ov{G}$;
\item[{\em (b)}] if $r$ is an odd prime divisor of $q-(\varepsilon1)$, where
$\varepsilon=+$, if $G$ is split and $\varepsilon=-$, if $G$ is twisted, then 
$N(\ov{G}_\sigma,\ov{T}_\sigma)$ contains a Sylow
$r$-sub\-gro\-up~$\ov{G}_\sigma$;
\item[{\em (c)}] if $r$ is a prime divisor of $q-(\varepsilon1)$, and $s$ is
an element of order $r$ of $G$ such that $C_{\ov{G}}(s)$ is connected, then, up
to conjugation by an element of $G$, an element $s$ is contained
in~$T=\ov{T}_\sigma\cap G$;
\end{itemize}
A torus $\ov{T}$ is unique, up to conjugation, in~$O^{p'}(\ov{G}_\sigma)$ and
$\vert\ov{T}_\sigma\vert=(q-\varepsilon1)^n$, where $n$ is a rank
of~$\ov{G}$.
\end{lem}

\begin{proof}
Since for every maximal torus $T$ of $\ov{G}_\sigma$ the equality 
$\ov{G}_\sigma=TO^{p'}(\ov{G}_\sigma)$ holds, without lost of generality we may
assume that $G=\ov{G}_\sigma$. If $G$ is split, then the lemma can be easily
proven. In this case $\ov{T}$ is a maximal torus such that $\ov{T}_\sigma$ is a
Cartan subgroup of $\ov{G}_\sigma$ (i.~e. $\ov{T}$ is a maximal split
torus) and (a) is evident. Point (b) follows from
\cite[(10.1)]{GorLyoLocalStructure}. More over, from
\cite[(10.2)]{GorLyoLocalStructure} it follows that order of $\ov{T}_\sigma$ is
uniquely defined and is equal to  $(q-1)^n$, where $n$ is a rank of $\ov{G}$. By
\cite[F, \S6]{BorCarSeminar} we have that each element of order 
$r$ of $\ov{T}$ is contained in $\ov{G}_\sigma$. Now there exists $g\in \ov{G}$
such that $s^g\in \ov{T}$, hence  $s^g\in G$. In view of connectedness of the
centralizer of $s$, elements $s$ and $s^g$ are conjugate in $\ov{G}$ if and
only if they are conjugate in $G$, so $s$ and $s^g$ are conjugate in $G$, whence
(c). The information about classes of maximal tori, given in
\cite[G]{BorCarSeminar} and  \cite{Ca3ConjClassesWeyl}, implies that, up 
to conjugation by an element from~$G$, there exists a unique torus $\ov{T}$
such that~${\vert\ov{T}_\sigma\vert=(q-1)^n}$.

Assume that $O^{p'}(G)\simeq {}^2A_n(q)$. Then  $\ov{T}$ is a maximal torus
such that  $\vert \ov{T}_\sigma\vert=(q+1)^n$. Note that $\ov{T}_\sigma$ can be
obtained from a maximal split torus by twisting by the element $w_0\sigma$.
Direct calculations by using
\cite[Proposition~3.3.6]{CarFiniteGrpsLieTypeConjClassesCharacters} show that 
$N(\ov{G}_\sigma,\ov{T}_\sigma)/\ov{T}_\sigma$ is isomorphic to
$W(\ov{G})$, which in turn is isomorphic to $\Sym_{n+1}$. The uniqueness
follows from \cite[Proposition~8]{Ca2CentsSemisimpleClassical}. Point (b)
follows from \cite[(10.1)]{GorLyoLocalStructure}. As to point (c) we show first
that each element of order $r$ from $\ov{T}$ is in $G$. Assume that $t$ is an
element of order $r$ in $\ov{T}$ (recall that in this case $r$ divides $q+1$).
Let $\ov{H}$ be a $\sigma$-stab\-le maximal split torus of $\ov{G}$. The torus
$\ov{T}_\sigma$ is obtained from $\ov{H}$ by twisting by $w_0\sigma$, where 
$w_0\in W(\ov{G})$ is a unique element mapping all positive roots into
negatives and $\ov{T}_\sigma\simeq \ov{H}_{\sigma w_0}$. Let
$r_1,\ldots,r_n$ be a set of fundamental roots of $A_n$. Then $t$ as an element
of $\ov{H}$, can be written as $h_{r_1}(\lambda_1)\cdot\ldots\cdot
h_{r_n}(\lambda_n)$. Now for each $i$ we have $\sigma
w_0:h_{r_i}(\lambda)\mapsto h_{-r_i}(\lambda^q)=h_{r_i}(\lambda^{-q})$, i.~e.,
$t^{\sigma w_0}=t^{-q}$. Since $r$ divides $q+1$, we obtain that 
$t^{q+1}=e$, i.~e., $t=t^{-q}$. Hence $t^{\sigma w_0}=t$ and~${t\in
\ov{T}_\sigma}$. Now as in nontwisted case, there exists $g\in\ov{G}$ such that
$s^g\in \ov{T}$, therefore $s^g\in\ov{T}_\sigma$. In view of connectedness of
$C_{\ov{G}}(s)$, elements $s$ and $s^g$ are conjugate in~$G$.

For  $O^{p'}(G)={}^2D_{2n+1}(q)$ we take $\ov{T}$ to be equal to a unique (up
to conjugation in $G$) maximal torus, which has order $\vert
\ov{T}_\sigma\vert$ equals  $(q+1)^{2n+1}$ (uniqueness follows from
\cite[Proposition~10]{Ca2CentsSemisimpleClassical}), and for 
$O^{p'}(G)={}^2 E_6(q)$ we take $\ov{T}$ to be equal to a unique
(again up to conjugation in $G$) maximal torus, which has order $\vert
\ov{T}_\sigma\vert$ equals  $(q+1)^6$ (uniqueness follows from \cite[Table~1,
p.~128]{der1CentSemisimpleExcept}). As in the case of $G={}^2
A_n(q)$ it is easy to show that $\ov{T}$ satisfies (a), (b), and (c) of the
lemma.
\end{proof}

\begin{lem}\label{NormOfRegularElementIsNotCentrLie}
Let $G$ be a finite group of Lie type and $\ov{G}$, $\sigma$  are chosen so
that $O^{p'}(\ov{G}_\sigma)\leq G\leq \ov{G}_\sigma$. Let $s$ be a regular
semisimple element of odd prime order  $r$ of~$G$.

\noindent Then $N_G(C_{\ov{G}}(s))\not= C_G(s)$.
\end{lem}

\begin{proof}
In view of \cite[Proposition~2.10]{Hu2ConjClasses} we have that
$C_{\ov{G}}(s)/C_{\ov{G}}(s)^0$ is isomorphic to a subgroup of
$\Delta(\ov{G})$. Now,
if the root system $\Phi$ of $\ov{G}$ is not equal to either $A_n$, or $E_6$,
then $\vert \Delta(\Phi)\vert$ is a power of $2$. Since $\Delta(\ov{G})$ is
a quotient of $\Delta(\Phi(\ov{G}))$, then Lemma \ref{CentrOfInvolution}
implies that $C_{\ov{G}}(s)=C_{\ov{G}}(s)^0=\ov{T}$ is a maximal torus  and
$C_G(s)=C_{\ov{G}}(s)\cap G=T$. Since $N_G(T)\geq N(G,T)\not=T$ we obtain the
statement of the lemma in this case. Thus we may assume that either $\Phi= A_n$,
or~${\Phi=E_6}$.

Assume first that $\Phi=A_n$, i.~e., $O^{p'}(G)=A_n^\varepsilon(q)$, where
$\varepsilon\in\{+,-\}$.   Clearly $T=C_{\ov{G}}(s)^0\cap G$ is a normal
subgroup of
$C_G(s)$, hence $C_G(s)\leq N(G,T)$.  Assume that $N_G(C_{\ov{G}}(s))=C_G(s)$.
Then
$C_G(s)=N_{N(G,T)}(C_G(s))$ and $C_G(s)/T$ is a self-normalizing subgroup of
$N(G,T)/T$. As we noted above $C_G(s)/T$ is isomorphic to a subgroup of
$\Delta(A_n)$,
i.~e., it is cyclic. By Lemma \ref{CentrOfInvolution}, we also have that
$C_G(s)/T$
is an $r$-gro\-up, thus $C_G(s)/T=\langle x\rangle$ for some $r$-ele\-ment $x\in
N(G,T)/T$. Thus
$\langle x\rangle$ is a Carter subgroup of $N(G,T)/T$.  Now, in view
of \cite[Proposition~3.3.6]{CarFiniteGrpsLieTypeConjClassesCharacters},
we have that $N(G,T)/T\simeq C_{W(\ov{G})}(y)$
for
some $y\in W(\ov{G})\simeq \Sym_{n+1}$. Clearly $C_{C_{W(\ov{G})}(y)}(x)$
contains $y$, thus $y$
must be an $r$-ele\-ment, otherwise $N_{C_{W(\ov{G})}(y)}(\langle x\rangle)$
would
contain an element of order coprime to $r$, i.~e., $N_{C_{W(\ov{G})}(y)}(\langle
x\rangle)\not=\langle x\rangle$. A contradiction with the fact that $\langle
x\rangle$ is a Carter subgroup of $C_{W(\ov{G})}(y)$.

Now let $y=\tau_1\cdot\ldots$ be the decomposition of
$y$ into the product of independent cycles and $l_1,\ldots$ be
the lengths of $\tau_1,\ldots$ respectively. Assume first that
$m_1$ cycles has the same length $l_1$,  $m_2$ cycles has the
length $l_2$ etc. Let $m_0=n+1-(l_1m_1+\ldots+l_km_k)$. Then
$$C_{\Sym_{n+1}}(y)\simeq \left(Z_{l_1}\wr\Sym_{m_1}\right)\times\ldots\times
\left(Z_{l_k}\wr\Sym_{m_k}\right) \times
\Sym_{m_0},$$ where $Z_{l_i}$ is a cyclic group of order $l_i$. If
$m_j>1$ for some $j\ge0$, then there exists a normal subgroup $N$ of
$C_{\Sym_{n+1}}(y)$ such that $C_{\Sym_{n+1}}(y)/N\simeq
\Sym_{m_j}\not=\{e\}$. By Lemma  \ref{notlCarter}, Carter subgroup in a group
$S$
satisfying $\mathrm{Alt}_\ell\leq S\leq \Sym_\ell$ are conjugate for
all $\ell\ge1$. Thus $C_{W(\ov{G})}(y)$ and $N$ satisfy {\bfseries
(C)} and $\langle x\rangle$ is the unique, up to conjugation, Carter
subgroup of $C_{W(\ov{G})}(y)$. By Lemma \ref{HomImageOfCarter} we
obtain that $\langle x\rangle$ maps onto a Carter subgroup of
$C_{W(\ov{G})}(y)/N\simeq \Sym_{m_j}$. By Lemma \ref{notlCarter} only a Sylow 
$2$-sub\-gro\-up of $\Sym_{m_j}$ can be a Carter
subgroup of $\Sym_{m_j}$. A contradiction with the fact that $x$ is
an $r$-ele\-ment and $r$ is odd.

Thus we may assume that
$C_{W(\ov{G})}(y)\simeq\left(Z_{l_1}\times\ldots\times Z_{l_k}\right)$
and $l_i\not=l_j$ if $i\not=j$. From the known structure of maximal
tori and their normalizers in $A_n^\varepsilon(q)$ (see 
\cite[Propositions~7,8]{Ca2CentsSemisimpleClassical}, for example) we obtain the
structure of $T$ and
$N(G,T)$, which we explain by using matrices. Below a group
$\GL_n^\varepsilon(q)$ is isomorphic to $\GL_n(q)$ if $\varepsilon=+$ and is
isomorphic to $\GU_n(q)$ if $\varepsilon=-$. For the decomposition
$l_1+l_2+\ldots+l_k=n+1$ in $\GL_{n+1}^\varepsilon(q)$ consider a subgroup $L$,
consisting of block-diagonal matrices of view
$$
\left(
\begin{array}{cccc}
A_1&0&\ldots&0\\
0&A_1&\ldots&0\\
\ldots&\ldots&\ldots&\ldots\\
0&0&\ldots&A_k
\end{array}
\right),
$$
where $A_i\in \GL_{l_i}^\varepsilon(q)$. Then $L\simeq
\GL_{l_1}^\varepsilon(q)\times\ldots\times \GL_{l_k}^\varepsilon(q)$. Denote,
for
brevity, $\GL_{l_i}^\varepsilon(q)$ by $G_i$. In every group $G_i$ consider a
Singer cycle $T_i$.  $N_{G_i}(T_i)/T_i$ is known to be a cyclic group of
order $l_i$ and $N(G_i,T_i)=N_{G_i}(T_i)$. There exists a subgroup $Z$ of
$Z(\SL_{n+1}^\varepsilon(q))$ such that $O^{p'}(G)\simeq
\SL_{n+1}^\varepsilon(q)/Z$. Then $T\simeq
\left(\left(T_1\times\ldots\times T_k\right)\cap
\SL_{n+1}^\varepsilon(q)\right)/Z$ and $N(G,T)\simeq
\left(\left(N(G_1,T_1)\times\ldots\times N(G_k,T_k)\right)\cap
\SL_{n+1}^\varepsilon(q)\right)/Z$. Since for every Singer cycle $T_i$ the group
$N(G_i,T_i)/T_i$ is cyclic, we may assume that
$N(G,T)=C_G(s)$ and $T$ is a Singer cycle, i.~e., is a cyclic group of order
$\frac{q^{n+1}-(\varepsilon
1)^{n+1}}{q-(\varepsilon1)}$ and $n+1=r^k$ for some $k\ge1$
(the last equality holds, since $N(G,T)/T$ is an $r$-gro\-up). But
$q^{r^k}\equiv q \pmod r$, hence, $r$ divides
$q-(\varepsilon1)$. By Lemma \ref{InvolutionsAndTori} we obtain that $s$ is
in  $N(G,H)$, where $H$ is a maximal torus such that the factor group
$N(G,H)/H$ is isomorphic to $\Sym_{n+1}$ and $\vert H\vert=(q-\varepsilon1)^n$.
In particular, $H$ is not a Singer cycle. If  $s\in H$, this immediately
implies a contradiction with the choice of~$s$. If $s\not\in H$, then, since
the order of $s$ is prime, the intersection $\langle s\rangle\cap H$ is trivial.
Hence, under the natural homomorphism $N(G,H)\rightarrow N(G,H)/H\simeq
\Sym_{n+1}$ the element $s$ maps on an element of order  $r$. But in
$\Sym_{n+1}$ every element of odd order is conjugate to its inverse. Thus,
there exists a $2$-ele\-ment $z$ of $G$, which normalizes, but not centralizes
$\langle s\rangle$. Therefore,  $z\leq N_{\ov{G}}(C_{\ov{G}}(s))\leq
N_{\ov{G}}(C_{\ov{G}}(s)^0)$ and $\vert N(G,T)/T\vert$ is divisible by $2$,
that contradicts the above proven statement that $N(G,T)/T$ is an $r$-gro\-up.
This final contradiction finish the case~${\Phi(\ov{G})=A_n}$.

In the remaining case~${\Phi=E_6}$ it is easy to see, that for every
$y\in W(E_6)$, the group  $C_{W(E_6)}(y)$ does not contain Carter subgroup of
order $3$. Indeed, if  $C_{W(E_6)}(y)$ has a Carter subgroup of order $3$, then
it is generated by $y$. But it is known (and can be easily checked
by using \cite[Table~9]{Ca3ConjClassesWeyl}), that in $W(E_6)$ there is no
elements of order $3$, which
centralizer has order $3$. Since  $\vert
C_G(s)/T\vert$ divides $3$ and the group $C_G(s)/T$ is a Carter subgroup of
$C_{W(E_6)}(y)$ for some $y$, we get a contradiction.
\end{proof}

By using Lemma \ref{NormOfRegularElementIsNotCentrLie} we can obtain a similar
result for semilinear groups.

\begin{lem}\label{NormOfRegularElementIsNotCentr}
Let $\langle G,\zeta g\rangle$ be a finite semilinear group of Lie type and
$\ov{G}$, $\sigma$ are chosen so that
$O^{p'}(\ov{G}_\sigma)\leq G\leq \ov{G}_\sigma$. Let $s$ be a regular
semisimple element of odd order of~$G$. Then $N_{\langle G,\zeta
g\rangle}(C_{\langle G,\zeta g\rangle}(s))\not=
C_{\langle G,\zeta g\rangle}(s)$.
\end{lem}

\begin{proof}
Since $s$ is semisimple, there exists $\sigma$-stab\-le maximal torus $\ov{S}$
of
$\ov{G}$ containing $s$. Since $\ov{G}_\sigma
=O^{p'}(\ov{G}_\sigma)\ov{S}_\sigma$ we may assume that
$g\in \ov{S}_\sigma$, i.~e. elements $g$ and $s$ commutes. If
$C_{\langle G,\zeta g\rangle}(s) G\not=\langle G,\zeta g\rangle$, then we can
substitute $\langle G,\zeta g\rangle$ by $C_{\langle
G,\zeta g\rangle}(s) G$ and prove the lemma for this group. More over, if
$C_{\langle G,\zeta g\rangle}(s)=C_G(s)$, then the lemma follows from Lemma
\ref{NormOfRegularElementIsNotCentrLie}, so we may assume that $\zeta$
centralizes $s$. If either  $G$ is not twisted, or $\vert\zeta\vert$ is
odd, then by  \cite[Proposition~2.5.17]{GorLySol} it follows that we may assume
$\sigma=\bar\zeta^k$ for some $k>0$. By Lemma
\ref{NormOfRegularElementIsNotCentrLie} there exists an element of $N_{G_{\zeta
g}}(C_{\ov{G}}(s))$, not contained in $C_{G_{\zeta g}}(s)$, and the lemma
follows.

Assume that $G$ is twisted and $\vert\zeta\vert$ is even. Then
$\sigma=\bar\gamma\bar{\varphi}^k$, $\bar\zeta=\bar{\varphi}^\ell$, where $k$
divides $\ell$. Therefore $s$ is in $\ov{G}_{\bar\gamma}$. Depending on the
root system  $\Phi(\ov{G})$, we obtain that $\ov{G}_{\bar\gamma}$ is isomorphic
to a simple algebraic group with root system equal to $B_m$ (for some $m>1$),
$C_m$ (for some $m>2$), or $F_4$. By Lemma  \ref{omnibus} the element
$s$ is conjugate with its inverse under
$O^{p'}((\ov{G}_{\bar\gamma})_{\sigma\bar\zeta g})\leq G_{\zeta g}$, so
$N_{\langle G,\zeta g\rangle}(C_{\langle G,\zeta g\rangle}(s))\not= C_{\langle
G,\zeta g\rangle}(s)$.
\end{proof}

\begin{lem}\label{Syl2InCentrOfFieldAut}
Let $G$ be a finite  group of Lie type over a field of odd
characteristic $p$. Assume that $\ov{G}$
and $\sigma$ are chosen so that $O^{p'}(\overline{G}_\sigma)\leq
G\leq\overline{G}_\sigma$. Let $\psi$ be a field
automorphism of odd order of $O^{p'}(\overline{G}_\sigma)$.

Then $\psi$ centralizes a Sylow $2$-sub\-gro\-up of $G$, and there exists a
$\psi$-stab\-le Cartan subgroup $H$  such that $\psi$ centralizes a Sylow
$2$-sub\-gro\-up of $H$. Moreover, if
$G\not\simeq {}^2G_2(3^{2n+1})$, ${}^3D_4(q^3)$, ${}^2D_{2n}(q^2)$, then there
exists a $\psi$-stab\-le torus $T$ of $G$ such that $\psi$ centralizes a Sylow
$2$-sub\-gro\-up of $T$ and  the factor group  $N(G,T)/T$ is isomorphic to
$N_{\ov{G}}(\ov{T})/\ov{T}$.
\end{lem}

\begin{proof}
Clearly we need to prove the lemma only for the case $G=\ov{G}_\sigma$.  Assume
that
$\vert \psi\vert=k$. Let $\F_q$ be the base field of $G$. Then $q=p^\alpha$ and
$\alpha=k\cdot m$. Now $\vert G\vert$ can be written as
$\vert
G\vert=q^N(q^{m_1}+\varepsilon_11)\cdot\ldots\cdot(q^{m_n} +\varepsilon_n 1)$
for
some $N$, where $n$ is the rank of $G$, $\varepsilon_i=\pm$ (see
\cite[Theorems~9.4.10 and~14.3.1]{CarSimpleGrpsLieType}). Similarly we have that
$\vert
G_\psi\vert=(p^m)^N((p^m)^{m_1}+\varepsilon_11)\cdot\ldots\cdot
((p^m)^{m_n}+\varepsilon_n1)$, i.~e., $\vert G\vert_2=\vert G_\psi\vert_2$ and a
Sylow $2$-sub\-gro\-up of $G_\psi$ is a Sylow $2$-sub\-gro\-up of $G$. By
\cite[Proposition~2.5.17]{GorLySol} there exists an automorphism $\psi_1$ of
$\ov{G}$ such that $\sigma=\psi_1^k$ and $\psi$ coincides with the restriction
of $\psi_1$ on $\ov{G}_\sigma$. Note that $\psi_1$, in general, is not equal to
$\bar\psi$ defined above. Consider a maximal split torus $\ov{H}_{\psi_1}$ of
$\ov{G}_{\psi_1}$. Then $H=\ov{H}_\sigma$ is a $\psi$-stab\-le Cartan subgroup
of
$G$. Since $\vert
H\vert=(q^{k_1}+\varepsilon1)\cdot\ldots\cdot(q^{k_l}+\varepsilon_l1)$, where
$\varepsilon_i=\pm$, the equality $\vert H\vert_2=\vert H_\psi\vert_2$ can be
proven in the same way.

Now assume that $G\not\simeq {}^2G_2(3^{2n+1})$, ${}^3D_4(q^3)$,
${}^2D_{2n}(q^2)$. By Lemma \ref{InvolutionsAndTori}, there exists a maximal
torus $T$ of $G_\psi$
such that $N(G_\psi,T)/T\simeq N_{\ov{G}}(\ov{T})/\ov{T}$ and $\vert
T_\psi\vert=(p^m-\varepsilon1)^n$. Since $\vert\psi\vert$ is odd and
$\ov{T}_{\psi_1}$  is obtained from a maximal splittorus $\ov{H}$ by twisting
by an element $w_0$, then $\ov{T}_\sigma$ is also obtained from a maximal
split torus $\ov{H}$ by twisting by element $w_0$ (see proof of Lemma 
\ref{Esyl2InnDiagExtension}). Therefore 
$\vert\ov{T}_\sigma\vert=(q-\varepsilon1)^n$, $\vert
\ov{T}_{\psi_1}\vert=(p^m-\varepsilon1)^n$, hence
$\vert\ov{T}_\sigma\vert_2=\vert T\vert_2=\vert T_\psi\vert_2$.
\end{proof}

\begin{lem}\label{ConjAutomorphisms}
{\em \cite[(7-2)]{GorLyoLocalStructure}} Let $\ov{G}$ be a connected simple
linear algebraic group over a field of characteristic $p$, $\sigma$ be a
Frobenius map of $\ov{G}$ and $G=\ov{G}_\sigma$ be a finite group of Lie
type. Let $\varphi$ be a field or a graph-fi\-eld automorphism of $G$ (if
$\varphi$ is graph-fi\-eld, then corresponding graph automorphism has order $2$)
and let $\varphi'$ be an element of
$(G\leftthreetimes\langle\varphi\rangle)\setminus G$ such
that~${\vert\varphi'\vert=\vert\varphi\vert}$ and
$G\leftthreetimes\langle\varphi\rangle=G\leftthreetimes\langle\varphi'\rangle$.

\noindent Then there exists $g\in G$ such that
$\langle\varphi\rangle^g=\langle\varphi'\rangle$. In particular, if
$G/O^{p'}(G)$ is a $2$-gro\-up and $\varphi$ is of odd order, then $g$ can be
chosen in~$O^{p'}(G)$.
\end{lem}

A particular case of the following lemma is proven
in~\cite[Theorem~A]{FeZuRealConjClasses}).

\begin{lem}\label{ConjInverseInGraph}
Let $G$ be a finite adjoint split group of Lie type,
$\overline{G}$, $\sigma$
are chosen so that $O^{p'}(\overline{G}_\sigma)\leq G\leq \overline{G}_\sigma$.
Assume that $\tau$ is  a graph automorphism of order $2$ of~$O^{p'}(G)$.

\noindent Then every semisimple element $s\in G$ is conjugate to its inverse
under~$\langle
O^{p'}(\overline{G}_\sigma), \tau a\rangle$, where $a$ is an element
of~$\ov{G}_\sigma$.
\end{lem}

\begin{proof}
If $\Phi(\ov{G})$ is not of type $A_n,D_{2n+1},E_6$,
then the
lemma follows from Lemma \ref{omnibus}, thus we need to consider groups of type
$A_n,D_{2n+1},E_6$.  Denote by $\bar\tau$ the graph automorphism of $\ov{G}$
such
that $\bar\tau\vert_G=\tau$. Let $\ov{T}$ be a maximal $\sigma$-stab\-le torus
of
$\ov{G}$ such that $\ov{T}_\sigma\cap G$ is a Cartan subgroup of $G$. Let
$r_1,\ldots, r_n$ be fundamental roots of $\Phi(\ov{G})$ and $\rho$ be the
symmetry
corresponding to $\bar\tau$. Denote $r_i^\rho$ by $\bar{r}_i$. Then
$\ov{T}=\langle
h_{r_i}(t_i)\mid 1\le i \le n,\ t_i\not=0\rangle$ and
$h_{r_i}(t_i)^{\bar\tau}=h_{\bar{r}_i}(t_i)$. Denote by $W$ the Weyl group of
$\ov{G}$. Let $w_0$ be the unique element of $W$ mapping all positive roots onto
negative roots and let $n_0$ be its preimage in $N_{\ov{G}}(\ov{T})$ under the
natural homomorphism $N_{\ov{G}}(\ov{T})\rightarrow
N_{\ov{G}}(\ov{T})/\ov{T}\simeq
W$. Since $\sigma$ acts trivially on $W=N(G,T)/T$ (see Lemma
\ref{InvolutionsAndTori}),  we can take $n_0\in G$, i.~e.,
$n_0^\sigma=n_0$. Then for all $r_i$ and $t$ we have that
$$h_{r_i}(t)^{n_0\bar\tau}=h_{r_i^{w_0\rho}}(t)=h_{-r_i}(t)=h_{r_i}(t^{-1}).$$
Thus $x^{n_0\bar\tau}=x^{-1}$ for all $x\in \ov{T}$.

Now let $s$ be a semisimple element of $G$. Then there exists a maximal
$\sigma$-stab\-le torus $\ov{S}$ of $\ov{G}$ containing $s$. Since all maximal
tori of
$\ov{G}$ are conjugate, we have that there exists $g\in \ov{G}$ such that
$\ov{S}^g=\ov{T}$. Since $\ov{G}_\sigma=O^{p'}(\ov{G}_\sigma)\ov{T}_\sigma$,
then we may assume that $a\in\ov{T}_\sigma$. Therefore $s^{gn_0\bar\tau a
g^{-1}}=s^{-1}$. Since
$n_0^\sigma=n_0$ and $\bar\tau^\sigma=\bar\tau$ we have that $(gn_0\bar\tau
ag^{-1})^\sigma=g^\sigma n_0\bar\tau a(g^{-1})^\sigma$. Moreover, since $\ov{S}$
is
$\sigma$-stab\-le, then for every $x\in\ov{S}$ we have that $x^{gn_0\bar\tau
ag^{-1}}=x^{g^\sigma n_0\bar\tau a(g^{-1})^\sigma}=x^{-1}$, i.~e., $gn_0\bar\tau
ag^{-1}\ov{S}=g^\sigma n_0\bar\tau a(g^{-1})^\sigma\ov{S}$. In particular, there
exists $t\in\ov{S}$ such that  $gn_0\bar\tau ag^{-1}t=g^\sigma n_0\bar\tau
a(g^{-1})^\sigma$. In view of Lang-Steinberg Theorem (Lemma
\ref{LangSteinbergTheorem}) there
exists $y\in\ov{S}$ such that $t=y\cdot (y^{-1})^\sigma$. Therefore, $g
n_0\bar\tau
ag^{-1}y=(g n_0\bar\tau ag^{-1}y)^\sigma$, i.~e., $gn_0\tau ag^{-1}y\in
\ov{G}_\sigma\leftthreetimes \langle\tau\rangle$, and $s^{gn_0\tau
ag^{-1}y}=s^{-1}$.
Since $O^{p'}(\ov{G}_\sigma)\ov{S}_\sigma=\ov{G}_\sigma$, and $\ov{S}_\sigma$ is
Abelian, we may find $z\in \ov{S}_\sigma$ such that $gn_0\tau ag^{-1}yz\in
\langle O^{p'}(\ov{G}_\sigma),\tau a\rangle$.
\end{proof}

\subsection{Carter subgroups of special type}

In this subsection we consider problems of structure and existence of Carter
subgroups in semilinear groups, containing a Sylow $2$-sub\-gro\-up or is
contained in the normalizer of a Borel subgroup.

\begin{lem}\label{Esyl2InnDiagExtension}
Let $G$ be a finite group of Lie type over a field of odd characteristic and
$\ov{G}$, $\sigma$ are chosen so that
$O^{p'}(\ov{G}_\sigma)\leq G\leq \ov{G}_\sigma$. If  $G$ satisfies
{\em\bfseries(ESyl2)}, then every subgroup $L$ with $G\leq L\leq
O^{p'}(\ov{G}_\sigma)$ satisfies~{\em\bfseries(ESyl2)}.
\end{lem}

\begin{proof}
Let $Q$ be a Sylow $2$-sub\-gro\-up of $\ov{G}_\sigma$ and
$Q^0=O^{p'}(\ov{G}_\sigma)\cap Q$ be a Sylow $2$-sub\-gro\-up of
$O^{p'}(\ov{G}_\sigma)$.  If $N_{\ov{G}_\sigma}(Q^0)=QC_{\ov{G}_\sigma}(Q)$,
then the statement of the lemma is clearly true. In view of
\cite[Theorem~1]{KondNormalizers}, for a classical group $\ov{G}_\sigma$ the
equality
$N_{\ov{G}_\sigma}(Q^0)=QC_{\ov{G}_\sigma}(Q)$ can fail to be true only if the
root system of $\ov{G}$ has type $A_1$ or $C_n$. If the root system of $\ov{G}$
has type $A_1$ or $C_n$, then $\vert
\ov{G}_\sigma:O^{p'}(\ov{G}_\sigma)\vert=2$ and the statement of the lemma
follows from Lemma~\ref{InhBy2-ext}.

Assume now that $G$ is a group of exceptional type. If
$\ov{G}_\sigma=O^{p'}(\ov{G}_\sigma)$, then the statement of the lemma is
clearly true. The equality $N_{\ov{G}_\sigma}(Q^0)=QC_{\ov{G}_\sigma}(Q)$ might
fail to be true only if the root system of $\ov{G}$ has type $E_6$ or $E_7$. If
the root system of $\ov{G}$ has type $E_7$, then $\vert
\ov{G}_\sigma:O^{p'}(\ov{G}_\sigma)\vert=2$ and the statement of the lemma
follows from Lemma~\ref{InhBy2-ext}.

Assume that the root system of $\ov{G}$ has type $E_6$. Then either
$\ov{G}_\sigma=O^{p'}(\ov{G}_\sigma)$ or $\vert
\ov{G}_\sigma:O^{p'}(\ov{G}_\sigma)\vert=3$. In the first case we have nothing
to prove, so assume that $\vert \ov{G}_\sigma:O^{p'}(\ov{G}_\sigma)\vert=3$.
Since the group $G$ coincides either with $\ov{G}_\sigma$, or with
$O^{p'}(\ov{G}_\sigma)$, and since in case $G=\ov{G}_\sigma$ there is nothing
to prove, we may assume that $G=O^{p'}(\ov{G}_\sigma)$. By
\cite[Theorem~4.10.2]{GorLySol} there exists a maximal torus $T$ of
$\ov{G}_\sigma$ such that $Q$ is contained in $N(\ov{G}_\sigma,T)$. Since
$\vert \ov{G}_\sigma:G\vert=3$, then $Q=Q^0\leq N(G,T\cap G)$. By
\cite[Theorem~6]{KonMaz} the equality
$N_G(Q)=Q\times R^0$ holds, where $R^0\leq T$ is a cyclic group of odd order.
Now
since $\ov{G}_\sigma=TG$, then $N_{\ov{G}_\sigma}(Q)=\langle
N_T(Q),N_G(Q)\rangle$. Indeed, $N(G,T\cap G)/(T\cap G)\simeq N(G,T)/T$.
Hence, a Sylow $2$-sub\-gro\-up $QT/T$ of $N(G,T)/T$ coincides with its
normalizer.
Since the factor group $\ov{G}_\sigma/G$ is cyclic of order $3$, then
$N_{\ov{G}_\sigma}(Q)=\langle tg, N_G(Q)\rangle$, where $t\in T$
and $g\in G$. Moreover, since $\vert \ov{G}_\sigma:G\vert=3$, we may assume
that $tg$ is an element of order $3^k$ for some $k>0$. Since
$t\in T\leq N(\ov{G}_\sigma, T)$, then $Q^t\leq N(G,T\cap G)$. So there exists
an element  $g_1\in N(G,T\cap G)$ such that $Q^t=Q^{g_1^{-1}}$. Therefore we
may assume that $tg=tg_1\in N(\ov{G}_\sigma,T)$. Under the natural epimorphism
 $\pi:N(\ov{G}_\sigma,T)\rightarrow N(\ov{G}_\sigma,T)/T$ the image of
$N_{N(\ov{G}_\sigma,T)}(Q)$ coincides with $Q$. Hence, $(tg)^\pi=e$,
so $tg\in T$. Thus each element of odd order of
$\ov{G}_\sigma$
normalizing $Q$ lies in $T$. Since $T$ is a torus, then $T$ is Abelian, hence
the set of elements of odd order of  $N_{\ov{G}_\sigma}(Q)$ forms a normal
subgroup $R\leq T$. Therefore $N_{\ov{G}_\sigma}(Q)=Q\times R$, i.~e.,
$\ov{G}_\sigma$ satisfies~{\bfseries(ESyl2)}.
\end{proof}

The following lemma is immediate from~\cite[теорема~1]{KondNormalizers}.

\begin{lem}\label{Norm2Syl}
Let $O^{p'}(\overline{G}_\sigma)= G$ be a canonical
finite group of Lie type and $\overline{G}$ is either of type $A_1$ or  of
type $C_n$, $p$ is odd, $q=p^\alpha$ is the order of the base field of
$G$. Then $G$ satisfies {\em\bfseries(ESyl2)} if and only
if~$q\equiv\pm1\pmod{8}$.
\end{lem}

Note that Lemma \ref{Esyl2InnDiagExtension}  together with
\cite[Theorem~1]{KondNormalizers} and \cite[Theorem~6]{KonMaz} implies that
every  group of Lie type over
a field of odd characteristic, distinct from a Ree group and groups from Lemma
\ref{Norm2Syl}, satisfies~{\bfseries(ESyl2)}.

\begin{lem}\label{ESyl2InhFieldGraph}
Let $G$ be a finite adjoint group of Lie type over a field of odd
characteristic,
$G\not\simeq {}^3D_4(q^3)$, and $\ov{G}$, $\sigma$ are chosen so that
$O^{p'}(\ov{G}_\sigma)\leq G\leq \ov{G}_\sigma$. Let $A$ be a
subgroup of $\Aut (O^{p'}(\ov{G}_\sigma))$ such that $A\cap \ov{G}_\sigma =G$.
If $O^{p'}(G)\simeq
D_4(q)$, assume also that $A$ is contained in the group generated by
inner-diagonal,
field automorphisms and a graph automorphism of order $2$. 

\noindent Then $A$ satisfies
{\em\bfseries(ESyl2)} if and only if $G$
satisfies~{\em\bfseries(ESyl2)}.
\end{lem}

\begin{proof}
Assume that $G$ satisfies {\bfseries(ESyl2)}. In the conditions of the lemma we
have that the factor group $A/G$ is Abelian, so $A/G=\ov{A}_1\times\ov{A}_2$,
where $\ov{A}_1$ is a Hall ${2'}$-sub\-gro\-up of $A/G$ and $\ov{A}_2$ is a
Sylow
$2$-sub\-gro\-up of $A/G$. Denote by $A_1$ the complete preimage of $\ov{A}_1$
in
$A$. If $A_1$ satisfies {\bfseries(ESyl2)}, then by Lemma \ref{InhBy2-ext} $A$
satisfies {\bfseries(ESyl2)} as well. Thus we may assume that the order $\vert
A/G\vert$ is odd. Since we are assuming that a graph automorphism of order $3$
is not contained in $A$, then $A/G$ is cyclic, hence $A=\langle G,\psi
g\rangle$, where $\psi$ is a field automorphism of odd order and
$g\in\ov{G}_\sigma$. Since $\vert A:G\vert=\vert \psi\vert$ is odd, we may
assume that $\vert\psi g\vert$ is also odd. By Lemma
\ref{Syl2InCentrOfFieldAut}, $\psi$ centralizes a Sylow $2$-sub\-gro\-up of
$\ov{G}_\sigma$, therefore $g$ is of odd order. Now the quotient
$\ov{G}_\sigma/G$ is Abelian and can be written as $\ov{L}\times\ov{Q}$, where
$\ov{L}$ is a Hall $2'$-sub\-gro\-up of $\ov{G}_\sigma/G$ and $\ov{Q}$ is a
Sylow
$2$-sub\-gro\-up of $\ov{G}_\sigma/G$. Let $L$ be the complete preimage of
$\ov{L}$ in $\ov{G}_\sigma$ under the natural homomorphism. Then $g\in L$.
Consider $L\leftthreetimes\langle\psi\rangle\geq A$. By construction, $\vert
L\leftthreetimes\langle\psi\rangle:A\vert=\vert L:G\vert$ is odd. By
Lemma \ref{Esyl2InnDiagExtension} the group $L$ satisfies {\bfseries(ESyl2)}.
By Lemma \ref{Syl2InCentrOfFieldAut} the field automorphism $\psi$ centralizes
a Sylow $2$-sub\-gro\-up $Q$ of $L$, Thus
$$N_{L\leftthreetimes\langle\psi\rangle}
(Q)=N_L(Q)\times\langle\psi\rangle=QC_L(Q)\times\langle
\psi\rangle=QC_{L\leftthreetimes\langle\psi\rangle}(Q),$$ i.~e., the group
$L\leftthreetimes \langle\psi\rangle$ satisfies {\bfseries(ESyl2)}. Since
$\vert L\leftthreetimes \langle\psi\rangle:A\vert$ is odd, then $A$
satisfies~{\bfseries(ESyl2)}.

Now assume that $A$ satisfies {\bfseries(ESyl2)}. If $G$ does not satisfies
{\bfseries(ESyl2)}, then  \cite[Theorem~1]{KondNormalizers} and
\cite[Theorem~6]{KonMaz} imply that the root system of $\ov{G}$ either has type
$A_1$, or has type
$C_n$. In both cases the factor group
$\Aut(O^{p'}(\ov{G}_\sigma)/\ov{G}_\sigma)$ is cyclic and is generated by a
field automorphism $\varphi$. Further, from  \cite[Theorem~1]{KondNormalizers}
it follows
that the order of the base field (that is equal to the field of definition in
this case, since $G$ is not twisted) is equal to $q=p^t$ and
$q\equiv\pm3\pmod8$. Therefore $t$ is odd and, so
$\vert\Aut(O^{p'}(\ov{G}_\sigma))/\ov{G}_\sigma\vert$ is odd. Thus  $\vert
A:G\vert$ is odd, hence $G$ satisfies~{\bfseries(ESyl2)}.
\end{proof}

\begin{lem}\label{CartBorel}
Let $\langle G,\zeta g\rangle$ be a finite semilinear group of Lie type over
a field of characteristic $p$ (we do not exclude the case $\langle G,\zeta
g\rangle=G$) and $G$ is of adjoint type (recall that $g\in\ov{G}_\sigma$, but
not necessary $g\in G$). Assume that $B=U\leftthreetimes H$, where $H$ is a
Cartan subgroup of $G$, is a
$\zeta g$-in\-va\-ri\-ant Borel subgroup of $G$ and
$\langle B, \zeta g\rangle$ contains a Carter subgroup $K$ of
$\langle G,\zeta g\rangle$. Assume that $K\cap U\not=\{e\}$. Then one of
the following statements holds:
\begin{itemize}
\item[{\em (a)}] either $\langle G,\zeta g\rangle=\langle{}^2A_2(2^{2t}), \zeta
g\rangle$, or
$\langle G,\zeta g\rangle=\widehat{{}^2A_2(2^{2t})}\leftthreetimes
\langle\zeta\rangle$; the
order $\vert \zeta\vert=t$ is odd and is not divisible by $3$, $C_G(\zeta)\simeq
\widehat{{}^2A_2(2^2)}$, $K\cap G$ is Abelian and has order~$2\cdot 3$;
\item[{\em (b)}] $\langle G,\zeta g\rangle=\langle {}^2A_2(2^{2t})\zeta
g\rangle$, the order $\vert\zeta\vert=t$ is odd, $C_G(\zeta)\simeq
{}^2A_2(2^2)$, the subgroup $K\cap G$ is a Sylow $2$-sub\-gro\-up
of~$G_\zeta$;
\item[{\em (c)}] either $\langle G,\zeta g\rangle=\langle A_2(2^{2t}),\zeta
g\rangle$, or $\langle G,\zeta g\rangle=\widehat{A_2(2^{2t})}\leftthreetimes
\langle\zeta\rangle$, $\zeta$ is a graph-field automorphism of order $2t$, $t$
is not divisible by $3$, and $C_G(\zeta)\simeq \widehat{{}^2A_2(2^2)}$, the
subgroup $K\cap G$ is Abelian and has order $2^{\vert\zeta_{2'}\vert}\cdot
3$;
\item[{\em (d)}] $\langle G,\zeta g\rangle=\langle A_2(2^{2t}),\zeta g\rangle$,
$\zeta$ is a graph-field automorphism and $C_G(\zeta)\simeq {}^2A_2(2^2)$, the
subgroup $K\cap G$ is a Sylow $2$-sub\-gro\-up of~$G_{\zeta_{2'}}$;
\item[{\em (e)}] $G$  is defined over
$\F_{2^t}$, $\langle
G,\zeta g\rangle =G\leftthreetimes\langle\zeta g\rangle$, $\zeta$ is either a
field automorphism of
order $t$ of $O^{2'}(G)$, if $O^{2'}(G)$ is split, or a graph automorphism of
order $t$, if $O^{2'}(G)$ is twisted,  and, up to conjugation in $G$,
$K=Q\leftthreetimes\langle\zeta g\rangle$, where $Q$ is a Sylow $2$-sub\-gro\-up
of~$G_{(\zeta g)_{2'}}$;
\item[{\em (f)}] $G$ is split and defined over $\F_{2^t}$, $\langle G,\zeta
g\rangle=G\leftthreetimes\langle \zeta g\rangle$, $\zeta$ is a product of a
field automorphism of odd order $t$ of $O^{2'}(G)$ and a graph automorphism of
order $2$, $\zeta$ and $\zeta g$ are conjugate under $\ov{G}_\sigma$, and, up
to a conjugation in $G$, $K=Q\leftthreetimes \langle \zeta g\rangle$, where $Q$
is a Sylow $2$-sub\-gro\-up of~$G_{(\zeta g)_{2'}}$;
\item[{\em (g)}] $G/Z(G)\simeq \P SL_2(3^{t})$, the order
$\vert\zeta\vert=t$ is odd (hence $\zeta\in\langle G,\zeta g\rangle$),
and $K$
contains a Sylow $3$-sub\-gro\-up of~$G_{\zeta_{3'}}$;
\item[{\em (h)}] $\langle G,\zeta
g\rangle={}^2G_2(3^{2n+1})\leftthreetimes\langle\zeta\rangle$,
$\vert\zeta\vert=2n+1$, $K\cap {}^2G_2(3^{2n+1})=Q\times P$, where $Q$ is of
order $2$ and~${\vert P\vert =3^{\vert\zeta\vert_3}}$.
\end{itemize}
\end{lem}

Note that in all points (a)--(h) of the lemma Carter subgroups, having given
structure, do exist. The existence of Carter subgroups in points
(a) and (c) follows from the existence of a Carter subgroup of order $6$ in
$\P \GU_3(2)$ (see \cite{DTZCartSbgrpsClassGrps}). The existence of Carter
subgroups in points
(b), (d)-(f) follows from the fact that a Sylow $2$-sub\-gro\-up in a group of
Lie
type defined over a field of order $2$, coincides with its normalizer. The
existence of Carter subgroups in point (g) follows from the fact that a Sylow
$3$-sub\-gro\-up of $\P \SL_2(3)$ coincides with its normalizer. The existence
of a
Carter subgroup, satisfying point (h) of the lemma, follows from the existence
of a Carter subgroup $K$ of order $6$ in a (non simple) group ${}^2G_2(3)$. The
existence of a Carter subgroup $K$ of order $6$ in ${}^2G_2(3)$ follows from
the results given in
\cite{levnuzReeGrps} and~\cite{warReeGrps}.

\begin{proof}
If $G$ is one of the  groups $A_1(q)$, $G_2(q)$, $F _4(q)$,
${}^2B_2(2^{2n+1})$, or ${}^2F_4(2^{2n+1})$, then the lemma follows from 
Table \ref{AlmSimpleNotCounterEx}. If $\langle
G,\zeta g\rangle=G$, then the lemma  follows from the results of section 3 and
Theorem \ref{CarterInClassicalGroups}. So we may assume that
$\langle G,\zeta g\rangle\not=G$,
i.~e., that $\zeta$ is a nontrivial field,  graph-field, or graph automorphism.
If $\Phi(\ov{G})=C_n$, the lemma follows from
Theorem~\ref{sympcarter} below, that does not use Lemma \ref{CartBorel}, so we
assume
that~${\Phi(\ov{G})\not=C_n}$. If $\Phi(\ov{G})=D_4$ and either a graph-field
automorphism $\zeta$ is a product of a field automorphism and a graph
automorphism of order $3$, or $G\simeq {}^3D_4(q^3)$, then the
lemma follows from
Theorem~\ref{CarterTriality} below, that does not use Lemma \ref{CartBorel}, so
we assume that $\langle  G,\zeta g\rangle$ is contained in the group $A_1$
defined in Theorem
\ref{CarterTriality}, and $G\not\simeq{}^3D_4(q^3)$. Since we shall use Lemma
\ref{CartBorel} in the proof of Theorem
\ref{CarterSemilinear}, after Theorems \ref{sympcarter} and
\ref{CarterTriality}, it
is possible to make such additional assumptions.

Assume that $q$ is odd and $\Phi(\ov{G})$ is one of the following types:
$A_n$ $(n\ge2)$, $D_n$ $(n\ge 4)$, $B_n$ $(n\ge 3)$, $E_6$, $E_7$ or
$E_8$. By Lemma \ref{HomImageOfCarter} we have that $KU/U$ is a Carter subgroup
of $\langle B,\zeta g\rangle/U\simeq\langle H, \zeta g\rangle$.
Since $\ov{G}_\sigma=G\ov{H}_\sigma$, where $\ov{H}$ is a maximal split torus of
$\ov{G}$ and $\ov{H}_\sigma\cap G=H$, then we may assume that
$g\in\ov{H}_\sigma$, in particular $g$ centralizes $H$. So  $H_\zeta\leq
Z(\langle H, \zeta g\rangle)$, and we obtain, up to
conjugation in $B$, that $H_\zeta\leq K$. By Lemma \ref{Syl2InCentrOfFieldAut},
the automorphism $\zeta_{2'}$ centralizes a Sylow $2$-sub\-gro\-up $Q$ of $H$.
Thus, each element of odd order of $\langle H,\zeta g\rangle$ centralizes $Q$
and Lemma \ref{CritSyl2Carter} implies, that, up to conjugation in $B$, the
inclusion $Q\leq K$ holds. By Lemma \ref{centUH} it follows that
$C_U(Q)=\{e\}$, a
contradiction with the fact that  ${K\cap U}$ is nontrivial.

Assume that $G\simeq{}^2G_2(3^{2n+1})$ and $\langle
G,\zeta g\rangle=G\leftthreetimes\langle\zeta\rangle$ (in this case
$O^{p'}(\ov{G}_\sigma)=\ov{G}_\sigma$). Again by Lemma \ref{HomImageOfCarter}
we have that $KU/U$ is a Carter subgroup of
$(B\leftthreetimes\langle
\zeta\rangle)/U\simeq H\leftthreetimes\langle\zeta\rangle$. By Lemma 
\ref{omnibus} every semisimple element of $G$ is conjugate to its
inverse. Since non-Abelian composition factors of every semisimple element of
$G$ can be isomorphic only to groups $A_1(q)$, by Table
\ref{AlmSimpleNotCounterEx} it follows
that the centralizer of every semisimple element of $G$ satisfies condition
{\bfseries(C)}. So Lemma \ref{power} implies that $KU/U\cap B/U$ is a
$2$-gro\-up. On the other hand, $\vert H\vert_2=2$ and $KU/U\geq Z(B/U)\geq
H_\zeta$, hence $\vert H_\zeta\vert=2$ and $\vert\zeta\vert=2n+1$. Thus $K
\cap G=(K\cap U)\times\langle t\rangle$,  where $t$ is an involution. It follows
that $K\cap U=C_G(t)\cap G_{\zeta_{3'}}$. Now the structure results from
\cite[Theorem~1]{levnuzReeGrps} and~\cite{warReeGrps} imply point (h) of the
lemma.

Assume now that $q=2^t$. Assume first that $\Phi(\ov{G})$ has one of the
types $A_n$ $(n\ge2)$, $D_n$ $(n\ge 4)$, $B_n$
$(n\ge 3)$, $E_6$, $E_7$ or $E_8$, $G$ is split, and $\zeta$ is a field
automorphism. Like above we obtain that $H_\zeta\leq
K$,  and $O^{2'}(G_\zeta)$ is a split group of Lie type with definition
field of order $q=2^{t/\vert\zeta\vert}$. By Hartley-Shute Lemma
\ref{Hartley-Shute}, for every $r\in\Phi(\ov{G})$ and for every $s\in
GF(2^{t/\vert\zeta\vert})^\ast$ there exists $h(\chi)\in H_\zeta \cap
O^{2'}(G_\zeta)$ such that $\chi(r)=s$. The same arguments as in Lemma
\ref{centUH} imply that if $\frac{t}{\vert\zeta\vert}\not=1$, inequality $K\cap
U\leq C_U(H_\zeta)=\{e\}$ holds, a contradiction. So $\vert\zeta\vert=t$ and
$H_\zeta=\{e\}$. Since $g$ can be chosen in $\ov{H}_\sigma$ and $\langle \zeta
g\rangle\cap \ov{G}_\sigma\leq \langle\zeta g\rangle\cap \ov{H}_\sigma\leq
H_\zeta=\{e\}$, then $\langle\zeta g\rangle\cap \ov{G}_\sigma=\{e\}$. By Lemma
\ref{ConjAutomorphisms} elements $\zeta g$ and $\zeta$ are conjugate under
$\ov{G}_\sigma$, and point (e) of the lemma follows.

Now assume that $\Phi(\ov{G})$ is  of type $A_n$ $(n\ge3)$, $D_n$
$(n\ge 4)$, or $E_6$;  and either $\zeta$ is a
graph-field automorphism and $G$ is split, or $ G$ is twisted. Let $\rho$
be the
symmetry of the Dynkin diagram of $\Phi(\ov{G})$ corresponding to $\gamma$
(recall that $\zeta=\gamma^\varepsilon\varphi^\ell$ by definition), and
$\bar{r}$ denotes $r^\rho$ for $r\in\Phi(\ov{G})$. Like above it is possible to
prove that, up to conjugation, $H_\zeta\leq K$. If $\vert\zeta\vert=2t$,
then $H_\zeta\not=\{e\}$, then by Hartley-Shute Lemma \ref{Hartley-Shute} we
obtain that $C_U(H_\zeta)=\{e\}$ that contradicts the condition $K\cap
U\not=\{e\}$. If $H_\zeta=\{e\}$, then either $G$ is twisted and
$\vert\zeta\vert=t$, that implies statement (e) of the lemma; or $G$ is
twisted, $\vert\zeta\vert=2t$, in particular, $t$ is odd, that implies point
(f) of the lemma.

Assume that $O^{2'}(G)\simeq A_2(2^{t})$, $\zeta$ is a graph-field automorphism
and $t$ is odd. If $\vert\zeta\vert\not=2t$, then arguments, using
Hartley-Shute Lemma  \ref{Hartley-Shute}, similar to the proof of
Lemma \ref{power} show that $C_U(H_\zeta)=\{e\}$, that contradicts to the
condition  $K\cap U\not=\{e\}$. If $\vert\zeta\vert=2t$, then we obtain point
(f) of the lemma.

Assume now that $O^{2'}(G)\simeq A_2(2^{2t})$ and $\zeta$ is a
graph-field automorphism. Again for $\vert\zeta\vert\not=2t$ from Hartley-Shute
Lemma \ref{Hartley-Shute} it follows that $C_U(H_\zeta)=\{e\}$, that
contradicts to the condition $K\cap U\not=\{e\}$. If $\vert\zeta\vert=2t$, then
either  $G_\zeta\simeq {}^2A_2(2^2)$, or $G_\zeta\simeq
\widehat{{}^2A_2(2^2)}$. If $G_\zeta\simeq {}^2A_2(2^2)$, then $H_\zeta=\{e\}$
and we obtain the statement (d) of the lemma. If $G_\zeta\simeq
\widehat{{}^2A_2(2^2)}$, then  $\vert H_\zeta\vert=3$, and so $KU/U\cap HU/U$
is a cyclic group  $\langle y\rangle$ of order $(2^{t_3}+1)_3=3^k$, where
$3^{k-1}=t_3$. If $k>1$, then Hartley-Shute Lemma \ref{Hartley-Shute} implies
that  $C_U(y)=\{e\}$, that is impossible. Thus  $t$ is nor divisible by $3$ and
$K\cap U$ is contained in the centralizer of an element $x$,
generating $H_\zeta$. Consider the homomorphism $GL_3(2^{2t})\rightarrow  \P
GL_3(2^{2t})$. Then some preimage of $x$ is similar to the matrix
$$
\left(
\begin{array}{ccc}
\lambda&0&0\\
0&\lambda^2&0\\
0&0&\lambda\\
\end{array}
\right),
$$
where $\lambda$ is the generating element of the multiplicative group of
$GF(2^2)$. The preimage of $U$ is similar with the set of upper
triangular matrices with the same elements on the diagonal. Direct
calculations show that $C_U(x)$ is isomorphic to the additive group of
$GF(2^{2t})$. The nilpotency of $K$ implies that $K\cap
U=(C_U(x))_{\zeta_{2'}}$, and point (c) of the lemma follows.

Assume now that $O^{2'}(G)\simeq {}^2A_2(2^{2t})$.  By Lemma
\ref{HomImageOfCarter}  $KU/U$ is a Carter subgroup of
$\langle B, \zeta g\rangle/U\simeq \langle H, \zeta g\rangle$ and, as above, we
may assume that  $H_\zeta\leq K$.  If
$\vert\zeta\vert=2t$, then $G_\zeta\simeq SL_2(2)$ and $H_\zeta=\{e\}$, and
point (e) of the lemma follows. Assume that  $t$ is even and
$\vert\zeta\vert\le t$. Then either $O^{2'}(G_\zeta)\simeq
SL_2(2^{2t/\vert\zeta\vert})$ (if the order $\vert\zeta\vert$ is even), or
$O^{2'}(G_\zeta)\simeq {}^2A_2(2^{2t/\vert\zeta\vert})$ (if the order
$\vert\zeta\vert$ is odd, hence  $\vert\zeta\vert<t$). Clearly
$H_\zeta$ contains an element $x$ such that ${K\cap U}\leq
C_U(H_\zeta)=\{e\}$, and this gives a contradiction with the condition ${K\cap
U}\not=\{e\}$. If $t$ is odd and $t\not=\vert\zeta\vert$, then
$O^{2'}(G_\zeta)\simeq
{}^2A_2(2^{2t/\vert\zeta\vert})$, and it follows that $H_\zeta$ contains
an element $x$ such that
$C_U(x)=\{e\}$. If $\vert\zeta\vert=t$ and  $t$ is odd, then the order
$\vert KU/U\cap B/U\vert$ can be divisible only by $3$ (otherwise by
Hartley-Shute Lemma \ref{Hartley-Shute} it again follows that
$C_U(H_\zeta)=\{e\}$). If $G_\zeta\simeq
{}^2A_2(2^{2t/\vert\zeta\vert})$, then $H_\zeta=\{e\}$ and we obtain point (b)
of the lemma. If $G_\zeta\simeq \widehat{{}^2A_2(2^{2t/\vert\zeta\vert})}$, then
$KU/U\cap HU/U$ is a cyclic group  $\langle y\rangle$ of order
$(2^{t_3}+1)_3=3^k$, where $3^{k-1}=t_3$. If $k>1$, then Hartley-Shute Lemma
\ref{Hartley-Shute} implies, that $C_U(y)=\{e\}$, that is impossible. Thus
$t$ is not divisible by $3$ and $K\cap U$ is contained in the centralizer of an
element  $x$, generating $H_\zeta$. As in the non-twisted case above, we obtain
that $C_U(x)$ is isomorphic to the additive group of
$GF(2^t)$. The nilpotency of $K$ implies that $K\cap
U=(C_U(x))_{\zeta_{2'}}$, and point (a) of the lemma follows.
\end{proof}

\section{Carter subgroups of semilinear groups}

\subsection{Brief review of results of this section}

In this section, by using notations and results obtained in section 4, we
classify Carter subgroups in groups of automorphisms of finite groups of Lie
type. First we give such a classification in the case, when a group of Lie
type has type $C_n$ or when a group of its automorphisms contains a triality
automorphism, since the arguments in these two cases differ from the
remainings. The we formulate the final theorem and we prove this theorem in two
subsections. In the last subsection we prove that in every finite group with
known composition factors Carter subgroups are conjugate.

\subsection{Carter subgroups of symplectic groups}

Consider a set $\mathcal{A}$ of almost simple groups $A$ such that a unique
non-Abelian composition factor  $S=F^\ast(A)$ is a canonical simple group of
Lie type and $A$ contains nonconjugate Carter subgroups. If the set
$\mathcal{A}$ is not empty, denote by $\textrm{\bfseries{Cmin}}$ the minimal
possible order of  $F^\ast(A)$, with $A\in\mathcal{A}$.  If the set
$\mathcal{A}$ is empty, then let $\textrm{\bfseries{Cmin}}=\infty$. We shall
prove that $\textrm{\bfseries{Cmin}}=\infty$, i.~e. that
$\mathcal{A}=\varnothing$. Note that if $A\in\mathcal{A}$ and
$G=F^\ast(A)$, then there exists a subgroup $A_1$ of $A$ such that
$A_1\in\mathcal{A}$ and $A_1=KG$ for a Carter subgroup $K$ of $A$. Indeed, if
for every nilpotent subgroup $N$ of $A$ Carter subgroups of $NG$ are conjugate,
then  $A$ satisfies {\bfseries (C)}, hence Carter subgroups of $A$ are
conjugate, that contradicts to the choice of $A$. So there exists a nilpotent
subgroup  $N$ of $A$ such that Carter subgroups of $NG$ are not conjugate. Let
$K$ be a Carter subgroup of $NG$. Then clearly $KG/G$ is a Carter subgroup of
$NG/G$, i.~e., coincides with $NG/G$. Therefore Carter subgroups of
$KG$ are not conjugate and~${KG=A_1\in\mathcal{A}}$. So the condition
$A=KG$ in Theorems \ref{sympcarter}, \ref{CarterTriality}, and
\ref{CarterSemilinear} is not a restriction and is used only to simplify
arguments.

In this section we consider Carter subgroups in an almost simple
group $A$ with simple socle $G=F^\ast(A)\simeq \P \Sp_{2n}(q)$. We
consider such groups in the separate section, since for groups of type $\P
\Sp_{2n}(q)$
Lemma \ref{centUH} is not true and we use arguments slightly
different from those that we use in the  proof of
Theorem~\ref{CarterSemilinear}.

We shall prove first two technical lemmas.

\begin{lem} \label{centUHsymp}
Let $O^{p'}(\overline{G}_\sigma)= G$ be a canonical adjoint finite group of
Lie type over a field of odd characteristic $p$ and $-1$ is not a square in the
base field of $G$. Assume that the root system $\Phi$ of $\overline{G}$
equals $C_n$.  Let $U$ be a maximal unipotent subgroup of $G$, 
$H$ be a Cartan subgroup of $G$, normalizing $U$, and $Q$ is a Sylow
$2$-sub\-gro\-up of $H$. 

\noindent Then 
$C_U(Q)=\langle X_r\mid r\textrm{ is a long root}\rangle$.
\end{lem}

\begin{proof}
If $r$ is a short root, then there exists a root $s$ with $<s,r>=1$. Thus 
$$x_{r}(t)^{h_{s}(-1)}=x_{r}((-1)^{<s,r>}t)=x_r(-t)$$ (see
\cite[Proposition~6.4.1]{CarSimpleGrpsLieType}). Therefore, if
$x\in C_U(Q)$ and $x_r(t)$ is a nontrivial multiplier in decomposition
\eqref{canonicalform} of $x$, then  $r$ is a long root.  Now if $r$ is a long
root, then for every root $s$ either $\vert <s,r>\vert=2$, or $<s,r>=0$, i.~e.,
$x_r(t)^{h_s(-1)}=x_r(t)$. Under the condition that $-1$ is not a square in the
base field of $G$ (i.~e., in the field $\F_{q}$) we obtain that 
$q\equiv-1\pmod4$, so $\langle h_s(-1)\mid s\in\Phi\rangle=Q$, and the lemma
follows.
\end{proof}

\begin{lem}\label{InducedAutomorphismsInCn}
Let $G=\P\Sp_{2n}(q)$ be a simple canonical group of Lie type, 
$J$ be a subset of the set of fundamental roots, containing the long
fundamental root $r_n$, $P_J$ be a parabolic subgroup, generated by a Borel
subgroup $B$ and by groups $X_r$ with $-r\in J$, $L$ be a Levi factor of $P_J$.
Denote by $S$ a quasisimple normal subgroup of $L$, isomorphic to
$\Sp_{2k}(q)$ (it always exists, since $r_n\in J$).

\noindent Then~$Aut_L(S/Z(S))=S/Z(S)$.
\end{lem}

\begin{proof}
This statement is known, it is proven in an unpublished paper by N.A.Vavilov.
We give a proof here for the completeness. As we noted above, $L$
is a reductive subgroup of maximal rank of $G$, ans so the following inclusions
hold $S/Z(S)\leq \Aut_L(S/Z(S))\leq \widehat{S/Z(S)}$. Since $\vert
\widehat{C_n(q)}:C_n(q)\vert=(2,q-1)$, then for $q$ even the statement is
evident. If $q$ is odd, then for $\Aut_L(S/Z(S))$ there can be only two
possibilities: either $\Aut_L(S/Z(S))=S/Z(S)$, or
$\Aut_L(S/Z(S))=\widehat{S/Z(S)}$. We shall show that the second equality is
impossible.

In our notations fundamental roots of the root system of $S$ are
$r_{n-k+1},\ldots,r_n$. If the equality $\Aut_L(S/Z(S))=\widehat{S/Z(S)}$
holds, then there exist elements $s_1,\ldots,s_k$ of $\Z\Phi=\Z C_n$ such that
$$<s_i,r_{n-k+j}>=\frac{(s_i,r_{n-k+j})}{(s_i,s_i)}=\left\{\begin{array}{lr}
1&\text{if }i=j,\\ 0&\text{if }i\not=j.
\end{array} \right.$$ (They generate the lattice of fundamental weights, thus
allow to obtain all diagonal automorphisms of $S$). But for each root $s$ of
$C_n$ we have that either $<s,r_n>=0$, or $<s,r_n>=\pm2$, i.~e., for each
element $s\in\Z\Phi$ the number $<s,r_n>$ is even, in particular is distinct
from $1$. Therefore  such a set of
elements $s_1,\ldots,s_k$ does not exists.
\end{proof}

\begin{ttt}\label{sympcarter}
Let $G$ be a finite adjoint group of Lie type (not necessary simple) over a
field of characteristic $p$,  and
$\overline{G}$, $\sigma$ are chosen so that  $\P \Sp_{2n}(p^t)\simeq
O^{p'}(\overline{G}_\sigma)\leq G\leq\overline{G}_\sigma$. Choose a
subgroup $A$ of $A\cap \ov{G}_\sigma=G$. Let
$K$ be a Carter subgroup of $A$. Assume also that~{\em$\vert \P
\Sp_{2n}(p^t)\vert\le \textbf{Cmin}$}
and~${A=KG}$.

\noindent Then exactly one of the following statements holds:

\begin{itemize}
\item[{\em (1)}] $G$ is defined over   $GF(2^t)$, a field
automorphism $\zeta$ is in $A$, $\vert\zeta\vert=t$, and, up to conjugation in
$G$, the equality $K=Q\leftthreetimes\langle\zeta\rangle$ holds, where $Q$ is a
Sylow $2$-sub\-gro\-up of~$G_{\zeta_{2'}}$.
\item[{\em (2)}] $G\simeq \P \SL_2(3^{t})\simeq \P \Sp_2(3^t)$, a field
automorphism $\zeta$ is in $A$, $\vert\zeta\vert=t$ is odd,  and, up to
conjugation in $G$, the equality $K=Q\leftthreetimes\langle\zeta\rangle$ holds,
where $Q$ is a  Sylow $3$-sub\-gro\-up of~$G_{\zeta_{3'}}$.
\item[{\em (3)}] $p$ does not divide $\vert K\cap G\vert$ and $K$ is contained
in the normalizer of a Sylow $2$-sub\-gro\-up of~$A$.
\end{itemize}
In particular, Carter subgroups of  $A$ are conjugate, i.~e., if
$A_1\in\mathcal{A}$ and {\em$F^\ast(A_1)=\textbf{Cmin}$}, then
$F^\ast(A_1)\not\simeq \P \Sp_{2n}(p^t)$.
\end{ttt}

\begin{proof}
Assume that the theorem is not true and  $A$ is a counter example
such that $\vert F^\ast(A)\vert$ is minimal. Note that no more than one
statement of the theorem can be fulfill, since if statement
(2) holds, then, by Lemmas  \ref{Norm2Syl} and \ref{ESyl2InhFieldGraph}, for a
Sylow $2$-sub\-gro\-up
$Q$ of $A$ the condition $N_G(S)=S C_G(S)$ is not true, i.~e.,
statement (3) of the theorem does not hold. Thus if $A_1$ is an
almost simple group with $F^\ast(A_1)$ being a simple group of Lie type of
order less, than $\vert F^\ast(A)\vert$, then Carter subgroups of $A_1$ are
conjugate. In view of Theorem \ref{CarterInClassicalGroups} we may assume that
$A\not=G$. Moreover, by Theorem \ref{main},  we may assume
that $q$ is odd, i.~e., that $\Aut(\P \Sp_{2n}(q))$ does not contain
a graph automorphism. Thus we may assume that $A=\langle G,\zeta
g\rangle$. 

Assume that $K$ is a Carter subgroup of $\langle G,\zeta g\rangle$
and $K$ does not satisfy to the statement of the theorem. Write $K=\langle
x,K\cap G\rangle$. If either $p\not= 3$ or
$t$ is even, then the theorem follows from Theorem \ref{main}. Thus we may
assume
that $q=3^t$ and $t$
is odd. Since $\vert \ov{G}_\sigma:O^{p'}(\ov{G_\sigma})\vert=2$ and the order
$\vert\zeta\vert$ is odd, we may assume that the order $\vert\zeta g\vert$ is
also odd and so $\zeta\in\langle G,\zeta g\rangle$, i.~e.,
$A=G\leftthreetimes\langle\zeta\rangle$. By Lemma \ref{omnibus}
every semisimple
element of odd order is conjugate to its inverse in $G$. Now, for
every semisimple element $t\in G$, each non-Abelian composition
factor of $C_G(t)$ is a simple group of Lie type (see
\cite{Car5CentSemisimpleLie}) of
order less, than $\textrm{\bfseries{Cmin}}$. Therefore, for every non-Abelian
composition factor $S$ of $C_A(t)$ and every nilpotent subgroup $N\leq C_A(t)$,
Carter subgroups of
$\langle \Aut_N(S),S\rangle$ are conjugate. It follows that $C_A(t)$ satisfies
{\bfseries (C)}. Hence, by Lemma \ref{power}, $\vert
K\cap G\vert=2^\alpha\cdot 3^\beta$ for some $\alpha,\beta\ge 0$.

If $G=\widehat{\P \Sp_{2n}(q)}$ then by \cite[Theorem~2]{WonConjByInvolution}
every unipotent element is conjugate to its inverse. Since $3$
is a good prime for $G$, then \cite[Theorems~1.2 and~1.4]{SeitzConjUnipElts}
imply that, for any element $u\in G$ of order $3$, all composition
factors of $C_G(u)$ are simple groups of Lie type of order less,
than $\textrm{\bfseries{Cmin}}$. Thus $C_A(u)$ satisfies {\bfseries (C)}, hence,
by Lemma \ref{power}, we obtain that $K\cap G$ is a $2$-gro\-up. By Lemmas
\ref{Syl2InCentrOfFieldAut} and \ref{ConjAutomorphisms} every
element  $x\in A\setminus G$ of odd order with $\langle x\rangle\cap
G=\{e\}$ centralizes some Sylow $2$-sub\-gro\-up of $G$. Hence $K$
contains a Sylow $2$-sub\-gro\-up of $G$, and hence of $A$, i.~e., $K$ satisfies
statement (3) of the theorem.

Thus we may assume that $G=\P Sp_{2n}(q)$ and $\beta\ge
1$, i.~e., a Sylow $3$-sub\-gro\-up $O_3(K\cap G)$ of $K\cap G$ is nontrivial.
By Lemma \ref{parabolic} we obtain that $K\cap G$ is contained in some
$K$-in\-va\-ri\-ant parabolic subgroup $P$ of $G$ with a Levi factor $L$
and, up to conjugation in $P$, a Sylow $2$-sub\-gro\-up $O_2(K\cap G)$ of
$K\cap G$ is contained in $L$. Note that all non-Abelian composition factors of
$P$ are simple groups of Lie type of order less, than {\bfseries Cmin}, so $P$
and each its homomorphic image satisfy {\bfseries(C)}. The group
$\widetilde{K}=KO_3(P)/O_3(P)$ is
isomorphic to $K/O_3(K\cap G)$ and, by Lemma
\ref{HomImageOfCarter}, $\widetilde{K}$ is a Carter subgroup of
$\langle \widetilde{K},P/O_3(P)\rangle$. Now $\widetilde{K}\cap
P/O_3(P)\simeq O_2(K\cap G)$
is a $2$-gro\-up and every element $x\in \langle
\widetilde{K},P/O_3(P)\rangle\setminus P/O_3(P)$ of odd order with $\langle
x\rangle\cap P/O_3(P)=\{e\}$  centralizes a Sylow $2$-sub\-gro\-up of
$P/O_3(P)\simeq L$ (see Lemmas \ref{Syl2InCentrOfFieldAut} and
\ref{ConjAutomorphisms}).
Therefore $O_2(K\cap G)$ contains a Sylow $2$-sub\-gro\-up of $L$, in
particular, contains a Sylow $2$-sub\-gro\-up $H_2$ of $H$. Since $K$ is
nilpotent,
Lemma \ref{centUHsymp}
implies that $O_3(K\cap G)\leq C_U(H_2)=\langle X_r\mid r
\textrm{ is a long root of }\Phi(G)^+\rangle$. Since for every two long positive
roots $r,s$ in $\Phi(G)^+$ we have that $r+s\not\in\Phi(G)$, Chevalley
commutator
formula \cite[Theorem~5.2.2]{CarSimpleGrpsLieType} (Lemma
\ref{ChevalleyFormulae}) implies that $\langle X_r\mid r
\textrm{ is a long root of }\Phi(G)^+\rangle$ is Abelian.

Since $\zeta$ is a field automorphism, it normalizes each parabolic subgroup of
$G$ containing a $\zeta$-stab\-le Borel subgroup. Thus for every subset $J$ of
the set of fundamental roots $\Pi=\{r_1,\ldots,r_n\}$ of $\Phi=\Phi(G)$ the
parabolic subgroup $P_J$ is $\zeta$-stab\-le. Therefore we may suppose that
$P=P_J$, where $J$ is a proper subset of the set of fundamental roots $\Pi$ of
$\Phi$. Choose the numbering of fundamental roots so that $r_n$ is a long
fundamental root, while the remaining fundamental roots $r_i$ are short roots.
If $r_n\in J$, then one of the components of the Levi factor, $G_1$ for
example, is isomorphic to $Sp_{2k}(q)$ for some $k<n$ (note that since
$A\not=G$ then $q\not=3$). By Lemma \ref{InducedAutomorphismsInCn} we obtain
that
$L/C_L(G_1)=\Aut_L(G_1/Z(G_1))=G_1/Z(G_1)$. By Lemma \ref{HomImageOfCarter}
$K_1=KC_L(G_1)O_3(P)/C_L(G_1)O_3(P)$ is a Carter subgroup of
$(P\leftthreetimes\langle\zeta\rangle)/C_L(G_1)O_3(P)$. Since $\vert K_1\cap
P/C_L(G_1)O_3(P)\vert$ is not divisible by $3$, and $\zeta$ centralizes a Sylow
$2$-sub\-gro\-up of $G_1/Z(G_1)$ (see Lemma \ref{Syl2InCentrOfFieldAut}), then
$K_1$ contains a Sylow $2$-sub\-gro\-up of $P/C_L(G_1)O_3(P)\simeq
G_1/Z(G_1)\simeq
\P Sp_{2k}(q)$. Moreover by Lemma \ref{Syl2InCentrOfFieldAut} a Sylow
$2$-sub\-gro\-up of $(P/C_L(G_1)O_3(P))_\zeta$ is a Sylow $2$-sub\-gro\-up of
$P/C_L(G_1)O_3(P)$. Thus $K_1\cap P/C_L(G_1)O_3(P)$ is a Sylow $2$-sub\-gro\-up
of
$(P/C_L(G_1)O_3(P))_\zeta\simeq \P Sp_{2k}(3)$. By Lemma \ref{Norm2Syl} there
exists an element $x$ of odd order of $\P Sp_{2k}(3)$ that normalizes but not
centralizes a Sylow $2$-sub\-gro\-up; a contradiction with the fact that $K_1$
is a
Carter subgroup of $(P\leftthreetimes\langle\zeta\rangle)/C_L(G_1)O_3(P)$. Thus
we may assume that~${r_n\not\in J}$.

Consider the set $J_n=\Pi\setminus\{r_n\}$ and the parabolic subgroup
$P_{J_n}$. From the above arguments it follows that $K\leq
P_J\leftthreetimes\langle \zeta\rangle\leq
P_{J_n}\leftthreetimes\langle\zeta\rangle$. Now the subgroup $\langle X_r\mid r
\textrm{ is a long root of } \Phi(G)^+\rangle$ is contained in $O_3(P_{J_n})$
and $O_3(K\cap G)$ is contained in $\langle X_r\mid r
\textrm{ is a long root of } \Phi(G)^+\rangle$, so $N_G(O_3(K\cap G))\leq
O_3(P_{J_n})$ and we may assume that $P=P_{J_n}$. By Lemma
\ref{HomImageOfCarter}, $\widetilde{K}=KO_3(P)/O_3(P)$ is a Carter subgroup of
$(P\leftthreetimes\langle\zeta\rangle)/O_3(P)$. Note that a unique
non-Abelian composition factor of $P\leftthreetimes\langle\zeta\rangle$ is
isomorphic to $A_{n-1}(q)\simeq \P SL_n(q)$. By
\cite[Theorem~1]{KondNormalizers} and \cite[Theorem~4]{KonMaz} we obtain that
$\widetilde{K}=R\times\langle\zeta\rangle$, where $R$ is a Sylow
$2$-sub\-gro\-up
of $P$ centralized by $\zeta$. Thus $O_3(K\cap G)\leq C_P(R)$. Consider
$Q=O_3(K\cap G)\cap P_\zeta$. Since $O_3(K\cap G)$ is nontrivial and $K$ is
nilpotent, then $Q=O_3(K\cap G)\cap P_\zeta=Z(K)\cap O_3(K\cap G)$ is
nontrivial. Therefore $N_G(Q)$ is a proper subgroup of $G$ and by Lemma
\ref{parabolic} $N_G(Q)$ is contained in a proper parabolic subgroup of $G$. On
the other hand, $K\leq N_G(Q)$ and $P=P_{J_n}$ is a maximal proper parabolic
subgroup of $G$. If $N_G(Q)$ is not contained in $P$, then $N_G(Q)$ and $K$ are
contained in a parabolic subgroup $P_J$ with $r_n\in J$. We have proved above
that $r_n\not\in J$, so $N_G(Q)$ is contained in~$P$.

We shall show that $R\times Q$ is a Carter subgroup of $G_\zeta$. Indeed,
assume that an element $x\in G_\zeta$ normalizes $R\times Q$. Then $x$
normalizes $Q$, so $x$ is in $P$ and normalizes $O_3(P)$. On the other hand $x$
normalizes $R$, therefore normalizes $C_P(R)$, so $x$ normalizes
$C_{O_3(P)}(R)$.  Moreover it is evident that $x$ and $\zeta$ commute. Thus $x$
normalizes $(R\times C_{O_3(P)}(P))\leftthreetimes \langle\zeta\rangle$. As we
noted above, $K\leq (R\times C_{O_3(P)}(P))\leftthreetimes \langle\zeta\rangle$
and $(R\times C_{O_3(P)}(P))\leftthreetimes \langle\zeta\rangle$ is solvable.
Lemma \ref{power}(a) implies that $(R\times C_{O_3(P)}(P))\leftthreetimes
\langle\zeta\rangle$ coincides with its normalizer in
$G\leftthreetimes\langle\zeta\rangle$, so $x\in R\times C_{O_3(P)}(R)$. The
group $C_{O_3(P)}(R)\leq \langle X_r\mid r \textrm{ is a long root of }
\Phi(G)^+\rangle$ is Abelian, so every element of $R\times C_{O_3(P)}(R)$
centralizes $C_{O_3(P)}(R)\geq O_3(K\cap G)$. Therefore $x$ normalizes $(R\times
C_{O_3(P)}(P))\leftthreetimes \langle\zeta\rangle=K$, i.~e., $x\in K$. By
construction $R\times Q=K\cap G_\zeta$, so $x\in R\times Q$ and $R\times Q$ is
a Carter subgroup of $G_\zeta$. On the other hand $O^{3'}(G_\zeta)\simeq \P
Sp_{2n}(3^{t/\vert\zeta\vert})$ and by induction groups  $\P
Sp_{2n}(3^{t/\vert\zeta\vert})$ and $\widehat{\P
Sp_{2n}(3^{t/\vert\zeta\vert})}$ does not contain Carter subgroups of order
divisible by $3$.
This final contradiction completes the proof of the theorem.
\end{proof}

\subsection{Группы с автоморфизмом тройственности}

\begin{ttt}\label{CarterTriality}
Let  $G$ be a finite adjoint  group of Lie type over a field of characteristic
$p$, $\ov{G}$, $\sigma$ are chosen so that $O^{p'}(\ov{G}_\sigma)\leq
G\leq \ov{G}_\sigma$, and $O^{p'}(\ov{G}_\sigma)$ is isomorphic to  either
$D_4(q)$, or $ {}^3D_4(q^3)$. Assume that $\tau$ is a graph automorphism of
order $3$ of $O^{p'}(G)$ (recall that for  $G\simeq {}^3D_4(q^3)$ $\tau$ is an
automorphism such that the set of its stable points is isomorphic to
$G_2(q)$). Denote by $A_1$ the
subgroup of $\Aut(D_4(q))$ generated by inner-diagonal and field automorphisms,
and also by a graph automorphism of order $2$. Let
$A\leq \Aut(G)$ be such that $A\not\leq A_1$ (if $O^{p'}(G)\simeq D_4(q)$),
and $K$ be a Carter subgroup of $A$.  Assume also that $\vert O^{p'}(G)\vert\le
\text{\em\bfseries{Cmin}}$, $G=A\cap \ov{G}_\sigma$ and $A=KG$. Then one of the
following statements holds:

\begin{itemize}
\item[{\em(a)}] $G\simeq {}^3D_4(q^3)$, $(\vert A:G\vert,3)=1$, $q$ is odd and
$K$ contains a Sylow $2$-sub\-gro\-up of~$A$;
\item[{\em(b)}] $(\vert A:G\vert,3)=3$, $q$ is odd, $\tau\in A$ and, up to
conjugation by an element of $G$, the subgroup $K$ contains a Sylow
$2$-sub\-gro\-up of $C_A(\tau)\in \Gamma G_2(q),$
and $\tau\in K$;
\item[{\em(c)}] $(\vert A:G\vert,3)=3$, $q=2^t$, ${\vert A:G\vert=3t}$,
$A=G\leftthreetimes\langle \tau,\varphi\rangle$, where  $\varphi$ is a field
automorphism of order $t$ commuting with $\tau$ and, up to conjugation by an
element of $G$,  the subgroup $K$ contains a Sylow $2$-sub\-gro\-up of
$C_G(\langle \tau,\varphi\rangle_{2'})\simeq G_2(2^{t_{2'}})$ and
$\tau\in K$;
\item[{\em(d)}] $O^{p'}(G)\simeq D_4(p^{3t})$, $p$ is odd, the factor group
$A/G$ is cyclic,  $\tau\not\in A$,
$A=G\leftthreetimes\langle \zeta\rangle$, where for some natural $m$,
$\zeta=\tau\varphi^m$ is a graph-field automorphism, and ,
up to conjugation by an element of $G$,
$K=Q\leftthreetimes\langle\zeta\rangle$, where $Q$ is a Sylow
$2$-sub\-gro\-up of $C_G(\zeta_{2'})\simeq
{}^3D_4(p^{3t/\vert\zeta_{2'}\vert})$.
\end{itemize}
In particular,  Carter subgroups of  $A$ are conjugate, i.~e., if  $A_2\in
\mathcal{A}$ and $\vert F^\ast(A_2)\vert=\text{\em\bfseries{Cmin}}$, then
$A_2$ does not satisfy to the conditions of the theorem, so
$F^\ast(A_2)\not\simeq{}^3D_4(q^3)$.
\end{ttt}

\begin{proof}
Assume that the theorem is not true and  $A$ is a counter example
such that  $\vert O^{p'}(G)\vert$ is minimal. In view of
\cite[Theorem~1.2(vi)]{TiepZalRealConjClasses} we have that every element of $G$
is conjugate to its inverse. By  \cite{Car5CentSemisimpleLie} and 
\cite[Theorems~1.2 and~1.4]{SeitzConjUnipElts} we obtain that for every element
$t\in G$ of odd prime
order, all non-Abelian composition factors of $C_G(t)$ are simple
groups of Lie type of order less, than $\textrm{\bfseries{Cmin}}$. Thus,
$C_A(t)$ satisfies {\bfseries (C)} and Lemma \ref{power} implies that
$K_G=K\cap G$ is a $2$-gro\-up. Now Lemma \ref{ConjAutomorphisms}
implies that all cyclic groups, generated by field automorphisms of the same
odd order of $G$, are conjugate under $G$. Since the centralizer of every field
automorphism in $G$ is a group of Lie type of order less, than {\bfseries
Cmin}, we again use Lemma \ref{power} and obtain the statement of the theorem
by induction. Lemma \ref{ConjAutomorphisms} implies also that if
$O^{p'}(G)\simeq D_4(q)$, then all cyclic groups generated by graph-field
automorphisms are conjugate. Since the centralizes of each graph-field
automorphism in  $G$ is a group of Lie type of order less than
$\textrm{\bfseries{Cmin}}$, we again use Lemma \ref{power} and obtain statement
(d) of the theorem by induction. Thus we may assume that $A$ does
not contain a field automorphism or a graph-field automorphism
of odd order. Therefore either $G\simeq{}^3D_4(q^3)$ and $A/G$ is a $2$-gro\-up,
or $K$ contains an element $s$ of order $3$ such that $\langle s\rangle\cap
A_1=\{e\}$ (for groups ${}^3D_4(q^3)$ the equality $\langle s\rangle \cap
G=\{e\}$ holds), $G\leftthreetimes \langle s\rangle=G\leftthreetimes \langle
\tau\rangle$ and $K\cap  G$  is a
$2$-gro\-up. 

In the first case we obtain the statement (a) of the theorem with condition
$(\vert
A:G\vert,3)=1$. In the second case there exists two non-conjugate
cyclic subgroups $\langle \tau\rangle$ and $\langle x\rangle$ of order $3$
of $A$ such
that $\langle\tau\rangle\cap A_1=\langle
x\rangle\cap
A_1=\{e\}$ and $G\leftthreetimes\langle x\rangle=G\leftthreetimes\langle
\tau\rangle$ (see
\cite[(9-1)]{GorLyoLocalStructure}). Hence, either $s=\tau\in K$, or $s=x\in K$.
Assume  that $q\not=3^t$. In the first case from the known structure of Carter
subgroups in a group from the set $\Gamma G_2(q)$, obtained in Theorem
\ref{main}, the statement (b) or (c)  of the
theorem follows,
in the second case we have that $K\leq C_A(x)$. By
\cite[(9-1)]{GorLyoLocalStructure} 
$C_G(x)\simeq\P
GL_3^\varepsilon(q)$, where $q\equiv \varepsilon1\pmod 3$, $\varepsilon=\pm$ and
$\P GL_3^{+}(q)=\P GL_3(q)$, $\P GL_3^{-}(q)=\P GU_3(q)$. Then $K=(K\cap
G)\times \langle y,\varphi\rangle$, where $\varphi$ is a field automorphism of
$O^{p'}(G)$ of order equal to a power of $2$ and $y$ is a graph automorphism
such that its order is a power of $3$ and $x\in \langle y\rangle$. By nilpotency
of
$K$ we obtain that $y\varphi=\varphi y$, it follows that
$C_{C_G(\varphi)}(x)=C_{C_G(x)}(\varphi)$. Now we have that
$$C_G(\varphi)=\left\{\begin{array}{ll}
D_4(q^{1/\vert\varphi\vert}),&\text{if }O^{p'}(G)\simeq D_4(q),\\
{}^3D_4(q^{3/\vert\varphi\vert}),&\text{if }G\simeq {}^3D_4(q^3).
                      \end{array}\right.$$
Hence $C_{C_G(x)}(\varphi)=C_{C_G(\varphi)}(x)\simeq \P
GL_3^\mu(q^{1/\vert\varphi\vert})$, with
$q^{1/\vert\varphi\vert}\equiv\mu1\pmod3$, where $\mu=\pm$ (note that
$\varepsilon$ and $\mu$ can be different). As we noted above, $K\cap G$ is a
$2$-gro\-up. On the other hand, by \cite[Theorem~4]{KonMaz} there exists an
element $z$ of order $3$ centralizing a Sylow $2$-sub\-gro\-up of $C_G(x)=\P
GL_3^\varepsilon(q)$ and belonging to $C_{C_G(x)}(\varphi)\simeq \P GL_3^\mu
(q^{1/\vert\varphi\vert})$. Thus $z$ centralizes $K$, hence is in $K$. But
$K\cap G$ does not contain elements of odd order, therefore this second case is
impossible.

Assume now that $q=3^t$. Then $C_G(\tau)\simeq G_2(q)$ and we obtain the
theorem. In
the second case $C_G(x)\simeq SL_2(q)\rightthreetimes U$, where $U$ is a
$3$-gro\-up
and $Z(C_G(x))\cap U\not=\{e\}$, a contradiction with
Lemma~\ref{power}.
\end{proof}

\subsection{Classification theorem}
\begin{ttt}\label{CarterSemilinear}
Let $G$ be a  finite adjoint group of Lie type  ($G$ is not necessary
simple) over a field of characteristic $p$  and
$\overline{G}$,
$\sigma$ are chosen so that $O^{p'}(\overline{G}_\sigma)\leq
G\leq\overline{G}_\sigma$. Assume also that $G\not\simeq
{}^3D_4(q^3)$. Choose a subgroup $A$ of $\Aut
(O^{p'}(\ov{G}_\sigma))$ with  $A\cap \ov{G}_\sigma=G$ and, if
$O^{p'}(G)=D_4(q)$, assume that $A$ is
contained in the subgroup $A_1$ defined in Theorem {\em\ref{CarterTriality}}.
Let $K$ be a Carter subgroup of $A$ and assume that~$A=KG$.

\noindent Then exactly one of the following statements holds:

\begin{itemize}
\item[{\em (a)}] $G$ is defined over a field of characteristic $2$,
$A=\langle\langle G,\zeta g,t\rangle$, where $t$ is a $2$-ele\-ment, $K$ is
contained in the normalizer of a $t$-stab\-le Borel subgroup of
$G$ $K\cap \langle
G,\zeta g\rangle$ satisfies to one of the statements {\em(a)--(f)} of
Lemma~{\em \ref{CartBorel}};
\item[{\em (b)}] $G\simeq \P SL_2(3^{t})$, a field automorphism $\zeta$ is in
$A$,  $\vert\zeta\vert=t$ is odd, and, up to conjugation in $G$, the equality
$K=Q\leftthreetimes\langle\zeta\rangle$ holds, where $Q$ is a   Sylow
$3$-sub\-gro\-up of~$G_{\zeta_{3'}}$;
\item[{\em (c)}] $A={}^2G_2(3^{2n+1})\leftthreetimes\langle\zeta\rangle$,
$\vert\zeta\vert=2n+1$, and, up to conjugation in $G$ the equality $K=(K\cap G)
\leftthreetimes\langle\zeta\rangle$ holds, and $K\cap
{}^2G_2(3^{2n+1})=Q\times P$, where
$Q$ is of order $2$ and~${\vert P\vert =3^{\vert\zeta\vert_3}}$.
\item[{\em (d)}] $p$ does not divide $\vert K\cap G\vert$ and $K$  contains  a
Sylow $2$-sub\-gro\-up of~$A$, in view of Lemma \ref{ESyl2InhFieldGraph} $A$
satisfies {\em\bfseries(ESyl2)} if and only if $G$
satisfies~{\em\bfseries(ESyl2)}.
\end{itemize}
In particular, Carter subgroups of $A$ are conjugate.
\end{ttt}

\begin{rem}
There exists a dichotomy for Carter subgroups in groups of automorphisms of
finite groups of Lie type, not containing a graph, or a graph-field
automorphism of order $3$. They either are contained in the normalizer of a
Borel subgroup, or the characteristic is odd and a Carter subgroup contains a
Sylow $2$-sub\-gro\-up of the hole group.
\end{rem}

Assume that the theorem is not true and $A$ is a counter example to
the theorem with $\vert F^\ast(A)\vert$ minimal. Among counterexamples with
$\vert F^\ast(A)\vert$ minimal take those, for which $\vert
A\vert$ is minimal. In this case for every
almost simple group $A_1$ such that $\vert F^\ast(A_1)\vert<\vert
F^\ast(A)\vert$, $F^\ast(A_1)$ is a finite simple group of Lie
type  and $A_1$ satisfies the conditions of Theorem
\ref{CarterSemilinear}, Carter subgroups are conjugate. Indeed, note that no
more, than one
statement of the theorem can be fulfill, since if either statement
(b), or statement (c) of the theorem holds, then the condition
$N_A(Q)=Q C_A(Q)$ for a Sylow $2$-sub\-gro\-up $Q$ of $A$ is not true,
i.~e., the statement (d) of the theorem does not hold (the fact that
other statements can not hold simultaneously is evident). Thus,
Carter subgroups of  $A_1$ are conjugate. Note also that from this
fact we immediately obtain the inequality $\vert F^\ast(A)\vert\le
\textrm{\bfseries{Cmin}}$. Indeed, if $A_2\in\mathcal{A}$ and
$F^\ast(A_2)=\textrm{\bfseries{Cmin}}$, then
either $A_2$ satisfies to the condition of Theorem
\ref{CarterTriality}, or $A_2$ satisfies conditions of Theorem
\ref{CarterSemilinear}. As we noted in Theorem \ref{CarterTriality},
the first case is impossible. The second case, as we just noted, is possible
only if $\vert F^\ast(A)\vert\le \vert
F^\ast(A_2)\vert=\textrm{\bfseries{Cmin}}$
(since $A$ is a counterexample to the statement of the theorem with $\vert
F^\ast(A)\vert$ is minimal).

We shall prove the theorem in the following way. If $F^\ast(A)\simeq
\P Sp_{2n}(q)$, then the theorem follows from Theorem \ref{sympcarter}. If
$A=G$, then the theorem follows from \cite{DT2CartSbgrpsPGL},
\cite{DTZCartSbgrpsClassGrps} and results from sections 3 and 4 of the present
paper. Thus we may assume, that
$A/(A\cap G)$ is nontrivial. Let $K$ be a Carter subgroup of $A$. We
shall prove first that if $p$ divides  $\vert K\cap G\vert$, then
one of the statements (a)--(c) of the theorem holds. Then we shall
prove that if $p$ does not divide $\vert K\cap G\vert$, then $K$
contains a Sylow $2$-sub\-gro\-up of $A$. Since both of these steps are quite
complicated, we divide them onto two sections. Note that, in view of
\cite{Car5CentSemisimpleLie}, for every semisimple element  $t\in G$,
all non-Abelian
composition factors of  $C_G(t)$, so of $C_A(t)$, are simple groups
of Lie type of order less, than  $\vert F^\ast(A)\vert$, and hence of order
less, than $\textrm{\bfseries{Cmin}}$. Therefore $C_A(t)$ satisfies {\bfseries
(C)}. In
order to apply Lemmas \ref{HomImageOfCarter} and \ref{power} we
shall use this fact without future references.

\subsection[Carter subgroups with unipotent radical]{Carter subgroups of order
divisible by characteristic}

Denote $K\cap G$ by $K_G$. For every group $A$, satisfying conditions of
Theorem \ref{CarterSemilinear},  the factor group $A/G$ is Abelian and, for
some natural $t$ is isomorphic to a subgroup of $\Z_2\times \Z_t$, where $\Z_t$
denotes a cyclic group of order $t$. If the factor group $A/G$ is not cyclic,
then the group $O^{p'}(G)$ is split and $A$ contains an element $\tau a$, where
$\tau$ is a graph automorphism of $O^{p'}(G)$ and $a\in\ov{G}_\sigma$. Then
every semisimple element of odd order is conjugate to its inverse in $A$ (cf.
Lemma \ref{ConjInverseInGraph}). By Lemma \ref{power} we obtain that
$\vert K_G\vert$ is divisible only by $2$ and $p$. If $p=2$, then we obtain that
$K_G$ is a $2$-gro\-up, it is contained in a proper $K$-in\-va\-ri\-ant
parabolic
subgroup $P$ of $G$  and by Lemma \ref{HomImageOfCarter} $KO_2(P)/O_2(P)$ is a
Carter subgroup of $KP/O_2(P)$. Since $K_G\leq O_2(P)$, then
$\left(KO_2(P)/O_2(P)\right)\bigcap \left(P/O_2(P)\right)=\{e\}$. Hence $P$ is
 a Borel subgroup of $G$, otherwise we would have
$C_{P/O_2(P)}(KO_2(P)/O_2(P))\not=\{e\}$, a contradiction with the fact that $
KO_2(P)/O_2(P)$ is a Carter subgroup of $KP/O_2(P)$. Thus $P$ is a Borel
subgroup and the theorem follows from Lemma \ref{CartBorel}. Now if $p\not=
2$, then again $K_G$ is contained in a proper parabolic subgroup $P$ of $G$
such that $O_p(K_G)\leq O_p(P)$ and $O_2(K_G)\leq L$. Then Lemmas
\ref{Syl2InCentrOfFieldAut} and \ref{ConjAutomorphisms} implies that
$H_2\leq O_2(K\cap G)\leq K$. Now Lemma \ref{centUH} implies that
$O_p(K_G)\leq
C_U(H_2)=\{e\}$. Therefore $K\cap G$ is a $2$-gro\-up. By Lemmas
\ref{Syl2InCentrOfFieldAut} and \ref{ConjAutomorphisms} every element  $x\in
A\setminus G$ of odd order such  that $\langle x\rangle\cap G=\{e\}$ centralizes
some Sylow $2$-sub\-gro\-up of $G$. Hence $K$ contains a Sylow $2$-sub\-gro\-up
of $A$,
i.~e.,  $K$ satisfies statement (d) of the theorem.
Therefore $A/G$ is cyclic and we may assume that~${A=\langle G,\zeta
g\rangle\in\Gamma G}$.

Recall that we are in the conditions of Theorem
\ref{CarterSemilinear}, $A=\langle G,\zeta g\rangle$ is supposed to be a
counter example
to the theorem with $\vert O^{p'}(G)\vert$ and $\vert A\vert$ minimal, and $K$
is a Carter
subgroup of $\langle G,\zeta g\rangle$ such that $p$ divides~$\vert K_G\vert$.
We have
that $K=\langle \zeta^kg, K_G\rangle$. Since $\vert O^{p'}(G)\vert
\le \textrm{\bfseries{Cmin}}$, Lemma \ref{HomImageOfCarter} implies that
 $KG/G$ is a Carter subgroup of $\langle G,\zeta g\rangle/G$. Therefore
$\vert\zeta^k\vert=\vert \zeta\vert$, and we may assume that~$k=1$ and
$K=\langle K_G,\zeta g\rangle$.

In view of Lemma \ref{parabolic} there exists a proper $\sigma$- and
$\bar{\zeta}g$- invariant parabolic subgroup $\overline{P}$ of
$\overline{G}$ such that $O_p(K_G)\leq R_u(\overline{P})$ and $K_G\leq
\overline{P}$. In particular, $\overline{P}$ and
$\overline{P}^{\bar{\zeta}}$ are conjugate in $\overline{G}$. Let $\Phi$ be
the root system of $\ov{G}$ and
$\Pi$ be a set of fundamental roots of $\Phi$. In view of
\cite[Proposition~8.3.1]{CarSimpleGrpsLieType},  $\overline{P}$ is conjugate to
some
$\overline {P}_J=\overline{B}\cdot\overline{N}_J\cdot \ov{B}$, where $J$ is a
subset of $\Pi$ and $\overline{N}_J$ is a complete preimage of $W_J$
in $\overline{N}$ under the natural homomorphism
$\overline{N}/\overline{T}\rightarrow W$. Now $\overline{P}_J$ is
$\bar\varphi$-in\-va\-ri\-ant, hence
$\overline{P}_J^{\bar\zeta}=\overline{P}_J^{\bar\gamma^\varepsilon}$ (recall
that $\bar\zeta=\bar\gamma^\varepsilon\varphi^k$ by definition). Consider
the symmetry $\rho$  of the
Dynkin diagram of $\Phi$ corresponding to $\bar\gamma$. Let $\overline{J}$ be
the image of $J$ under
$\rho$. Clearly
$\overline{P}_J^{\bar\gamma}=\overline{P}_{\overline{J}}$. Since
$\overline{P}$ and $\overline{P}^{\bar\zeta}$ are conjugate in
$\overline{G}$ we obtain that $\overline{P}_J$ and
$\overline{P}_J^{\bar\zeta}$ are conjugate in $\overline{G}$. By
\cite[Theorem~8.3.3]{CarSimpleGrpsLieType} it follows that either
$\varepsilon=0$,
or $J=\overline{J}$; i.~e.,  $\overline{P}_J$ is
$\bar\zeta$-in\-va\-ri\-ant.

Now we have that $\overline{P}^{\bar{y}}=\overline{P}_J$ for some $\bar y\in
\ov{G}$. So
$\langle\bar{\zeta}g,\overline{P}\rangle^{\bar{y}}=\langle
(\bar{\zeta}g)^{\bar{y}},
\overline{P}_J\rangle$ and
$\overline{P}_J^{(\bar{\zeta}g)^{\bar{y}}}=\overline{P}_J$. It follows
$$(\bar{\zeta}g)^{\bar{y}}=\bar{y}^{-1}\bar{\zeta}g\bar{y}=
\bar{\zeta}\left(\bar{\zeta}^{-1}\bar{y}^{-1}\bar{\zeta}g\bar{y}\right)=
\bar{\zeta}\cdot h,$$ where
$h=\left(\bar{\zeta}^{-1}\bar{y}^{-1}\bar{\zeta}g\bar{y}\right)\in\overline{G}$.
Since
$\overline{P}_J^{\bar\zeta}=\overline{P}_J=\overline{P}_J^{h^{-1}}$
we obtain that $h\in N_{\overline{G}}(\overline{P}_J)$. By
\cite[Theorem~8.3.3]{CarSimpleGrpsLieType}, $N_{\overline{G}}(\overline{P}
_J)=\overline{P}_J$, thus
$\langle\bar{\zeta}g,\overline{P}\rangle^y=\langle
\bar\zeta,\overline{P}_J\rangle$. Now both $\overline{P}$ and
$\overline{P}_J$ are $\sigma$-in\-va\-ri\-ant. Hence
$\bar{y}\sigma(\bar{y}^{-1})\in
N_{\overline{G}}(\overline{P})=\overline{P}$. Therefore, by
Lang-Steinberg Theorem (Lemma \ref{LangSteinbergTheorem}), we may assume that
$\bar{y}=\sigma(\bar{y})$, i.~e., $\bar{y}\in\overline{G}_\sigma$.
Since $\overline{G}_\sigma=\overline{T}_\sigma\cdot
O^{p'}(\overline{G}_\sigma)$ and $\overline{T}\leq \overline{P}_J$,
then we may assume that $\bar{y}\in O^{p'}(\overline{G}_\sigma)$.
Thus, up to conjugation in $G$, we may assume that $\overline{K}\leq
\langle \bar\zeta,\overline{P}_J\rangle=
\overline{P}_J\leftthreetimes\langle\bar\zeta\rangle$ and $$K\leq\langle
(\overline{P}_J\cap G),\zeta g\rangle=\langle P_J,\zeta g\rangle,$$
in particular, $g\in
(\ov{P}_J)_\sigma$. Further if
$\overline{L}_J=\langle\overline{T}, \overline{X}_r\vert r\in J\cup
-J\rangle$, then $\overline{L}_J$ is a $\sigma$- and
$\bar\zeta$- invariant Levi factor of $\overline{P}_J$ and
$L_J=\overline{L}_J\cap G$ is a $\zeta$-in\-va\-ri\-ant Levi factor
of~$P_J$. Then $L_J^g$ is a $\zeta g$-stab\-le factor Levi of $P_J$. Since all
Levi factors are conjugate under $O_p(P_J)$, we may assume that $L_J$ is a
$\zeta g$-stab\-le Levi factor. Lemma \ref{HomImageOfCarter} implies that
$$KO_p(P_J)/O_p(P_J)=X$$ is a Carter subgroup of
$\langle P_J, \zeta g\rangle/O_p(P_J)$ and
$$K Z(L_J)O_p(P_J)/Z(L_J)O_p(P_J)=\widetilde{X}$$ is a Carter subgroup
of~$\langle P_J,  \zeta g\rangle/Z(L_J)O_p(P_J)$. Recall that
$K=\langle \zeta g,K_G\rangle$, hence, if $v$ and $\tilde{v}$ are
the images of $g$ under the natural homomorphisms
$$\omega:\langle P_J,\zeta g\rangle\rightarrow
\langle L_J,\zeta g\rangle\simeq\langle P_J,\zeta g\rangle/O_p(P_J),$$
$$\tilde\omega:\langle P_J,\zeta g\rangle\rightarrow
\langle P_J,\zeta g\rangle/Z(L_J)O_p(P_J)\simeq
\langle L_J,\zeta g\rangle/Z(L_J),$$ then
$X=\langle\zeta v,K_G^{\omega}\rangle$
and~$\widetilde{X}=\langle\zeta
\tilde{v},K_G^{\tilde\omega}\rangle$. Note that $O_p(P)$ and
$Z(L_J)$ are characteristic subgroups of $P$ and $L_J$ respectively,
hence we may consider $\zeta$ as an automorphism of $L_J\simeq
P/O_p(P)$ and $\widetilde{L}=L_J/Z(L_J)$. Note also that all
non-Abelian composition factors of $P$ are simple groups of Lie type
of order less than $\textrm{\bfseries{Cmin}}$, hence
$\langle P,\zeta g\rangle$ satisfies {\bfseries (C)}.
Thus we may apply Lemma \ref{HomImageOfCarter}
to $\langle \widetilde{L},\zeta g\rangle$, $\langle L,\zeta
g\rangle$, and $\langle P,\zeta g\rangle$.

If $P_J$ is a Borel subgroup of $G$, then the statement of the theorem follows
from Lemma \ref{CartBorel}. So we may assume
that $L_J\not=Z(L_J)$, i.~e., that $P_J$ is not a Borel subgroup of $G$. Then
$L_J=H(G_1\ast\ldots\ast G_k)$, where $G_i$ are subsystem subgroups of $G$,
$k\ge1$, and $H$
is a Cartan subgroup of $G$. Let $\zeta g=(\zeta_2 g_2)\cdot(\zeta_{2'} g_{2'})$
be the product of
$2$- and $2'$- parts of $\zeta g$ (with $g_2,g_{2'}\in(\ov{P}_J)_\zeta$). Now
$\zeta_{2'}=\varphi^k$, for some $k$, is a
field
automorphism (recall that we do not consider the triality automorphism) and it
normalizes each $G_i$, since $\varphi$ normalizes each $G_i$. Moreover, in
view of
Lemma \ref{Syl2InCentrOfFieldAut}, we have that $\zeta_{2'}$ centralizes a Sylow
$2$-sub\-gro\-up of $H$. In particular, it centralizes a Sylow $2$-sub\-gro\-up
of
$Z(L_J)\leq H$. Therefore, every element of odd order of $\langle L_J,
\zeta_{2'}v_{2'}\rangle$ centralizes a Sylow $2$-sub\-gro\-up
of~$Z(L_J)$ (here $v_{2'}$ is the image of $g_{2'}$ under~$\omega$).

Now $\widetilde{L}=(\P G_1\times\ldots\times \P
G_k)\widetilde{H}$, where
$\widetilde{H}=H^{\omega}$ and $\P G_1,\ldots,\P G_k$ are canonical finite
groups of
Lie type with trivial center. Set $M_i=C_{\widetilde{L}}(\P G_i)$, clearly
$M_i=(\P G_1\times
\ldots\times \P G_{i-1}\times \P
G_{i+1}\times\ldots\times \P G_k)C_{\widetilde{H}}(\P G_i)$; denote by $L_i$ the
factor
group
$\widetilde{L}/M_i$ and by $\pi_i$ corresponding natural
homomorphism. Then $L_i$ is a finite group of Lie type and $\P G_i\leq
L_i\leq \widehat{\P G}_i$.

Set $M_{i,j}=C_{\widetilde{L}}(\P G_i\times \P
G_j)$, then $$M_{i,j}=(\P G_1\times
\ldots\times \P G_{i-1}\times \P
G_{i+1}\times\ldots\times \P G_{j-1}\times \P G_{j+1}\times\ldots\times\P
G_k)C_{\widetilde{H}}(\P G_i\times\P G_j);$$ denote by $\pi_{i,j}$
corresponding natural homomorphism $\widetilde{L}\rightarrow
\widetilde{L}/M_{i,j}$. If $M_i$ (respectively $M_{i,j}$) is
$\zeta$-in\-vari\-ant,
then $M_i$ (resp. $M_{i,j}$) is normal in
$\langle\widetilde{L},\zeta\tilde{v}\rangle$ and we denote by $\pi_i$
(resp. $\pi_{i,j}$) the natural homomorphism
$\pi_i:\langle\widetilde{L},\zeta\tilde{v}\rangle\rightarrow
\langle\widetilde{L},\zeta\tilde{v}\rangle/M_i$
($\pi_{i,j}:\langle\widetilde{L},\zeta\tilde{v}\rangle\rightarrow
\langle\widetilde{L},\zeta\tilde{v}\rangle/M_{i,j}$).

Now consider $\zeta$. Since $\zeta^2$ is a field automorphism, there can be
two cases: either $\zeta$
normalizes $\P G_i$, or $\zeta^2$ normalizes $\P G_i$ and
$\P G_i^{\zeta}=\P G_j$ for some $j\not=i$. Consider these two cases
separately.

Let $\zeta$ normalizes $\P G_i$. Then $\zeta$ normalizes $M_i$,
and Lemma \ref{HomImageOfCarter} implies that
$\widetilde{X}^{\pi_i}=K_i$ is a Carter subgroup of
$\langle L_i, (\zeta \tilde{v})^{\pi_i}\rangle$. Since
$\langle L_i, (\zeta \tilde{v})^{\pi_i}\rangle$ is a semilinear group of
Lie type satisfying the conditions of Theorem \ref{CarterSemilinear}
(by definition, $\zeta^2$ is a field automorphism, so we are not in
the conditions of Theorem \ref{CarterTriality}), $\vert L_i
\vert<\vert G\vert$, and $p$ does not divide $\vert K_i\vert$, we
have that $K_i$ contains a Sylow $2$-sub\-gro\-up $Q_i$ of
$\langle L_i, (\zeta \tilde{v})^{\pi_i}\rangle$ (in particular, $p\not=2$) and,
by Lemma \ref{CritSyl2Carter}, the group $\langle L_i, (\zeta
\tilde{v})^{\pi_i}\rangle$ satisfies~{\bfseries(ESyl2)}.

Let $\zeta^2$ normalizes $\P G_i$ and $\P G_i^{\zeta}=\P G_j$.
Then $M_{i,j}$ is normal in
$\langle\widetilde{L},\zeta\tilde{v}\rangle$. We want to show
that $\langle\widetilde{L},\zeta\tilde{v}\rangle^{\pi_{i,j}}$ satisfies
{\bfseries(ESyl2)}.
Since $M_{i,j}$ is a normal subgroup of
$\langle\widetilde{L},\zeta\tilde{v}\rangle$, then,  by Lemma
\ref{HomImageOfCarter}, $(\widetilde{X})^{\pi_{i,j}}$ is a Carter
subgroup of
$\langle\widetilde{L},\zeta\tilde{v}\rangle^{\pi_{i,j}}$. Consider the subgroup
$$\langle (\P G_i)^{\pi_{i,j}}\times(\P
G_j)^{\pi_{i,j}},\widetilde{X}^{\pi_{i,j}}\rangle$$
of $\langle\widetilde{L},\zeta\tilde{v}\rangle^{\pi_{i,j}}$ (note that $(\P
G_i)^{\pi_{i,j}}\simeq\P G_i$ and $(\P G_j)^{\pi_{i,j}}\simeq \P G_i$, and till
the end of this paragraph for brevity we shall identify these groups). Now we
are in the conditions of Lemma
\ref{CarterInGroupOFInducedAutomorphismsWithoutCenter}, namely we have a
finite group $\widetilde{G}=(\widetilde{X})^{\pi_{i,j}}(\P G_i\times\P
G_j)$, where $\P G_i\simeq \P G_j$ has trivial center. Then
$\Aut_{(\widetilde{X})^{\pi_{i,j}}}(\P G_i)\simeq \Aut_{\widetilde{X}}(\P G_i)$
is a Carter subgroup of $\Aut_{\widetilde{G}}(\P G_i)$. Now $\P G_i$ is a
canonical finite group of Lie type and
$$\P G_i\leq \Aut_{\widetilde{G}}(\P G_i)\leq \Aut(\P G_i),$$ i.~e.,
$\Aut_{\widetilde{G}}(\P G_i)$ satisfies to the conditions of Theorem
\ref{CarterSemilinear} (by construction $\zeta^2$ is a field automorphism and
so we are not in the conditions of Theorem \ref{CarterTriality}) and
$(\widetilde{X})^{\pi_{i,j}}\cap (\P G_i\times \P G_j)$ is not divisible by the
characteristic. By induction, $\Aut_{(\widetilde{X})^{\pi_{i,j}}}(\P G_i)$
contains a Sylow $2$-sub\-gro\-up of $\Aut_{\widetilde{G}}(\P G_i)$ (in
particular,
$p\not=2$). The same arguments show that $\Aut_{\widetilde{X}}(\P G_j)$
contains a Sylow $2$-sub\-gro\-up of $\Aut_{\widetilde{G}}(\P G_j)$. Therefore,
$\Aut_{\widetilde{G}}(\P G_i)$ and $\Aut_{\widetilde{G}}(\P G_j)$ satisfy
{\bfseries(ESyl2)}. Since $\Aut_{\widetilde{G}}(\P G_i)\leq
\Aut_{\langle\widetilde{L},\zeta\tilde{v}\rangle^{\pi_{i,j}}} (\P G_i)$ and
$\Aut_{\widetilde{G}}(\P G_j)\leq
\Aut_{\langle\widetilde{L},\zeta\tilde{v}\rangle^{\pi_{i,j}}}(\P G_j)$, Lemmas
\ref{Esyl2InnDiagExtension} and \ref{ESyl2InhFieldGraph} imply that groups of
induced automorphisms
$\Aut_{\langle\widetilde{L},\zeta\tilde{v}\rangle^{\pi_{i,j}}}(\P G_i)$ and
$\Aut_{\langle\widetilde{L},\zeta\tilde{v}\rangle^{\pi_{i,j}}}(\P G_j)$ satisfy
{\bfseries(ESyl2)}. Consider
$N_{\langle\widetilde{L},\zeta\tilde{v}\rangle^{\pi_{i,j}}}(\P G_i)$ and
$N_{\langle\widetilde{L},\zeta\tilde{v}\rangle^{\pi_{i,j}}}(\P G_j)$. Since
$$\vert{\langle\widetilde{L},\zeta\tilde{v}\rangle^{\pi_{i,j}}}:
N_{\langle\widetilde{L},\zeta\tilde{v}\rangle^{\pi_{i,j}}} (\P G_i)\vert=
\vert{\langle\widetilde{L},\zeta\tilde{v}\rangle^{\pi_{i,j}}}:
N_{\langle\widetilde{L},\zeta\tilde{v}\rangle^{\pi_{i,j}}} (\P G_j)\vert=2,$$
it is easy to see that for every element $h$ of
$\langle\widetilde{L},\zeta\tilde{v}\rangle^{\pi_{i,j}}$ the equality of cosets
$h N_{\langle\widetilde{L},\zeta\tilde{v}\rangle^{\pi_{i,j}}}(\P G_i)=h
N_{\langle\widetilde{L},\zeta\tilde{v}\rangle^{\pi_{i,j}}}(\P G_j)$ holds, it
follows that $N_{\langle\widetilde{L},\zeta\tilde{v}\rangle^{\pi_{i,j}}}(\P
G_i)=N_{\langle\widetilde{L},\zeta\tilde{v}\rangle^{\pi_{i,j}}}(\P G_j)$. By
construction $C_{\langle\widetilde{L},\zeta\tilde{v}\rangle^{\pi_{i,j}}} (\P
G_i)\cap
C_{\langle\widetilde{L},\zeta\tilde{v}\rangle^{\pi_{i,j}}}(\P
G_j)=\{e\}$, so Lemma \ref{Syl2centrcomposit} (with
$C_{\langle\widetilde{L},\zeta\tilde{v}\rangle^{\pi_{i,j}}}(\P G_i)$ and
$C_{\langle\widetilde{L},\zeta\tilde{v}\rangle^{\pi_{i,j}}}(\P G_j)$ as normal
subgroups) implies that the normalizer
$N_{\langle\widetilde{L},\zeta\tilde{v}\rangle^{\pi_{i,j}}}(\P G_i)$ satisfies
{\bfseries(ESyl2)}. Now $\vert
{\langle\widetilde{L},\zeta\tilde{v}\rangle^{\pi_{i,j}}}:
N_{\langle\widetilde{L} ,\zeta\tilde{v}\rangle^{\pi_{i,j}}} (\P
G_1)\vert=2$, thus Lemma \ref{InhBy2-ext} implies that
$\langle\widetilde{L},\zeta\tilde{v}\rangle^{\pi_{i,j}}$
satisfies~{\bfseries(ESyl2)}.

Now we shall show that $\langle L_J,\zeta v\rangle$ satisfies
{\bfseries(ESyl2)}. Since $\widetilde{L}\not=\{e\}$, then, as we noted above,
$p\not=2$. Let $Q$ be a Sylow $2$-sub\-gro\-up of $\langle L_J,\zeta v\rangle$.
Consider an element $x\in N_{\langle L_J,\zeta v\rangle}(Q)$ of odd order. We
need to prove that $x$ centralizes $Q$. As we noted above, every element of odd
order of $\langle L_J,\zeta v\rangle$ centralizes $Q\cap Z(L_J)$, hence, if
$\tilde{x}=x^{\tilde\omega}$
centralizes $\widetilde{Q}=Q^{\tilde\omega}\simeq Q/(Q\cap Z(L_J))$, then $x$
centralizes $Q$. Now either $M_i$ is
normal in $\langle \widetilde{L},\zeta \tilde{v}\rangle$, or $M_{i,j}$ is
normal in $\langle \widetilde{L},\zeta \tilde{v}\rangle$ and
$\left(\cap_iM_i\right)\bigcap\left(\cap_{i,j}M_{i,j}\right)=\{e\}$. Moreover,
as we proved above, $x^{\pi_i}$ centralizes $\widetilde{Q}M_i/M_i$, and
$x^{\pi_{i,j}}$ centralizes $\widetilde{Q} M_{i,j}/M_{i,j}$. By Lemma
\ref{Syl2centrcomposit} (with normal subgroup $M_i$ and $M_{i,j}$) we obtain
that $\tilde{x}$ centralizes~$\widetilde{Q}$.

Thus $\langle L_J,\zeta v\rangle$ satisfies {\bfseries(ESyl2)} and by
Lemma \ref{CritSyl2Carter} there exists a Carter subgroup $F$ of $\langle
L_J,\zeta v\rangle$ containing $Q$.  Since $\langle L_J,\zeta v\rangle$
satisfies {\bfseries(C)}, Theorem \ref{ConjugacyCriterion} implies
that $X=K^\omega$
and $F$ are conjugate, i.~e., $X$ contains a Sylow $2$-sub\-gro\-up of $\langle
L_J,\zeta v\rangle$ and, up to conjugation in $\langle P_J,\zeta v\rangle$, $K$
contains a Sylow $2$-sub\-gro\-up of $\langle P_J,\zeta v\rangle$. In
particular, a
Sylow $2$-sub\-gro\-up $Q_1$ of a Cartan subgroup $H$ is in $K$ and $Q_1$
centralizes $K\cap O_p(P_J)\not=\{e\}$; a
contradiction with Lemma~\ref{centUH}.

\subsection[Carter subgroups without unipotent radical]{Carter
subgroups of order not divisible by characteristic}

Again we are in the conditions of Theorem \ref{CarterSemilinear}. As we noted
in the previous section, for every group $A$ satisfying conditions of Theorem
\ref{CarterSemilinear}, the factor group $A/G$ is Abelian and, for some natural
$t$ is isomorphic to a subgroup of $\Z_2\times\Z_t$. If the factor group $A/G$
is not cyclic, then $O^{p'}(G)$ is split and $A$ contains an element $\tau a$,
where $\tau$ is a graph automorphism of $O^{p'}(G)$ and $a\in\ov{G}_\sigma$.
Thus, if $A/G$ is not cyclic, or $\Phi(\ov{G})\not=A_n,D_{2n+1}, E_6$, then
by Lemmas \ref{omnibus} and \ref{ConjInverseInGraph} every semisimple
element of $G$ is conjugate to its inverse. By Lemma \ref{power} we obtain that
$K_G=K\cap G$ is a $2$-gro\-up. In the conditions of Theorem
\ref{CarterSemilinear} the group $A/G$ is Abelian and, if $\ov{A}_1$ is a Hall
$2'$-sub\-gro\-up of $A/G$, then $\ov{A}_1$ is cyclic. Let $x$ be the preimage
of
the generating element of $\ov{A}_1$ taken in $K$. Then $\langle x\rangle\cap
G\leq \langle x\rangle\cap\ov{G}_\sigma\leq K\cap \ov{G}_\sigma=K\cap (A\cap
\ov{G}_\sigma)=K\cap G$. As we noted above, $K\cap G$ is a $2$-gro\-up, hence
$\langle x\rangle \cap \ov{G}_\sigma=\{e\}$. By Lemma \ref{ConjAutomorphisms},
the element $x$ under $\ov{G}_\sigma$ is conjugate to a field automorphism of
odd order and by Lemma \ref{Syl2InCentrOfFieldAut}, the element $x$ centralizes
a Sylow $2$-sub\-gro\-up of $G$ (in particular, $p\not=2$) and, since $A/G$ is
Abelian, Lemma \ref{Syl2centrcomposit} implies that $K$ contains a Sylow
$2$-sub\-gro\-up of $A$. Thus Theorem \ref{CarterSemilinear} is true in this
case.
So we may assume that
$A=\langle G, \zeta g\rangle$  is a semilinear group of Lie type,
$K=\langle \zeta^kg,K_G\rangle$ is a Carter
subgroup of $A$, and $\Phi(\ov{G})\in\{A_n, D_{2n+1},
E_6\}$. Like in the
previous section we may assume
that~$k=1$. Since $G_\zeta$ is nontrivial, then the centralizer $C_G(\zeta g)$
is also nontrivial, we have that $K_G$ is also
nontrivial. Since $G_\zeta$ is nontrivial, then the centralizer $C_G(\zeta g)$
is nontrivial, so $K_G$ is also nontrivial. Therefore $Z(K)\cap K_G$ is
nontrivial. Consider an element
$x\in Z(K)\cap K_G$ of prime order. Then $K\in C_{A}(x)=\langle \zeta
g,C_G(x)\rangle$. Now $C_{\overline{G}}(x)^0=\overline{C}$ is a
connected $\sigma$-stab\-le reductive subgroup of maximal rank of
$\overline{G}$. Moreover
$\overline{C}$ is a characteristic subgroup
of $C_{\overline{G}}(x)$ and $C_{\ov{G}}(x)/\ov{C}$ is isomorphic to a
subgroup of $\Delta$ (see \cite[Proposition~2.10]{Hu2ConjClasses}). Thus $K$ is
contained
in $\langle K,
C\rangle$, where $C=\ov{C}\cap G$. Moreover, by Lemma \ref{equivnormalizer}, the
subgroup $C=\ov{C}\cap
G=T(G_1\ast\ldots\ast G_m)$ is normal in $C_{A}(x)$ and $K_GC/C$ is isomorphic
to a subgroup of $\Delta$. Assume that $\vert
K_G\vert$ is not divisible by~$2$.

If $m=0$, then $C=T=Z(C)$ is a maximal torus. Then $\ov{T}$
is $\bar\zeta g$-stab\-le. In view of Lemma \ref{NormOfRegularElementIsNotCentr}
we obtain that $N_A(C_A(x))\not=C_A(x)$. Since $C_A(x)$ is solvable in this case
this gives a contradiction with Lemma~\ref{power}.

If $m\ge1$, then $Z(C)$  and $G_1\ast\ldots\ast G_m$ are normal subgroups of
$\langle K, C\rangle$. Hence we may consider $\widetilde{G}=\langle
K,G_1\ast\ldots\ast G_m\ast Z(C)\rangle/Z(C)\leq \langle K,C\rangle/Z(C)$. Then
$\widetilde{G}=\widetilde{K}(\P G_1\times\ldots\times \P G_m)$,
where $\widetilde{K}=KZ(C)/Z(C)$ is a Carter subgroup of
$\widetilde{G}$ (see Lemma \ref{HomImageOfCarter}) and $Z(\P G_i)$
is trivial. Now $\widetilde{K}$ acts by conjugation on $\{\P
G_1,\ldots, \P G_ m\}$ and without lost of generality we may assume that $\{\P
G_1,\ldots,\P G_m\}$ is a $\widetilde{K}$-or\-bit. Thus we are in the
condition of Lemma \ref{CarterInGroupOFInducedAutomorphismsWithoutCenter} and
$\Aut_{\widetilde{K}}(\P G_1)$ is a Carter subgroup of
$\Aut_{\widetilde{G}}(\P G_1)$. Moreover $\vert\widetilde{K}\cap \P
G_1\times\ldots\times\P G_m\vert$ is not divisible by the
characteristic. By induction we have that either $\Aut_{\widetilde{K}}(\P
G_1)$ contains a Sylow $2$-sub\-gro\-up of $\Aut_{\widetilde{G}}(\P
G_1)$, or $\Aut_{\widetilde{G}}(\P G_1)$ satisfies to the
conditions of Theorem \ref{CarterTriality} and $\Aut_{\widetilde{G}}(\P
G_1)\cap \P G_1$ is a nontrivial $2$-gro\-up, in particular $p$ is odd. In any
case $\vert K\cap G\vert$ is divisible by $2$ that contradicts our
assumption. Therefore the order $\vert K_G\vert$ is even and we may assume that
$x\in Z(K)\cap K_G$ is an involution.

Write $\zeta g=\zeta_2 g_1\cdot \zeta_{2'} g_2$, where
$\zeta_2g_1$ is the $2$-part and $\zeta_{2'}g_2$ is  the $2'$-part of $\zeta g$.
By Lemma \ref{Syl2InCentrOfFieldAut} the element
$\zeta_{2'}$ centralizes a Sylow $2$-sub\-gro\-up $Q_G$ of $G$, so we may assume
that the order of  $g_2$ is odd. Up to conjugation in $G$ we may assume that
$\zeta_{2'}$ centralizes a Sylow $2$-sub\-gro\-up of $K_G$. In particular,
$\zeta_{2'}$ centralizes $x$. Let $Q$ be a Sylow $2$-sub\-gro\-up of $C_G(x)$.
Then
there exists $y\in G$ such that $Q^y\leq Q_G$. Substituting the subgroup
$K$ by its conjugate $K^y$, we may assume that $\zeta_{2'}$ centralizes a Sylow
$2$-sub\-gro\-up of~$C_G(x)$. Since $\zeta_{2'} g_2$ centralizes $x$, we obtain
that $g_2\in C_{\ov{G}_\sigma}(x)$. Moreover, by Lemma \ref{CentrOfInvolution}
it follows that $g_2\in C_{\ov{G}}(x)^0$.  In particular,  $g_2$ normalizes
each $G_i$ and centralizes $Z(C)$ and~$Z(C_G(x))$.

Note that $\zeta_{2'}$ normalizes each $G_i$ and centralizes a Sylow
$2$-sub\-gro\-up of $Z(C_G(x))$ (recall that $\zeta_{2'}$ centralizes a
Sylow $2$-sub\-gro\-up of $C_G(x)$). Indeed,  $\zeta_{2'}$ normalizes $C$, hence
normalizes characteristic subgroups $O^{p'}(C)=G_1\ast\ldots\ast G_m$ and $Z(C)$
of $C$. So we may consider the induced automorphism $\zeta_{2'}$ of
$$O^{p'}(C)/(Z(C)\cap O^{p'}(C)=\P G_1\times\ldots\times \P G_m.$$ Since
each $\P G_i$ has trivial center and can not be written as a direct product
of proper subgroups, corollary from Krull-Remak-Schmidt
Theorem  \cite[3.3.10]{Rob} implies that $\zeta_{2'}$ permutes distinct  $\P
G_i$. Since $\zeta_{2'}$ centralizes a Sylow $2$-sub\-gro\-up of $C_G(x)$ and
$C\unlhd C_G(x)$, then $\zeta_{2'}$ centralizes a Sylow $2$-sub\-gro\-up of 
$C$,
hence centralizes a Sylow  $2$-sub\-gro\-up $Q_1\times\ldots\times Q_m$ of
$\P G_1\times\ldots\times \P G_m$, where  $Q_i$ is a Sylow $2$-sub\-gro\-up
of $\P G_i$. If $\zeta_{2'}$ would induce a nontrivial permutation on
the set $\{\P G_1,\ldots,\P G_m\}$, then in would induce a nontrivial
permutation on $\{Q_1,\ldots,Q_m\}$. Since each $Q_i$ is nontrivial, this is
impossible. Thus every element of odd order of
$\langle K,C\rangle$ centralizes a Sylow $2$-sub\-gro\-up of~$Z(C)$ and
normalizes
each~$G_i$.

If $\Phi(\ov{G})=E_6$, then by Lemma \ref{CentrOfInvolution} the centralizer of
every involution of $G$ in  $\ov{G}$ is connected. By Lemma
\ref{InvolutionsAndTori} every involution of $G$ is contained in a maximal torus
$T$ such that $N(G,T)/T\simeq W$, where $W$ is a Weyl group of
$\ov{G}$.  $\ov{C}$ is welknown to be generated by the torus $\ov{T}$ and
$\overline{T}$-ro\-ot subgroups. Write $\ov
C=\ov{T}(\ov{G}_1\ast\ldots\ast\ov{G}_k)$. Since $\ov{T}_\sigma$ either is
obtained from a maximal split torus $\ov{H}$ by twisting with an element $w_0$
of order $2$, or is equal to $\ov{H}$, and each field automorphism acts
trivially on the factor group  $N_{\ov{G}}(\ov{H})/\ov{H}$, then
$\bar{\zeta}_{2'}$ normalizes every subgroup $\ov{G}_i$. So, if
$\Phi(\ov{G}_i)=D_4$, then  $\bar{\zeta}_{2'}$ induces a field (but not a graph
or a graph-field) automorphism of $\ov{G}_i$. Moreover, since $\sigma$ acts
trivially on the factor group  $N_{\ov{G}}(\ov{T})/\ov{T}$ (see Lemma
\ref{Syl2InCentrOfFieldAut}), then  \cite[Proposition~6]{Car5CentSemisimpleLie}
implies
that $\sigma$ normalizes each $\ov{G}_i$. Therefore, none of $G_i$ is
isomorphic to ${}^3D_4(q^3)$. If $\Phi(\ov{G})$ coincides with $A_n$ or  $D_n$,
then  \cite[Propositions 7, 8, 10]{Ca2CentsSemisimpleClassical} imply that
none of $G_i$ is isomorphic
to ${}^3D_4(q^3)$. Therefore in any case none of $G_i$ is isomorphic
to~${{}^3D_4(q^3)}$. Moreover Lemma \ref{CentrOfInvolution} implies that
$\vert K_G:(K_G\cap C)\vert$ divides $\vert C_{\ov
G}(x)/C_{\ov{G}}(x)^0\vert$ and $C_{\ov{G}}(x)/C_{\ov{G}}(x)^0$ is a
$2$-gro\-up. In  \cite{Ca2CentsSemisimpleClassical} it is proven that if a root
system
$\Phi$ has type  $D_n$ and $\Psi$ is its subsystem of type $D_4$, then none
element from $N_{W(\Phi)}(W(\Psi))$ induces a symmetry of order $3$ of the
Dynkin diagram of $\Psi$. Since $\zeta^2$ is a field automorphism, lack of a
symmetry of order  $3$ together with 
\cite[Proposition~6]{Car5CentSemisimpleLie}
implies that for
each $G_i$ the automorphism $\zeta_{2'}$ is field (but not graph o
graph-field). Therefore the group of induced automorphisms  $\langle
\Aut_{\widetilde{K}}(\P G_i),\P G_i\rangle$ satisfies to the conditions of
Theorem~\ref{CarterSemilinear} for all~$i$.

Now consider $\widetilde{G}=\widetilde{K}(\P G_1\times\ldots\times\P
G_m)\leq \langle K,C\rangle/Z(C)$ (probably, $m=0$), where
$\widetilde{K}=KZ(C)/Z(C)$ is a Carter subgroup of $\widetilde{G}$ (see Lemma
\ref{HomImageOfCarter}) and, for all $i$, $Z(\P G_i)=\{e\}$. By Lemma
\ref{CarterInGroupOFInducedAutomorphismsWithoutCenter} we have that
$\Aut_{\widetilde{K}}(\P G_1)$ is a Carter subgroup of
$\Aut_{\widetilde{G}}(\P G_1)$. Since $\P G_1$ is a finite group of
Lie type satisfying Theorem \ref{CarterSemilinear}, by induction we
obtain that $\Aut_{\widetilde{G}}(\P G_1)$ satisfies {\bfseries(ESyl2)}.
Similarly we have
that $\Aut_{\widetilde{G}}(\P G_i)$ satisfies {\bfseries(ESyl2)} for all $i$.
Since $$\Aut_{\langle K,C\rangle/Z(C)}(\P G_i)\geq \Aut_{\widetilde{G}}(\P
G_i),$$ Lemmas \ref{Esyl2InnDiagExtension} and \ref{ESyl2InhFieldGraph} imply
that $\Aut_{\langle K,C\rangle/Z(C)}(\P G_i)$ satisfies {\bfseries(ESyl2)}.
Since
$C_{\langle K,C\rangle/Z(C)}(\P G_1\times\ldots\times \P G_m)=\{e\}$, Lemma
\ref{Syl2centrcomposit} with normal subgroups $C_{\langle K,C\rangle/Z(C)}(\P
G_1)\cap N_{\langle K,C\rangle/Z(C)}(\P
G_1),\ldots,\linebreak C_{\langle K,C\rangle/Z(C)}(\P G_m)\cap N_{\langle
K,C\rangle/Z(C)}(\P G_1)$ implies that
$N_{\langle K,C\rangle/Z(C)}(\P G_1)$ satisfies {\bfseries(ESyl2)}. Now $$\vert
\langle K,C\rangle/Z(C):N_{\langle K,C\rangle/Z(C)}(\P G_1)\vert=2^t,$$ and
each element of odd order of $\langle K,C\rangle/Z(C)$ normalizes $\P G_1$,
thus, by Lemma
\ref{InhBy2-ext}, we obtain that the factor group $\langle K,C\rangle/Z(C)$
satisfies {\bfseries(ESyl2)} and, by Lemma \ref{Syl2centrcomposit} $\langle
K,C\rangle$ satisfies {\bfseries(ESyl2)}. Since $\vert\P G_i\vert<
\textrm{\bfseries{Cmin}}$, then $\langle K,C\rangle$ satisfy {\bfseries (C)}. By
Lemma \ref{CritSyl2Carter} we obtain that there exists a Carter
subgroup  $F$ of $\langle K,C\rangle$ containing  a Sylow
$2$-sub\-gro\-up of $\langle K,C\rangle$. By Theorem
\ref{ConjugacyCriterion} subgroups $F$ and $K$ are conjugate in $\langle
K,C\rangle$, thus
$K$ contains a Sylow $2$-sub\-gro\-up $Q$ of $\langle K,C\rangle$. Since $\vert
C_G(x):C\vert$ is a power of $2$ and $\langle K,C\rangle$ normalizes $C_G(x)$,
we obtain that $\vert \langle K,C_G(x)\rangle:\langle K,C\rangle\vert$ is a
power of $2$. Moreover by construction each element of odd order of $\langle
K,C_G(x)\rangle$ is in $\langle K,C\rangle$. Thus by Lemma \ref{InhBy2-ext}
$\langle K, C_G(x)\rangle$
satisfies {\bfseries(ESyl2)} and $K$ contains a Sylow $2$-sub\-gro\-up $Q$
of~${\langle K,C_G(x)\rangle}$.

Let $\Gamma  Q$ be a Sylow $2$-sub\-gro\-up of $\langle G,\zeta g\rangle$
containing $Q$ and $t\in Z(\Gamma Q)\cap G$. Then $t\in C_G(x)$, hence, $t\in
Z(Q)$ and $t\in Z(K)$. Thus we may substitute
$x$ by $t$ in arguments above and obtain that $Q=\Gamma Q$, i.~e., $K$ contains
a Sylow $2$-sub\-gro\-up of~$\langle G,\zeta g\rangle$, than completes the proof
of
Theorem~\ref{CarterSemilinear}.

\subsection{Carter subgroups of finite groups are conjugate}

Before we formulate the main theorem, note a corollary of
Theorem~\ref{CarterSemilinear}.

\begin{cor}\label{CminInfinity}
{\em$\textrm{\bfseries{Cmin}}=\infty$}, i.~e.~${\mathcal{A}=\varnothing}$.
\end{cor}

\begin{proof}
Indeed, let $\mathcal{A}\not=\varnothing$ and $A\in\mathcal{A}$ is such that
the equality $\vert F^\ast(A)\vert=\textrm{\bfseries{Cmin}}$ holds. Since 
$F^\ast(A)=O^{p'}(\ov{G}_\sigma)$ for an adjoint simple connected linear
algebraic group $\ov{G}$ and a Frobenius map $\sigma$, denote the intersection
$A\cap \ov{G}_\sigma$ by $G$. As we noted in the beginning of subsection  6.1,
we may assume that  $A=K F^\ast(A)=KG$. Therefore $A$  satisfies either to
the conditions of Theorem \ref{CarterTriality}, or to the conditions of Theorem
\ref{CarterSemilinear}. In both cases we have proved that Carter subgroups of
$A$ are conjugate, that contradicts to the choice of~$A$.
\end{proof}

In order to state the main theorem without using of the classification of
finite simple groups, we give the following definition. A finite group is said
to be a {\em $K$-gro\-up}\index{groupК@$K$-group} if all its non-Abelian
composition factors are known
simple groups.

\begin{ttt}\label{FinalConjugacy}
{\em (Main  Theorem)} Let $G$ be a finite $K$-gro\-up. Then Carter subgroups of
$G$ are conjugate.
\end{ttt}

\begin{proof}
By Theorems \ref{main}, \ref{CarterInClassicalGroups},
\ref{sympcarter}, \ref{CarterTriality}, and  \ref{CarterSemilinear} of the
present paper, and also by \cite{DT2CartSbgrpsPGL} we obtain that for each
known simple group  $S$ and each nilpotent subgroup $N$ of a group of its
automorphisms, Carter subgroups of
$\langle N,S\rangle$ are conjugate. So $G$ satisfies
{\bfseries (C)}. Hence by Theorem \ref{ConjugacyCriterion}, Carter subgroups of
$G$ are conjugate.
\end{proof}

From Lemma \ref{HomImageOfCarter} and  Main Theorem \ref{FinalConjugacy}
it follows that a homomorphic image of a Carter subgroup is a Carter subgroup.

\begin{ttt}\label{HomImageOfCarterFinal}
Let $G$ be a finite $K$-gro\-up, $H$ a Carter subgroup of $G$, and
$N$ a normal subgroup of $G$. Then  $HN/N$ is a Carter subgroup of~$G/N$.
\end{ttt}

\section{Existency criterion}

\subsection{Brief review Краткий обзор результатов параграфа}

In this section we shall obtain a criterion of existence of Carter subgroups in
a finite group in terms of its normal series. Note that there exist
finite groups without Carter subgroups, a minimal counter example
is~$\mathrm{Alt}_5$. We shall construct an example showing that an essential
improvement of the criterion is impossible. At the end of the section, for
convenience of the reader, we assemble the classification of Carter subgroups
in finite almost simple groups, that is obtained in the present paper.

Recall that in view of Theorem \ref{FinalConjugacy} in every almost simple
group with known simple socle Carter subgroups are conjugate. Thus, modulo the
classification of finite simple groups, in every finite group Carter subgroups
are conjugate. In this section by a finite group we always mean a finite group
satisfying  {\bfseries(C)}, thus the results of the section do not depend on
the classification of finite simple groups.

\begin{df}\label{DefinitionOfExistence}
Let 
$G=G_0\geq G_1\geq\ldots\geq G_n=\{e\}$ be a chief series of
$G$ (recall that $G$ is assumed to satisfy {\bfseries(C)}). Then
$G_i/G_{i+1}=T_{i,1}\times\ldots\times T_{i,k_i}$, where
$T_{i,1}\simeq\ldots\simeq T_{i,k_i}\simeq T_i$ and $T_i$ is a simple group.
If $i\ge 1$, then denote by $\ov{K}_i$ a Carter subgroup of $G/G_i$
(if it exists) and by $K_i$ its complete preimage in
$G/G_{i+1}$. If  $i=0$, then $\ov{K}_0=\{e\}$ and $K_0=G/G_1$ (note that
$\ov{K}_0$ always exists). A finite group $G$ is said to satisfy  {\bfseries
(E)}\index{condition!of existence {\bfseries (E)}}\glossary{E@{\bfseries(E)}},
if for each $i,j$, either $\ov{K}_i$ does not exists, or $\Aut_{K_i}(T_{i,j})$
contains a Carter subgroup.
\end{df}

By Theorem \ref{CriterionOfExistence} and Theorem 
\ref{HomImageOfCarterFinal} it follows that if a finite group satisfies
{\bfseries (E)}, then, for every $i$, subgroup $\ov{K}_i$ exists, so the first
part of condition {\bfseries (E)} is never satisfied. Recall that by Theorem 
\ref{HomImageOfCarterFinal} a homomorphic image of a Carter subgroup 
is a Carter subgroup. We shall constantly use this fact.

\subsection{Criterion}

Below we shall need an additional information on the structure of Carter
subgroups in groups of special type. Let   $A'$ be a group with a normal
subgroup $T'$. Consider the direct product 
$A_1\times\ldots\times A_k$, where  $A_1\simeq\ldots\simeq A_k\simeq A'$, and
its normal subgroup  $T=T_1\times\ldots\times T_k$, where $T_1\simeq
\ldots\simeq T_k\simeq T'$. Consider the symmetric group $\Sym_k$, acting on
$A_1\times\ldots\times A_k$ by $A_i^s=A_{i^s}$, for all $s\in S$ and define
$X$ to be equal to a semidirect product 
$\left(A_1\times\ldots\times A_k\right)\leftthreetimes \Sym_k$ (permutation
wreath product of $A'$ and $\Sym_k$). Denote by $A$ the direct product
$A_1\times\ldots\times A_k$ and by $\pi_i$ the projection
$\pi_i:A\rightarrow A_i$. In these notations the following lemma holds.

\begin{lem}\label{CarterInSubDirectProductOfSoluble}
Let $G$ be a subgroup of $X$ such that $T\leq G$, $G/(G\cap T)$ is nilpotent,
and $(G\cap A)^{\pi_i}=A_i$. Assume also that $A$ is solvable. Let $K$ be a 
Carter subgroup of~$G$.

\noindent Then $(K\cap A)^{\pi_i}$ is a Carter subgroup of~$A_i$.
\end{lem}

\begin{proof}
Assume that the statement is false and let $G$ be a counter example of
minimal order with $k$ minimal. Then $S=G/(G\cap A)$ is transitive and
primitive. Indeed, if $S$ is not transitive, then $S\leq
\Sym_{k_1}\times\Sym_{k-k_1}$, hence $G\leq G_1\times G_2$. If we denote
by $\psi_i:G\rightarrow G_i$ the natural homomorphism, then $G^{\psi_i}=G_i$
satisfies conditions of the lemma and $K^{\psi_i}=K_i$ is a Carter subgroup of
$G_i$. Clearly $(G\cap A)^{\pi_j}=(G_i\cap A^{\psi_i})^{\pi_j}$, where 
$i=1$ if $j\in\{1,\ldots,k_1\}$ and $i=2$ if $j\in\{k_1+1,\ldots,k\}$, i.~e.,
the following diagrams are commutative:
$$
\xymatrix{
G\cap A\ar[rr]^{\textstyle{\pi_j}}\ar[dr]^{\textstyle{\psi_1}}&&A_j,&G\cap
A\ar[rr]^{\textstyle{\pi_j}}\ar[dr]^{\textstyle{\psi_2}}&&A_j.\\
&G_1\cap A^{\psi_1}\ar[ur]^{\textstyle{\pi_j}}&&&G_2\cap
A^{\psi_2}\ar[ur]^{\textstyle{\pi_j}}&
}
$$
Thus we obtain the statement by induction. If $S$ is transitive, but
is not primitive, let  $$\Omega_1=\{T_1,\ldots,T_m\},
\Omega_2=\{T_{m+1},\ldots,T_{2m}\}, \ldots,
\Omega_l=\{T_{(l-1)m+1},\ldots,T_{lm}\}$$ be a system of imprimitivity. Then it
contains a nontrivial nontransitive normal subgroup 
$$F'\leq \underbrace{\Sym_{m}\times\ldots\times \Sym_{m}}_{l\text{ times}},$$
where $k=m\cdot l$. Consider a complete preimage $F$ of $F'$ in $X$. Then
$G\cap F\leq F_1\times\ldots\times F_l$. Denote by $\psi_i:F\rightarrow 
F_i$ a natural projection, then $(G\cap F)^{\psi_i}=F_i$. Note that all 
$F_i$-s satisfy conditions of the lemma and, if we define
$T_i'=T_{(i-1)m+1}\times\ldots\times T_{im}$, then   $G$ satisfies 
conditions of the lemma with  $T'=T'_1\times\ldots\times T'_l$ and $A'=F$. By
induction we have that $(K\cap F)^{\psi_i}$ is a Carter subgroup of
$F_i$ and, if $j\in\{m\cdot(i-1)+1,\ldots m\cdot i\}$, then  $\left((K\cap
F)^{\psi_i}\cap A^{\psi_i}\right)^{\pi_j}$ is a Carter subgroup of $A_j$. Since
$(G\cap A)^{\pi_j}=\left((K\cap F)^{\psi_i}\cap A^{\psi_i}\right)^{\pi_j}$ (for
suitable $i$), we get the statement by induction.

Let $Y'$ be a minimal normal subgroup of $G$, contained in $T$ (if $Y'$
is trivial, then $T$ is trivial and we have nothing to prove, since $G$ is
nilpotent in this case). Thus $Y'$ is a normal elementary Abelian $p$-gro\-up.
Let $Y_i=(Y')^{\pi_i}$, then $Y=Y_1\times\ldots\times Y_k$ is a nontrivial
normal subgroup of $G$ ($Y$ is a subgroup of $G$ since $T\leq G$). Let 
$\bar{\pi}_i:(G\cap A)\rightarrow A_i/Y_i=\ov{A}_i$ be a projection,
corresponding to $\pi_i$. Denote by $\ov{K}=KY/Y$ a corresponding Carter
subgroup of  $\ov{G}=G/Y$. Then $\ov{G}$ satisfies conditions of the Lemma. By
induction,  $(\ov{K}\cap\ov{A})^{\bar{\pi}_i}$ is a Carter subgroup of
$\ov{A}_i$. Let $K_1$ be a complete preimage of $\ov{K}$ in $G$, and let $Q$ be
a Hall $p'$-sub\-gro\-up of $K_1$. Then $(Q\cap A)^{\pi_i}$ is a Hall 
$p'$-sub\-gro\-up of $(K_1\cap A)^{\pi_i}$. In view of the proof of
\cite[Theorem~20.1.4]{KarMer}, we obtain that $K=N_{K_1}(Q)$ is a Carter
subgroup of $G$ and $(N_{K_1\cap A}(Q\cap A))^{\pi_i}$ is a Carter subgroup of
$A_i$. Thus we need to show that $(N_{K_1\cap A}(Q\cap
A))^{\pi_i}=(N_{K_1\cap S}(Q))^{\pi_i}$. By induction the equality 
$(N_{\ov{K}\cap \ov{A}}(\ov{A}\cap \ov{Q}))^{\bar{\pi}_i}=(N_{\ov{K}\cap
\ov{G}}(\ov{Q}))^{\bar{\pi}_i}$ holds. Thus we need to prove that 
$(N_Y(Q\cap A))^{\pi_i}=(N_Y(Q))^{\pi_i}$. Note also that 
$(N_Y(Q\cap A))^{\pi_i}\leq N_{Y_i}((Q\cap A)^{\pi_i})$.

Since $S$ is transitive and primitive subgroup of $\Sym_k$, then $k=r$ is a
prime and $S=\langle s\rangle$ is cyclic. If $r=p$, then $Q\cap A=Q$ and we
have nothing to prove. Otherwise let $h$ be an $r$-ele\-ment of $K$, generating
$S$ modulo $K\cap A$. Clearly $Q=(Q\cap A)\langle h\rangle$. Let
$t\in Y_i$ be an element of $N_{Y_i}((Q\cap A)^{\pi_i})$. Then
$(t\cdot t^h\cdot\ldots\cdot t^{h^{r-1}})\in N_Y(Q)$ and
$t^{\pi_i}= (t\cdot t^h\cdot\ldots\cdot
t^{h^{r-1}})^{\pi_i}$, hence $(N_Y(Q\cap A))^{\pi_i}\leq
N_{Y_i}((Q\cap A)^{\pi_i})\leq (N_Y(Q))^{\pi_i}\leq {(N_Y(Q\cap A))^{\pi_i}}$.
\end{proof}

\begin{ttt}\label{CriterionOfExistence}
Let $G$ be a finite group. Then $G$ contains a Carter subgroup if and only if
$G$ satisfies~{\em\bfseries (E)}.
\end{ttt}

\begin{proof}
We prove the part ``only if'' first. Let $H$ be a minimal normal subgroup of
$G$. Then $H=T_1\times\ldots\times T_k$, where $T_1\simeq\ldots\simeq T_k\simeq
T$ are simple groups.

If $H$ is elementary Abelian (i.~e., $T$ is cyclic of prime order), then 
$\Aut(T)$ is solvable and contains a Carter subgroup. Assume that $T$ is a
non-Abelian simple group. Clearly $K$ is a Carter subgroup of $KH$.  By Lemma 
\ref{CarterInGroupOFInducedAutomorphismsWithoutCenter} we obtain that 
$\Aut_{KH}(T_i)$ contains a Carter subgroup for all~$i$. Induction by the order
of the group completes the proof of necessity.

Now we prove the ``if'' part. Again assume by contradiction that $G$ is a
counter example of minimal order, i.~e., that  $G$ does not contain a Carter
subgroup, but $G$ satisfies {\bfseries (E)}. Let $H$ be a minimal normal
subgroup of $G$. Then  $H=T_1\times\ldots\times T_k$, where $T_1\simeq
\ldots\simeq T_k\simeq T$, and $T$ is a finite simple group.

By definition,  $G/H$ satisfies {\bfseries (E)}, thus, by induction, there
exists a Carter subgroup  $\ov{K}$ of $\ov{G}=G/H$. Let $K$ be a complete
preimage of $\ov{K}$, then  $K$ satisfies {\bfseries (E)}. If $K\not=G$, then
by induction $K$ contains a Carter subgroup $K'$. Note that $K'$ is a Carter
subgroup of $G$. Indeed, assume that $x\in N_G(K')\setminus K'$. Since 
$K'H/H=\ov{K}$ is a Carter subgroup of $\ov{G}$, we obtain that 
$x\in K$. But $K'$ is a Carter subgroup of $K$, thus~${x\in K'}$. Hence 
$G=K$, i.~e., $G/H$ is nilpotent.

If $H$ is Abelian, then $G$ is solvable, therefore, $G$ contains a Carter
subgroup. So assume that $T$ is a non-Abelian finite simple group. We shall
show first that $C_G(H)$ is trivial. Assume that $C_G(H)=M$ is nontrivial.
Since $T$ is a non-Abelian simple group, it follows that $M\cap H=\{e\}$, so $M$
is nilpotent. By Lemma \ref{CarterSubgroupsInCompFactorsUnderHomomorphism} we
obtain that $G/M$ satisfies {\bfseries (E)}. By induction we obtain that $G/M$
contains a Carter subgroup $\ov{K}$. Let  $K'$ be a complete preimage of
$\ov{K}$ in $G$. Then $K'$ is solvable, therefore contains a Carter subgroup 
$K$. Like above we obtain that $K$ is a Carter subgroup of $G$, a
contradiction. Hence~${C_G(H)=\{e\}}$.

Since $H$ is a minimal normal subgroup of $G$, we obtain that $\Aut_G(T_1)\simeq
\Aut_G(T_2)\simeq\ldots\simeq \Aut_G(T_k)$. Thus there exists a monomorphism
$$\varphi:G\rightarrow 
\left(\Aut_G(T_1)\times\ldots\times\Aut_G(T_k)\right)\leftthreetimes
\Sym_k=G_1$$ and we identify $G$ with  $G^\varphi$. Denote by $K_i$ a Carter
subgroup of $\Aut_G(T_i)$ and by $A$ the subgroup 
$\Aut_G(T_1)\times\ldots\times\Aut_G(T_k)$. Since $G/H$ is nilpotent then
$K_iT_i=\Aut_G(T_i)$ and $G_1=\left(K_1T_1\times\ldots\times
K_kT_k\right)\leftthreetimes \Sym_k$. Let  $\pi_i:G\cap
A\rightarrow (G\cap A)/C_{(G\cap A)}(T_i)$ be canonical projections. Since 
$G/(G\cap A)$ is transitive, we obtain that 
$(G\cap A)^{\pi_i}=K_iT_i$.

Since  $\Aut_{G\cap A}(T_i)=K_iT_i$, then $G\cap A$ satisfies
{\bfseries (E)}. By induction it contains a Carter subgroup $M$. By Lemma 
\ref{CarterInGroupOFInducedAutomorphismsWithoutCenter} we obtain that 
$M^{\pi_i}$ is a Carter subgroup of $K_iT_i$, therefore we may assume that
$M^{\pi_i}=K_i$. In particular, if   $R=(K_1\cap T_1)\times\ldots\times(K_k\cap
T_k)$, then  $M\leq N_G(R)$. In view of Theorems 
\ref{ConjugacyCriterion} and \ref{FinalConjugacy}, Carter subgroups in each
finite group are conjugate. Since $(G\cap A)/H$ is nilpotent, we get that $G\cap
A=MH$, so $G=N_G(M)H$. More over $N_G(M)\cap A=M$, hence $N_G(M)$ is solvable.
Since $M$ normalizes $R$, and $M^{\pi_i}=K_i$, we obtain that $N_G(M)$
normalizes $R$, so $N_G(M)R$ is solvable. Therefore it contains a 
Carter subgroup  $K$. By Lemma~\ref{CarterInSubDirectProductOfSoluble}, 
$(K\cap A)^{\pi_i}$ is a Carter subgroup of $(N_G(M)R\cap A)^{\pi_i}$ ($R$
play the role of subgroup $T$ from Lemma~\ref{CarterInSubDirectProductOfSoluble}
in this case), so  $(K\cap A)^{\pi_i}=K_i$. Assume that 
$x\in N_G(K)\setminus K$. Since $G/H=N_G(M)H/H=KH/H$, it follows that 
$x\in H$. Therefore
$x^{\pi_i}\in (N_G(K)\cap A)^{\pi_i}\leq N_{T_i}((K\cap A)^{\pi_i})=K_i$.
Since $\bigcap_i \mathrm{Ker}(\pi_i)=\{e\}$, it follows that $x\in
R\leq N_G(M)R$. But $K$ is a Carter subgroup of $N_G(M)R$, hence $x\in
K$. This contradiction completes the proof.
\end{proof}

\subsection{Example}

In this subsection we shall construct an example, showing that we cannot
substitute condition {\bfseries (E)} by a weaker condition: for each
composition factor $S$ of $G$, $\Aut_G(S)$ contains a Carter subgroup. This
example also shows that an extension of a group containing a Carter subgroup, by
a group, containing a Carter subgroup, may fail to contain a Carter subgroup.

Consider $L=\P\SL_2(3^3)\leftthreetimes \la\varphi\ra$, where $\varphi$ is a
field automorphism of $\P\SL_2(3^3)$. Let $X=(L_1\times
L_2)\leftthreetimes \Sym_2$, where $L_1\simeq L_2\simeq L$ and if
$\sigma=(1,2)\in\Sym_2\setminus\{e\}$, $(x,y)\in L_1\times L_2$, then 
$\sigma(x,y)\sigma=(y,x)$ (permutation wreath product of $L$ and $\Sym_2$).
Denote by $H=\P\SL_2(3^3)\times \P\SL_2(3^3)$ a minimal normal subgroup of
$X$ and by $M=L_1\times L_2$. Let 
$G=(H\leftthreetimes \la(\varphi,\varphi^{-1})\ra)\leftthreetimes\Sym_2$ be a
subgroup of~$X$. Then the following statements hold:
\begin{itemize}
\item[1.] For every composition factor $S$ of $G$, $\Aut_G(S)$ contains a Carter
subgroup.
\item[2.] $G\cap M\unlhd G$ contains a Carter subgroup.
\item[3.] $G/(G\cap L)$ is nilpotent.
\item[4.] $G$ does not contain a Carter subgroup.
\end{itemize}

1. Clearly we need to check the statement only for non-Abelian composition
factors. Every non-Abelian composition factor $S$ of $G$ is isomorphic to
$\P\SL_2(3^3)$ and $\Aut_G(S)=L$. By Theorem  \ref{CarterSemilinear}, $L$
contains a Carter subgroup (that is equal to a Sylow $3$-sub\-gro\-up).

2. Since $(G\cap M)/H$ is nilpotent, and from the previous statement
we obtain that $G\cap M$ satisfies {\bfseries (E)}, so it contains a Carter
subgroup (it is easy to see that a Sylow $3$-sub\-gro\-up of $G\cap M$ is a
Carter subgroup of~${G\cap M}$).

3. Evident.

4. Assume that $K$ is a Carter subgroup of $G$. Then $KH/H$ is a Carter
subgroup of $G/H$. But $G/H$ is a non-Abelian group of order $6$, hence 
$G/H\simeq \Sym_3$ and $KH/H$ is a Sylow $2$-sub\-gro\-up of $G/H$. By Lemma
\ref{CarterSubgroupsInCompFactorsUnderHomomorphism} $\Aut_K(\P\SL_2(3^3))$
is a Carter subgroup of $\Aut_{KH}(\P\SL_2(3^3))=\P\SL_2(3^3)$.
But $\P\SL_2(3^3)$ does not contain Carter subgroups in view of
Theorem~\ref{CarterSemilinear}.

\subsection{Classification of Carter subgroups}

In view of condition  {\bfseries (E)} and Theorem \ref{CriterionOfExistence},
the description of Carter subgroups in finite groups is reduced to the
classification of Carter subgroups in almost simple groups $A$ with
$A/F^\ast(A)$ nilpotent. The classification of Carter subgroups in groups with
this condition is obtained in previous sections and we give it here for
convenient usage.

We prove first the following theorem, showing that if, for a subgroup 
$S$ of $\Aut(G)$, there exists a Carter subgroup, then it exists in every
larger group $S\leq A\leq \Aut(G)$ (here $G$ is a known simple group).

\begin{ttt}\label{AlmSimpleInherit}
Let $G$ be  a finite simple group and $G\leq A\leq \Aut(G)$ an almost simple
group with simple socle $G$. Assume that $A$ contains a subgroup
$S$ such that $G\leq S$ and $S$ contains a Carter subgroup.

\noindent Then $A$ contains a Carter subgroup.
\end{ttt}

\begin{proof}
Let $K$ be a Carter subgroup of $S$. Clearly we may assume that~${S=KG}$.

Assume that either $G\simeq \Alt_n$ for some $n\ge5$, or $G$ is sporadic. Since
by Lemma \ref{notl} each element of odd prime order of $G$ is conjugate 
to its inverse, and since $\vert \Aut(G):G\vert$ is a $2$-po\-w\-er, Lemmas
\ref{power} and \ref{InhBy2-ext} imply that if some group 
$G\leq S\leq \Aut(G)$ contains a Carter subgroup $K$, then $K$, is a Sylow 
$2$-sub\-gro\-up of $S$. Since $\vert A:S\vert$ is a $2$-po\-w\-er, the
statement of the theorem in this case follows from Lemma~\ref{InhBy2-ext}.

Assume that $G={}^3D_4(q)$. By \cite[Theorem~1.2(vi)]{TiepZalRealConjClasses}
each element of $G$ is conjugate to its inverse. If $q$ is odd, then Lemma
\ref{Syl2InCentrOfFieldAut} implies that $K$ is a Sylow $2$-sub\-gro\-up of
$S$. So by Lemmas \ref{InhBy2-ext} and \ref{Syl2InCentrOfFieldAut} it follows
that $A$ satisfies {\bfseries (ESyl2)}, i.~e., contains a Carter subgroup. If
$q=2^t$ is even, then by Theorems \ref{CarterTriality} and
\ref{CarterSemilinear} it follows that $S=\Aut(G)$ and we have nothing to prove.

Assume that $G$ is a group of Lie type, $G\not\simeq {}^3D_4(q)$ and, if
$G\simeq D_4(q)$, then  $S\leq A_1$, where $A_1\leq \Aut(D_4(q))$ is defined in
Theorem \ref{CarterTriality}. Then $S$ satisfies one of conditions (a)--(d) of
Theorem \ref{CarterSemilinear}. Consider all these cases separately.

Assume that $S$ satisfies (a). In this case we have
$\vert \Aut(G):S\vert\le 2$ and so, for each $A$ such that $S\leq A\leq
\Aut(G)$, either $A=S$, or $A=\Aut(G)$. In any case $A$ satisfies to statement
(a) of Theorem \ref{CarterSemilinear} and contains a Carter subgroup.

Assume that $S$ satisfies to statement (b). Then $\vert \Aut(G):S\vert=2$ and
either  $A=S$, or $A=\Aut(G)$. In the first case we have nothing to prove. In
the second case $\widehat{G}=\P GL_2(3^t)$ satisfies {\bfseries(ESyl2)}, hence
by Lemma \ref{ESyl2InhFieldGraph} group $A$ also satisfies {\bfseries(ESyl2)}
and by Lemma \ref{CritSyl2Carter} contains a Carter subgroup.

Assume that  $S$ satisfies statement (c) of Theorem \ref{CarterSemilinear}.
Then $S=\Aut(G)$ and we have nothing to prove.

Assume that $S$ satisfies condition (d) of Theorem \ref{CarterSemilinear}. By
Lemma \ref{Esyl2InnDiagExtension}, $S\cap\widehat{G}$ satisfies 
{\bfseries(ESyl2)}. By Lemma \ref{ESyl2InhFieldGraph} every subgroup 
$A$ of $\Aut{G}$, containing $S\cap \widehat{G}$ also satisfies 
{\bfseries(ESyl2)}, hence by Lemma \ref{CritSyl2Carter}, it contains a Carter
subgroup.

Now assume that $G=D_4(q)$ and $S$ satisfies conditions of Theorem 
\ref{CarterTriality}. Since graph automorphisms of orders $2$ and $3$ do not
commute, only one of them can be contained in a nilpotent subgroup. Thus we may
assume that only one of them is contained in  $A$. Then every subgroup $A$
containing $S$, either satisfies to Theorem \ref{CarterTriality}, or satisfies
to Theorem \ref{CarterSemilinear}, condition (a), if $q$ is even and condition
(d), if  $q$ is odd, i.~e., it contains a Carter subgroup.
\end{proof}

Note that from Theorem \ref{AlmSimpleInherit} and \cite{KondNormalizers} an
interesting corollary follows.

\begin{lem}\label{CarterInAlmSimpleAlwaysExists}
Let $S$ be a known finite simple group, $S\not\simeq J_1$ and
$G=\Aut(S)$. Then $G$ contains a Carter subgroup.
\end{lem}

\begin{proof}
By \cite[Theorems~2 and~3]{KondNormalizers}, if $S$ is not of Lie type and is
not equal to $J_1$, then group of its automorphisms $\Aut(S)$ satisfies
{\bfseries (ESyl2)} and, by Lemma \ref{CritSyl2Carter}, contains a Carter
subgroup. Now, if $S$ is of Lie type in even characteristic, then 
$\Aut(S)$ contains a Carter subgroup in view of Theorem
\ref{CarterSemilinear}(a). If $S$ is of Lie type in odd characteristic and 
$S\not\simeq {}^2G_2(3^{2n+1})$, then $\widehat{S}$ satisfies 
{\bfseries (ESyl2)}, so contains a Carter subgroup by Lemma
\ref{CritSyl2Carter}. By Theorem \ref{AlmSimpleInherit}, $\Aut(S)$
contains a Carter subgroup. Finally, if  $S\simeq {}^2G_2(3^{2n+1})$, then
$\Aut(S)$ contains a Carter subgroup in view of
Theorem~\ref{CarterSemilinear}(c).
\end{proof}

Tables given below are arranged in the following way. In the first column
is given a simple group $S$ such that Carter subgroups of $\Aut(S)$ are
classified. In the second column we give conditions for a subgroup  $A$ of its
group of automorphisms for  $A$ to contain a Carter subgroup. In the third
column we give the structure of a Carter subgroup $K$. In every subgroup of
$\Aut(S)$ lying between $S$ and $A$ Carter subgroups does non exist. By
$P_r(G)$\glossary{PrG@$P_r(G)$}  a Sylow $r$-sub\-gro\-up of $G$ is denoted. By
$\varphi$ we denote
 a field automorphism of a group of Lie type $S$, by $\tau$ we denote a graph
automorphism of a group of Lie type $S$ contained in $K$ (since graph
automorphisms of order  $2$ and $3$ of $D_4(q)$ does not commute, only one of
them can be in~$K$).  If $A$ does not contains a graph
automorphism, then we
suppose  $\tau=e$. By $\psi$ we denote a field automorphism of $S$ of maximal
order contained in $A$ (it is a power of $\varphi$, but
$\langle\psi\rangle$ can be different from $\langle\varphi\rangle$). By
$K(U_3(2))$ a Carter subgroup of order $2\cdot 3$ of $\widehat{{}^2A_2(2)}$
is denoted. If
$G$ is solvable, then by $K(G)$ we denote a Carter subgroup of
$G$. In Table \ref{CarterSubgroupsClassic} by $\zeta$ is denoted a
graph-fi\-eld automorphism of order~$2t$  of~$A_2(2^{2t})$.

\begin{longtable}{|c|c|c|}
\caption{Groups of automorphisms of alternating groups, containing Carter
subgroups \label{CarterSubgroupsAlt}}\\ \hline
Group $S$&Conditions on $A$&Structure of $K$\\ \hline
$\Alt_5$&$A=\Sym_5$&$K=P_2(\Sym_5)$\\ \hline
$\Alt_n$, $n\ge6$& none&$K=P_2(S)$\\ \hline
\end{longtable}

\begin{longtable}{|c|c|c|}
\caption{Groups of automorphisms of sporadic groups, containing Carter
subgroups \label{CarterSubgroupsCporadic}}\\ \hline
Group $S$&Conditions on $A$&Structure of $K$\\ \hline
$J_2,J_3,Suz, HN$&$A=\Aut(S)$&$K=P_2(A)$\\ \hline
$\not\simeq J_1,J_2,J_3,Suz,HN$&none&$K=P_2(A)$\\ \hline
\end{longtable}

\begin{longtable}{|c|c|c|}
\caption{Groups of automorphisms of exceptional groups of Lie type,
containing Carter subgroups \label{CarterSubgroupsExcept}}\\ \hline
Group $S$&Conditions on $A$&Structure of $K$\\ \hline
${}^2B_2(2^{2n+1})$, $n\ge1$&$A=\Aut(S)$&$K=\langle\varphi\rangle\times
P_2({}^2B_2(2))$\\ \hline
$({}^2F_4(2))'$&none&$K=P_2(A)$\\ \hline
${}^2F_4(2^{2n+1})$, $n\ge1$&$A=\Aut(S)$&$K=\langle\varphi\rangle\times
P_2({}^2F_4(2))$\\ \hline
${}^2G_3(3^{2n+1})$&$A=\Aut(G)$&$\langle \varphi\rangle\rightthreetimes
(2\times P)$,\\
&&where $\vert P\vert=3^{\vert\varphi\vert_3}$\\ \hline remaining, $q$
odd&none&$K=P_2(A)\times K(O(N_A(P_2(A))))$\\ \hline
remaining, $q=2^t$&$\varphi g\in
A$, $g\in\widehat{S}$&$\langle\tau,\varphi\rangle\rightthreetimes
P_2(S_{\varphi_{2'}})$\\ \hline
\end{longtable}

{\footnotesize
\begin{longtable}{|c|c|c|}
\caption{Groups of automorphisms of classical groups, containing Carter
subgroups\label{CarterSubgroupsClassic}}\\ \hline
Group $S$&Conditions on $A$&Structure of $K$\\ \hline
$A_1(q)$, $q\equiv\pm1\pmod8$&none& $K=N_A(P_2(S))$\\ \hline
$A_1(q)$, $q\equiv\pm3\pmod8$&$\widehat{S}\leq A$&$K=N_A(P_2(\widehat{S}))$\\
\hline $A_n(2^t)$, $t\ge2$, if $n=1$&$\varphi g\in A$,
$g\in\widehat{S}$& $K=\langle \varphi,\tau\rangle\rightthreetimes
S_{\varphi_{2'}}$\\ \hline $A_2(2^{2t})$, $3\nmid t$&$\la S,\zeta
g\ra\leq A\leq S\leftthreetimes\la\zeta\ra$,&$K=\langle\zeta
g\rangle\times K(\P\GU_3(2))$\\&$C_{A\cap
\widehat{S}}(\varphi_{2'})\simeq\P\GU_3(2)$&\\ \hline $A_n(q)$, $q$
odd, $n\ge2$&none&$K=P_2(A)\times K(O(N_A(P_2(A))))$\\ \hline
${}^2A_2(2^{t})$, $t$ odd, $3\nmid t$ & $\langle S, \varphi_{2'}
g\rangle\leq A\leq
\widehat{S}\leftthreetimes\langle\varphi_{2'}\rangle $&\\
&$C_{A\cap
\widehat{S}}(\varphi_{2'})\simeq\P\GU_3(2)$&$K=\langle\varphi_{2'}\rangle\times
K(\P\GU_3(2))$\\ &$C_{A\cap
\widehat{S}}(\varphi_{2'})\simeq\P\SU_3(2)$&$K=\langle\varphi_{2'}\rangle\times
P_2(\P\SU_3(2))$\\ \hline
${}^2A_2(2^{t})$&$A=\Aut(S)$&$K=\langle\varphi
 \rangle\rightthreetimes P_2(S_{\varphi_{2'}})$\\ \hline
${}^2A_n(q)$, $q$ odd&none&$K=P_2(A)\times
K(O(N_A(P_2(A))))$\\ \hline ${}^2A_n(2^{t})$,
$n\ge3$&$A=\Aut(S)$&$K=\langle\varphi\rangle\rightthreetimes
P_2(S_{\varphi_{2'}})$\\ \hline $B_2(q)$,
$q\equiv\pm1\pmod8$&none&$K=P_2(A)\times K(O(N_A(P_2(A))))$\\
\hline $B_2(2^t)$, $t\ge 2$&$\varphi\in A$&$K=\langle
\varphi,\tau\rangle\rightthreetimes P_2((S_\tau)_\varphi)$\\ \hline
$B_2(q)$, $q\equiv\pm3\pmod8$&$\widehat{S}\leq A$&$K=P_2(A)\times
K(O(N_A(P_2(A))))$\\
\hline $B_n(q)$, $q$ odd, $n\ge3$&none&$K=P_2(A)\times
K(O(N_A(P_2(A))))$\\
\hline $C_n(q)$, $q\equiv\pm1\pmod8$& none&$K=P_2(A)\times
K(O(N_A(P_2(A))))$\\
\hline
$C_n(q)$, $q\equiv\pm3\pmod8$& $\widehat{S}\leq A$&$K=P_2(A)\times
K(O(N_A(P_2(A))))$\\
\hline $C_n(2^t)$,
$n\ge3$&$A=\Aut(S)$&$K=\langle\varphi\rangle\times
P_2(S_{\varphi_{2'}})$\\ \hline
$D_4(q)$, $q$ odd&none&if $\vert\tau\vert\le2$, then\\
&&$K=P_2(A)\times K(O(N_A(P_2(A))))$;\\
&&if $\vert\tau\vert=3$, then\\ &&
$K=\langle\tau,\psi\rangle\rightthreetimes P_2(S_\tau)$\\ \hline
$D_4(2^t)$&$\varphi\in A$&if $\vert\tau\vert\le2$, then\\
&&$K=\langle \tau,\varphi\rangle\rightthreetimes P_2(S_{\varphi_{2'}})$;\\
&&if $\vert\tau\vert=3$, then\\ &&$K=\langle
\tau,\varphi\rangle\rightthreetimes P_2((S_\tau)_{\varphi_{2'}})$\\
\hline $D_n(q)$, $q$ odd, $n\ge5$&none&$K=P_2(A)\times
K(O(N_A(P_2(A))))$\\ \hline $D_n(2^t)$, $n\ge5$&$\varphi\in
A$&$K=\langle\tau,\varphi\rangle\rightthreetimes
P_2(S_{\varphi_{2'}})$\\ \hline ${}^2D_n(q)$, $q$
odd&none&$K=P_2(A)\times K(O(N_A(P_2(A))))$\\ \hline
${}^2D_n(2^{t})$&$A=\Aut(S)$&$K=\langle\varphi\rangle\rightthreetimes
P_2(S_{\varphi_{2'}})$\\ \hline
\end{longtable}}




\addcontentsline{toc}{section}{List of Tables}
\listoftables

\end{document}